\documentclass[11pt,a4paper,notitlepage,oneside]{amsart}
\usepackage{xcolor}
\usepackage{etoolbox}
\usepackage{microtype}
\usepackage{amsmath}
\usepackage{amsthm}
\usepackage{braket}
\usepackage[psamsfonts]{amssymb}
\usepackage[T1]{fontenc}
\usepackage[english]{babel}
\usepackage[utf8]{inputenc}
\usepackage{enumitem}
\usepackage{mathtools}
\usepackage{caption,subfig}
\usepackage{rotating}
\usepackage[makeroom]{cancel}
\usepackage[normalem]{ulem}

\usepackage{scalerel}
\usepackage{relsize}

\usepackage[autostyle,italian=guillemets]{csquotes}
\usepackage[style=alphabetic,firstinits=true,maxbibnames=99,backend=bibtex]{biblatex}
\bibliography{bibliography}
\allowdisplaybreaks 

\usepackage{pgf,tikz,pgfplots}

\usepackage{mathrsfs}
\usetikzlibrary{arrows}

\usepackage{fancyhdr}
\usepackage{tikz-cd}
\usepackage{amscd}
\usepackage[a4paper, bottom=3cm, left=3cm,right=3cm]{geometry}
\usepackage{hyperref}
\hypersetup{hidelinks}

\usepackage{subfiles}
\usepackage{my_shortcuts}

\author{Riccardo Zuffetti}
\title{Cones of special cycles of codimension~2 on orthogonal Shimura varieties}
\address{
Institut f\"ur Mathematik, Goethe--Universit\"at Frankfurt,
Robert-Mayer-Str. 6--8,
60325 Frankfurt am Main, Germany\\
}
\email{zuffetti@math.uni-frankfurt.de, riccardo.zuffetti@gmail.com}

\newcommand{\myblack}{black}

\begin{document}
	\maketitle
	\begin{abstract}
	Let $X$ be an orthogonal Shimura variety associated to a unimodular lattice.
	We investigate the polyhedrality of the cone $\mathcal{C}_X$ of special cycles of codimension 2 on~$X$.
	We show that the rays generated by such cycles accumulate towards infinitely many rays, the latter generating a non-trivial cone.
	We also prove that \textcolor{\myblack}{such an} accumulation cone is polyhedral.
	The proof relies on analogous properties satisfied by the cones of Fourier coefficients of Siegel modular forms.
	We show that the accumulation rays of~$\mathcal{C}_X$ are generated by combinations of Heegner divisors intersected with the Hodge class of~$X$.
	\textcolor{\myblack}{As a result of the classification of the accumulation rays,} we implement an algorithm in SageMath to certify the polyhedrality of~$\mathcal{C}_X$ in some cases.
	\end{abstract}
	\tableofcontents
	\section{Introduction}\label{sec;intropolyhed}

	The cones of divisors and the cones of curves on complex (quasi-)projective varieties have been intensely studied.
	In recent years the cones of cycles in higher codimension have attracted increasing interest.
	Up to now, few examples of cones of cycles have been computed;
	see e.g.\ the cases of projective bundles over curves~\cite{fulgconeprojbund}, moduli spaces of curves~\cite{chencosk}~\cite{effcone2cymodc}, spherical varieties~\cite{lisphericalvar}, point blow-ups of projective spaces~\cite{coleotch}, symmetric products of curves~\cite{effsymprod} and moduli spaces of principally polarized abelian threefolds~\cite{gruhu}.
	
	In this paper we illustrate properties of certain cones of \emph{codimension~2 cycles}.
	We focus on orthogonal Shimura varieties, studying the geometric properties of the cones generated by the codimension~$2$ \emph{special cycles} via the arithmetic properties of the Fourier coefficients of genus~$2$ Siegel modular forms.
	
	The case of codimension~$1$ has already been treated in~\cite{brmo}, in which the cones of \emph{special divisors} (also known as \emph{Heegner divisors}) is proven to be rational and polyhedral, whenever these varieties arise from lattices that split off a hyperbolic plane.
	In analogy with~\cite{brmo}, it is tempting to prove the following conjecture.
	Part of this paper aims to provide supporting results to it, together with a reformulation in terms of Jacobi forms and a proof of the very first cases in low dimension.
	
	\begin{conj}\label{conj;ourchow}
	Let~$X$ be an orthogonal Shimura variety arising from an even unimodular lattice of signature~$(b,2)$, where~$b>2$.
	The cone of special cycles of codimension~$2$ on~$X$ is polyhedral.
	\end{conj}
	To study the properties of cones of special cycles, it is useful to find all rays arising as limits of rays generated by sequences of pairwise different special cycles.
	For instance, the properties of the cone of special divisors studied in~\cite{brmo} are deduced by showing that the rays generated by such divisors accumulate towards a \emph{unique ray}, and that the latter lies in the interior of the cone.
	
	In this paper we show that the situation of cones generated by special cycles of codimension~$2$ is much more interesting, although more complicated, and in fact that the rays generated by such cycles \emph{accumulate towards a non-trivial cone}.
	The latter, which in our terminology will be called the \emph{accumulation cone}, is the main subject of investigation of this article.
	We illustrate the geometric properties of \textcolor{\myblack}{this} cone, and classify all rays arising as limits as explained above.
	Eventually, thanks to such classification, we translate Conjecture~\ref{conj;ourchow} in terms of Jacobi forms.
	
	To illustrate the main achievements of the article, we need to introduce some notation.	
	Let~$V$ be a finite dimensional vector space over~$\QQ$, and let~$V_\RR$ be the vector space~$V\otimes\RR$ endowed with the Euclidean topology.
	As recalled above, to study the properties of a (convex) cone~$\cone$ generated by some subset~$\gencone\subseteq V$, in short~$\cone=\langle\gencone\rangle_{\QQ_{\ge0}}$, it is useful to find all rays in~$V_\RR$ arising as limits of rays generated by sequences of elements in~$\gencone$.
	This motivates the following definition.
	
	\begin{defi}
	A ray $r$ of~$V_\RR$ is said to be an \emph{accumulation ray of~$\cone$ with respect to the set of generators~$\gencone$} if there exists a sequence of pairwise different generators~$\seq{g_j}{j\in\NN}$ in~$\gencone$ such that
	\bes
	\RR_{\ge0}\cdot g_j \longrightarrow r,\qquad \text{when~$j\longrightarrow\infty$},
	\ees	
	where we denote by~$\RR_{\ge0}\cdot g_j$ the ray generated by~$g_j$ in~$V_\RR$.
	The \emph{accumulation cone of $\cone$ with respect to~$\gencone$} is defined as the cone in~$V_\RR$ generated by~$0$ and the accumulation rays of~$\cone$ with respect to~$\gencone$.
	\end{defi}
	
	By what we recalled above on~\cite{brmo}, the accumulation cone of the cone generated by special divisors is of dimension~$1$, in fact simply a ray.
	The cones generated by special cycles of codimension~$2$ have a much more interesting geometry, as we are going to explain.
	
	Let $X$ be a Shimura variety of orthogonal type, arising from an even \emph{unimodular} lattice of signature $(b,2)$, with $b> 2$.
	The special cycles of~$X$ are suitable sums of orthogonal Shimura subvarieties of $X$, and are parametrized by half-integral positive semi-definite matrices of order~$2$.
	We denote by $\halfint_2$ the set of these matrices, and by $\halfint^+_2$ the subset of the ones whose determinant is positive. 
	If~$T\in\halfint^+_2$, we denote by $Z(T)$ the special cycle associated to $T$, and by~$\{Z(T)\}$ its rational class in the Chow group~$\CH^2(X)$.
	If $T$ is singular, it is still possible to define a special cycle in~$\CH^2(X)$ by intersecting with (the dual of) the rational class of the Hodge bundle~$\omega$ of~$X$.
	We refer to Section~\ref{sec;conesspeccycles} for further details.
	\begin{defi}\label{def;introconecynewbis}
	The \emph{cone of special cycles (of codimension~$2$) on $X$} is the cone defined in~${\CH^2(X)\otimes\QQ}$ as
	\bes
	\cone_{X}=\langle \{Z(T)\}\,:\, T\in\halfint_2^+\rangle_{\QQ_{\ge 0}}.
	\ees
	The \emph{cone of rank one special cycles (of codimension~$2$) on $X$} is
	\bes
	{\cone'}_{\!\! X}=\langle \{Z(T)\}\cdot\{\omega^*\}\,:\, \text{$\T\in\halfint_2$ and $\rank(T)=1$} \rangle_{\QQ_{\ge 0}}.
	\ees
	\end{defi}
	Whenever we refer to the \emph{accumulation cones} or \emph{accumulation rays} of~$\cone_{X}$ and~${\cone'}_{\!\! X}$, we implicitly consider them with respect to the set of generators of~$\cone_{X}$ and~${\cone'}_{\!\! X}$ used in Definition~\ref{def;introconecynewbis}.
	All these cones are of finite dimension.
	
	We briefly recall some properties of cones, referring to Section~\ref{sec;backgrcones} for a more detailed explanation.
	Let~$\cone$ be a cone in a finite dimensional vector space~$V$ over~$\QQ$.
	We say that~$\cone$ is \emph{pointed} if it contains no lines.
	The~\emph{$\RR$-closure}~$\overline{\cone}$ is the topological closure of~$\cone$ in~$V_\RR$.
	We say that~$\cone$ is \emph{rational} if~$\overline{\cone}$ may be generated over~$\RR$ by a subset of the rational space~$V$.
	The cone~$\cone$ is \emph{polyhedral} if~$\cone=\langle\gencone\rangle_{\QQ_{\ge0}}$, for some finite set of generators~$\gencone$.
	
	It is a direct consequence of~\cite{brmo} that the cone of rank one special cycles~${\cone'}_{\!\! X}$ has a unique accumulation ray, and that such ray is internal in the $\RR$-closure of~${\cone'}_{\!\! X}$.
	The cone~$\cone_{X}$ has a more interesting, although intricate, geometry, as shown by the following result.
	We denote by~$M^k_1$ the space of weight~$k$ elliptic modular forms.
	\begin{thm}\label{thm;mainintrobis}
	Let $X$ be an orthogonal Shimura variety associated to an even unimodular lattice of signature~$(b,2)$, with~$b>2$.
	\begin{enumerate}[label=(\roman{*})]
	\item The cone~$\cone_{X}$ has \emph{infinitely many} accumulation rays if and only if~$b\ge 34$.
	They are all generated by~$\QQ$-cycles, in particular~$\cone_{X}$ is rational.
	\item The accumulation cone of the cone of special cycles~$\cone_{X}$ is pointed, rational, polyhedral, and of dimension~$\dim {{M_1}^{\!1+b/2}}$.
	\end{enumerate}
	\end{thm}
	Along the proof of Theorem~\ref{thm;mainintrobis} we will deduce the following relation between~$\cone_{X}$ and~${\cone'}_{\!\! X}$.
	\begin{prop}
	The cones~$\cone_{X}$ and~${\cone'}_{\!\! X}$ intersect only at the origin.
	Moreover, the accumulation cone of~$\cone_{X}$ contains~$-{\cone'}_{\!\! X}$.
	\end{prop}
	In this article we also classify all accumulation rays of~$\cone_{X}$.
	The following result provides the ones that may be extremal in the accumulation cone of~$\cone_{X}$; see Section~\ref{sec;accrays} and Section~\ref{sec;acrayschow} for a detailed description of all remaining accumulation rays.
	We denote by~$\mu$ the Möbius function.
	\begin{thm}\label{thm;introacconeshimura}
	The accumulation cone of~$\cone_{X}$ is generated by the rays
	\be\label{eqv2;introacconeshimura}
	\textcolor{\myblack}{\RR_{\ge0}\cdot\bigg(\sum_{t^2|m}\mu(t)
	\{H_{m/t^2}\}\cdot\{\omega\}
	\bigg),}
	\ee
	indexed over~$m\in\ZZ_{>0}$, where~$H_n$ denotes the~$n$-th Heegner divisor of~$X$.
	\end{thm}
	\textcolor{\myblack}{If the arithmetic group giving rise to the orthogonal Shimura variety~$X$ is large enough, then the ray~\eqref{eqv2;introacconeshimura} is generated by the intersection of~$\{\omega\}$ with an irreducible component of the Heegner divisor~$\{H_m\}$.
	This is explained in Section~\ref{sec;acrayschow} by means of the so-called \emph{primitive} Heegner divisors.}
	
	The key result to deduce Theorem~\ref{thm;mainintrobis} is Kudla's modularity conjecture, recently proved by Bruinier and Raum~\cite{brra;modconj}, which enables us to deduce geometric properties of $\cone_X$ via arithmetic properties of the Fourier coefficients of genus~$2$ Siegel modular forms, as we briefly recall.

	Let $k$ be a positive even integer and let $M^k_2(\QQ)$ be the vector space over~$\QQ$ of weight~$k$ and genus~$2$ Siegel modular forms with \emph{rational} Fourier coefficients. For every~$F$ in~$M^k_2(\QQ)$, we denote the Fourier expansion of $F$ by
	\bes
	F(Z)=\sum_{T\in\halfint_2}c_T(F)e^{2\pi i\trace(TZ)},
	\ees
	where $Z$ lies in the Siegel upper-half space $\HH_2$, and $c_T(F)$ is the rational Fourier coefficient of~$F$ associated to the matrix $T\in\halfint_2$. The dual space~$M^k_2(\QQ)^*$ is generated by the \emph{coefficient extraction functionals}~$c_T$, defined for every~$T\in\halfint_2$ as
	\bas
	c_T\colon M^k_2(\QQ)\longrightarrow\QQ,\quad F\longmapsto c_T(F).
	\eas
	The main result of \cite{brra;modconj} implies that the linear map
	\bas
	\psi_X\colon M^{1+b/2}_2(\QQ)^*\longrightarrow \CH^2(X)\otimes\QQ,\quad c_T\longmapsto\{Z(T)\}\cdot\{\omega^*\}^{2-\rank(T)},
	\eas
	is well-defined.
	Note that $1+b/2$ is an even integer, in fact~$1+b/2\equiv 2$ mod~$4$.
	This follows from the well-known classification of even indefinite unimodular lattices.
	\begin{defi}
	The \emph{modular cone of weight $k$} is the cone in $M^k_2(\QQ)^*$ defined as
	\bes
	\cone_k=\langle c_T\,:\,T\in\halfint^+_2 \rangle_{\QQ_{\ge0}}.
	\ees
	\end{defi}
	The key idea of this paper is to deduce the properties of the cone of special cycles~$\cone_X$ appearing in Theorem~\ref{thm;mainintrobis} proving analogous properties of the associated modular cone~$\cone_{1+b/2}$.
	In fact, such properties are preserved \textcolor{\myblack}{under} the linear map~$\psi_X$, as we prove in Section~\ref{sec:concoefextrfun}.
	\begin{thm}\label{thm;intromodcone}
	Suppose that $k\equiv 2$ mod $4$ and $k>4$.
	\begin{enumerate}[label=(\roman{*})]
	\item The modular cone~$\cone_k$ is pointed and full-dimensional in~$M^k_2(\QQ)^*$.
	Its accumulation rays are all generated over~$\QQ$, and are infinitely many if and only if~$k\ge 18$.
	\item The accumulation cone of~$\cone_k$ is pointed, rational, polyhedral, and of the same dimension as~$M^k_1(\QQ)$.
	\end{enumerate}
	\end{thm}
	If the map~$\psi_X$ is injective, then also the pointedness and the dimension of~$\cone_{1+b/2}$ are preserved, implying that~$\cone_X$ has such properties too;
	see Remark~\ref{rem;nonpointandinjpsigamma} for more information on the injectivity of~$\psi_X$.

	In Section~\ref{sec;accrays} we provide a complete classification of all accumulation rays of the modular cone~$\cone_k$, for every integer~${k>4}$ such that~${k\equiv 2}$ mod~$4$, from which we deduce Theorem~\ref{thm;introacconeshimura}.
	
	The results on the accumulation rays of the modular cone $\cone_k$ are deduced via estimates of the growth of the Fourier coefficients of genus~$2$ Siegel modular forms, and via the values assumed by certain ratios of Fourier coefficients of the weight $k$ Siegel Eisenstein series.
	The main difficulty arising in genus~$2$ is the presence of the so-called \emph{Klingen Eisenstein series}, which do not appear if only elliptic modular forms are considered, as in \cite{brmo}.
	The main resource we use to treat this issue is the recent paper \cite{boda;petnorm}, where the growth of the coefficients of the Klingen Eisenstein series is clarified; we refer to Sections~\ref{sec;elmodandjac} and~\ref{sec;siegmodforms} for the needed background.
	
	In Section~\ref{sec;accumulcone} and Section~\ref{sec;propC} we prove Theorem~\ref{thm;intromodcone}.
	We furthermore provide a sufficient condition to the polyhedrality of~$\cone_k$ (hence also of~$\cone_X$); see Theorem \ref{thm;polyhedwithhyp}.
	We implemented an algorithm~\cite{zufprog} in SageMath to check that such condition is surely satisfied if~$k\le 38$.
	In fact, such program provides empirical evidences to the polyhedrality of~$\cone_k$ also for~$k>38$.
	This leads us to the following conjecture.
	\begin{conj}\label{conj;our}
	Suppose that $k\equiv 2$ mod $4$ and $k>4$.
	The cone $\cone_k$ is polyhedral.
	\end{conj}
	\begin{prop}
	If~$k\le 38$, then Conjecture~\ref{conj;our} is satisfied.
	In particular, the cone of special cycles~$\cone_X$ of an orthogonal Shimura variety~$X$ arising from an even unimodular lattice of signature~$(b,2)$, where~$2<b\le74$, is polyhedral.
	\end{prop}
	We conclude Section~\ref{sec;propC} reducing the problem of the polyhedrality of~$\cone_k$ to ``how a sequence of rays~$\seq{\RR_{\ge0}\cdot c_{T_j}}{j\in\NN}$ converges towards the accumulation cone of $\cone_k$'', with a translation of Conjecture~\ref{conj;our} into a conjecture on Fourier coefficients of Jacobi cusp forms.
	
	\subsection*{Acknowledgments}
	We are grateful to Martin Möller for suggesting this project and for his encouragement.
	His office door was always open for discussions.
	We would like to thank also Siegfried Böcherer, Jan Bruinier, Jolanta Marzec, Martin Raum and Brandon Williams for useful conversations and/or emails they shared with us\textcolor{\myblack}{, as well as the anonymous referee for a thorough report and helpful suggestions}.
	This work is a result of our PhD~\cite{zufthesis}, which was founded by the LOEWE research unit ``Uniformized Structures in Arithmetic and Geometry'', and by the Collaborative Research Centre TRR~326 ``Geometry and Arithmetic of Uniformized Structures'', project number~444845124.
	
	\section{Elliptic and Jacobi modular forms}\label{sec;elmodandjac}
	
	To fix the notation, in this section we recall the definitions of elliptic and Jacobi modular forms.
	Eventually, we illustrate some properties about positive linear combinations of coefficients extraction functionals associated to these forms.
	Such properties will be essential in Section \ref{sec;accumulcone} to prove that certain accumulation rays of $\cone_k$ lies in the interior of its accumulation cone.
	
	For the purposes of this paper, we do not need to consider congruence subgroups, hence all modular forms here treated are with respect to the full modular groups. Introductory books are e.g.~\cite{1-2-3} and \cite{eiza;theoryjacobi}.\\
	
	We begin with elliptic modular forms. The modular group $\SL_2(\ZZ)$ acts on the upper-half plane $\HH$ via the \emph{Möbius transformation} as
	\bes
	\gamma=\left(\begin{smallmatrix}
	a & b\\
	c & d
	\end{smallmatrix}\right):\HH\longrightarrow\HH,\quad \tau\longmapsto\gamma\cdot\tau=\frac{a\tau +b}{c\tau +d},
	\ees
	where $\gamma\in \SL_2(\ZZ)$.
	Let $k>2$ be an even integer and let $f\colon\HH\to\CC$ be a holomorphic function on the upper-half plane.
	We say that $f$ is an \emph{elliptic} (or \emph{genus~$1$}) \emph{modular form of weight $k$} if $f$ satisfies $f(\gamma\cdot\tau)=(c\tau+d)^kf(\tau)$ for all $\tau\in\HH$ and all $\gamma=\left(\begin{smallmatrix}
	a & b\\
	c & d
	\end{smallmatrix}\right)\in \SL_2(\ZZ)$, and if it admits a Fourier expansion of the form
	\bes
	f(\tau)=\sum_{n=0}^\infty c_n(f)q^n,\quad\text{where $q=e^{2\pi i\tau}$.}
	\ees
	The complex number~$c_n(f)$ is the \emph{$n$-th Fourier coefficient} of $f$.
	We denote the finite-dimensional complex vector space of weight $k$ elliptic modular forms by $M^k_1$.
	We put the subscript $1$ to recall that these are modular forms of genus~$1$, avoiding confusion with the Siegel modular forms we are going to define in Section \ref{sec;siegmodforms}. 
	The first examples of such functions are the (normalized) \emph{Eisenstein series}
	\be\label{eq;normeisserelliptic}
	E^k_1(\tau)=1+\frac{2}{\zeta(1-k)}\sum_{n=1}^\infty\sigma_{k-1}(n)q^n,
	\ee
	where $\zeta(s)$ is the Riemann zeta function, and $\sigma_{k-1}(n)$ is the sum of the $(k-1)$-powers of the positive divisors of $n$.
	
	An \emph{elliptic cusp form of weight $k$} is a modular form $f\in M^k_1$ such that its first Fourier coefficient is trivial, namely $c_0(f)=0$. We denote by $S^k_1$ the subspace of cusp forms of weight~$k$.
	It is well-known that the space of elliptic modular forms decomposes as~${M^k_1=\langle E^k_1\rangle_\CC\oplus S^k_1}$.
	
	We denote by $M^k_1(\QQ)$ (resp.\ $S^k_1(\QQ)$) the space of elliptic modular forms (resp.\ cusp forms) with \emph{rational} Fourier coefficients. Since $S^k_1$ admits a basis of cusp forms with rational coefficients, it turns out that the dimension of	$M^k_1(\QQ)=\langle E^k_1\rangle_\QQ\oplus S^k_1(\QQ)$ over $\QQ$ is equal to the complex dimension of $M^k_1$.
	The dual space $M^k_1(\QQ)^*$ is generated by the \emph{coefficient extraction functionals} $c_n$, defined as
	\bes
	c_n\colon M^k_1(\QQ)\longrightarrow\QQ,\quad f\longmapsto c_n(f),
	\ees
	for every $n\ge0$. In \cite{brmo}, the authors proved that whenever $k\equiv 2$ mod $4$, the cone generated by the functionals $c_n$ with $n\ge1$ is rational polyhedral in $M^k_1(\QQ)$.
	A key result used in the cited paper is~\cite[Proposition 3.3]{brmo}, here stated in our setting.
	\begin{lemma}\label{lemma;propbmc}
	Suppose that $k\equiv 2$ mod $4$. There exist a positive integer $A$ and positive rational numbers $\eta_j$, with $j=1,\dots,A$, such that
	\bes
	\sum_{j=1}^A\eta_j \cdot c_{j}|_{S^k_1(\QQ)}=0.
	\ees
	Furthermore, the constant $A$ can be chosen arbitrarily large such that the restrictions $c_{j}|_{S^k_1(\QQ)}$ generate~$S^k_1(\QQ)^*$.
	\end{lemma}\noindent
	For the purposes of this paper, we need a slight generalization of Lemma~\ref{lemma;propbmc} to Jacobi forms, as we are going to illustrate.\\
	
	Jacobi forms play an important role in the study of the Fourier coefficients of Siegel modular form.
	As we will recall in the next sections, the Fourier series of a Siegel modular form can be rewritten in terms of Jacobi forms.
	This \emph{arithmetic} property will be translated into a \emph{geometric} property of the cone $\cone_k$ we defined in the introduction. Namely, the convergence of certain sequences of rays in $\cone_k$ will be deduced from results on the growth of Fourier coefficients of Jacobi forms. This is one of the goals of Section \ref{sec;accrays}.
	
	Let $k>2$ be an even integer, and let~$m\in\ZZ_{\ge0}$. A holomorphic function $\phi:\HH\times \CC\to\CC$ is said to be a \emph{Jacobi form of weight $k$ and index $m$} if
	\bas
	\phi\left(\frac{a\tau+b}{c\tau+d},\frac{z}{c\tau+d}\right)&=(c\tau+d)^ke^{\frac{2\pi imcz^2}{c\tau+d}}\phi(\tau,z),\quad \text{for every $\left(\begin{smallmatrix}
	a & b\\ c & d
	\end{smallmatrix}\right)\in \SL_2(\ZZ)$},\\
	\phi(\tau,z+\lambda\tau+\mu)&=e^{-2\pi im(\lambda^2\tau+2\lambda z)}\phi(\tau,z),\quad \text{for every $(\lambda,\mu)\in\ZZ^2$},
	\eas
	and if $\phi$ admits a Fourier expansion of the form
	\be\label{eq;fourdecphi}
	\phi(\tau,z)=\sum_{n=0}^\infty\sum_{\substack{r\in\ZZ\\ 4nm-r^2\ge0}}c_{n,r}(\phi) q^n\zeta^r,\quad\text{where $q=e^{2\pi in\tau}$ and $\zeta=e^{2\pi irz}$}.
	\ee
	The complex numbers $c_{n,r}(\phi)$ are the \emph{Fourier coefficients} of $\phi$.
	We denote by $J_{k,m}$ the finite-dimensional complex vector space of such functions.
	If in the Fourier expansion~\eqref{eq;fourdecphi} the coefficients~$c_{(n,r)}(\phi)$ such that $4nm=r^2$ are zero, then~$\phi$ is said to be a \emph{Jacobi cusp form}.
	We denote the space of these forms by $J_{k,m}^{\text{cusp}}$.
	
	First explicit examples of Jacobi forms are the \emph{Jacobi Eisenstein series}.
	We avoid to define them explicitly in this paper, we refer instead to~\cite[Section 2]{eiza;theoryjacobi} for a detailed introduction.
	The subspace generated by the Jacobi Eisenstein series is denoted by $J_{k,m}^{\text{Eis}}$.
	By \cite[Theorem 2.4]{eiza;theoryjacobi}, the space of Jacobi forms of even weight $k>2$ decomposes into
	\be\label{eq;strthmJac}
	J_{k,m}=J_{k,m}^{\text{Eis}}\oplus J_{k,m}^{\text{cusp}}.
	\ee
	
	In analogy with the case of elliptic modular forms, the spaces $J_{k,m}^{\text{Eis}}$ and $J_{k,m}^{\text{cusp}}$ admit a basis of Jacobi forms with rational Fourier coefficients.
	We denote the associated spaces of Jacobi forms with rational coefficients by $J_{k,m}^{\text{Eis}}(\QQ)$ and $J_{k,m}^{\text{cusp}}(\QQ)$, respectively. An analogous decomposition as~\eqref{eq;strthmJac} holds also over $\QQ$.
	
	The dual space $J_{k,m}(\QQ)^*$ is generated by the \emph{Jacobi coefficient extraction functionals}~$c_{n,r}$, defined as
	\bes
	c_{n,r}\colon J_{k,m}(\QQ)\longrightarrow\QQ,\quad \phi\longmapsto c_{n,r}(\phi),
	\ees
	for every $n\ge0$ and $r\in\ZZ$ such that $4nm-r^2\ge0$.
	
	The slight generalization of Lemma \ref{lemma;propbmc} previously announced is the following.
	\begin{lemma}\label{lem:relcoefextjac}
	Suppose that $k\equiv 2$ mod $4$. For every positive integer $m$ there exist a positive integer $A$ and positive rational numbers $\mu_{n,r}$ such that 
	 \be\label{eq:relcoefextjac}
	 \sum_{1\le n\le A}\sum_{\substack{r\in\ZZ\\ 4nm-r^2>0}}\mu_{n,r}c_{n,r}|_{J_{k,m}^{\text{\rm cusp}}(\QQ)}=0.
	 \ee
	 Furthermore, the constant~$A$ can be chosen arbitrarily large such that the restricted functionals~$c_{n,r}|_{J_{k,m}^{\text{\rm cusp}}(\QQ)}$ generate~$J_{k,m}^{\text{\rm cusp}}(\QQ)^*$.
	 \end{lemma}
	 \begin{proof}
	 If $\phi\in J_{k,m}(\QQ)$, then the map on $\HH$ defined as $\phi(\tau,0)$ lies in $M^k_1(\QQ)$; see e.g.~\cite[Section 3]{eiza;theoryjacobi}.
	 Its Fourier expansion is
	 \bes
	 \phi_m(\tau,0)=\sum_{n=0}^\infty\Big(\sum_r c_{n,r}(\phi)\Big) q^n.
	 \ees
	 The previous sums over $r$ are finite, because $c(n,r)\ne 0$ implies $r^2\le 4nm$. Since the finite sum $\sum_r c_{n,r}|_{J_{k,m}^{\text{cusp}}(\QQ)}$ extracts the $n$-th Fourier coefficient of the \emph{elliptic} modular form $\phi(\tau,0)$ for any Jacobi cusp form~$\phi$ and any $n\ge 1$, it is enough to apply Lemma~\ref{lemma;propbmc} to such sum of functionals to conclude the proof.
	 \end{proof}	
	
	\section{Siegel modular forms of genus~$2$}\label{sec;siegmodforms}
	
	We briefly recall Siegel modular forms, which are the counterpart of elliptic modular forms in several variables.
	For the aim of this paper, we treat only the genus~$2$ case.
	
	The \emph{Siegel upper-half space} $\HH_2$ is the set of~$2\times 2$ symmetric matrices over~$\CC$ with positive definite imaginary part. It is a simply connected open subset of~$\CC^3$. The symplectic group~$\Sp_4(\RR)$ acts on~$\HH_2$ as a group of automorphisms by
	\bes
	g\colon Z\longmapsto g\cdot Z=(AZ+B)(CZ+D)^{-1},
	\ees
	for every $Z\in\HH_2$, where we decompose $g\in\Sp_4(\RR)$ in $2\times 2$ matrices as $g=\left(\begin{smallmatrix}
	A & B\\ C & D
	\end{smallmatrix}\right)$.
	Let $k\ge 4$ be an even integer. The symplectic group $\Sp_4(\RR)$ acts also on the space of complex-valued functions $F\colon\HH_2\to\CC$ via the so-called \emph{$|_k$-operator}, defined as
	\bes
	(F|_kg)(Z)=\det(CZ+D)^{-k}F(g\cdot Z),
	\ees
	for every $g\in \Sp_4(\RR)$.
	A \emph{Siegel modular form of weight $k$ (and genus~$2$)} is a holomorphic function $F\colon\HH_2\to\CC$ that satisfies the transformation law
	\bes
	F|_k\gamma=F,\quad\text{for every $\gamma\in\Sp_4(\ZZ)$}.
	\ees
	We denote the finite-dimensional complex vector space of these forms by~$M^k_2$. By the Koecher Principle, every Siegel modular form admits a Fourier expansion.
	We denote by~$\halfint_2$ the set of symmetric half-integral positive semi-definite matrices of order~$2$, namely
	\bes
	\halfint_2=\big\{T=\left(\begin{smallmatrix}
	n & r/2\\
	r/2 & m \end{smallmatrix}\right)\,:\, \text{$n,r,m\in\ZZ$ and $T\ge0$}\big\},
	\ees
	and by~$\halfint_2^+$ the subset of matrices which are \emph{positive definite}.
	The Fourier expansion of any~${F\in M^k_2}$ is indexed over~$\halfint_2$ as
	\be
	F(Z)=\sum_{T\in\halfint_2}c_T(F) e^{2\pi i\trace(TZ)}.
	\ee
	The complex numbers~$c_T(F)$ are the \emph{Fourier coefficients} of~$F$.
	If the Fourier expansion is supported on $\halfint_2^+$, then $F$ is called a \emph{Siegel cusp form}.
	We denote the subspace of cusp forms in $M^k_2$ by $S^k_2$.
	
	The group $\GL_2(\ZZ)$ acts on $\halfint_2$ via the action $T\mapsto u^t\cdot T\cdot u$, where~$u\in\GL_2(\ZZ)$ and~${T\in\halfint_2}$, preserving $\halfint_2^+$.
	The Fourier coefficients of Siegel modular forms of \emph{even} weight are invariant with respect to this action, namely $c_T(F)=c_{ u^t\cdot T\cdot u}(F)$ for every $F\in M^k_2$.
	We say that a matrix $T=\big(\begin{smallmatrix}
	n & r/2\\
	r/2 & m
	\end{smallmatrix}\big)\in\halfint_2$ is \emph{reduced} if $0\le r\le m\le n$.
	\begin{rem}\label{rem:reducedmat}
	The orbit of the subset of reduced matrices via the action of $\GL_2(\ZZ)$ is the whole $\halfint_2$.
	For this reason, the study of the Fourier coefficients of Siegel modular forms~(of even weight) restricts to the ones associated to reduced matrices.
	\end{rem}
	Our definition of reduced matrix is slightly different from the one in the literature.
	In fact, the reduced matrices are usually constructed to be representatives in~$\halfint_2$ with respect to the action of~$\SL_2(\ZZ)$, fulfilling the weaker condition~$|r|\le m\le n$.
	In our case, in virtue of Remark~\ref{rem:reducedmat}, we may consider the action of the whole~$\GL_2(\ZZ)$ on~$\halfint_2$.
	In particular, we may suppose $r$ to be non-negative.
	
	In analogy with the case of elliptic and Jacobi modular forms, the spaces $M^k_2$ and~$S^k_2$ admit a basis of Siegel modular forms with rational Fourier coefficients.
	We denote the~$\QQ$-vector spaces generated by these bases by $M^k_2(\QQ)$ and $S^k_2(\QQ)$, respectively.
	
	The dual space $M^k_2(\QQ)^*$ is generated by the \emph{Siegel coefficient extraction functionals}~$c_T$, defined for every~$T\in\halfint_2$ as
	\bes
	c_T\colon M^k_2(\QQ)\longrightarrow\QQ,\quad F\longmapsto c_T(F).
	\ees
	
	An important feature of the Siegel modular forms is that their Fourier expansions can be rewritten via Jacobi modular forms.
	That is, every $F\in M^k_2$ admits a \emph{Fourier--Jacobi expansion}
	\be\label{eq:FJdec}
	F(Z)=\sum_{m=0}^\infty \phi_m(\tau,z) e^{2\pi im\tau'},
	\ee
	where $Z=\left(\begin{smallmatrix}\tau & z\\ z & \tau'\end{smallmatrix}\right)\in\HH_2$, and $\phi_m\in J_{k,m}$ is the \emph{$m$-th Fourier--Jacobi coefficient} of $F$.
	Whenever we want to highlight that $\phi_m$ is a coefficient of $F\in M^k_2$, we write $\phi^F_m$.
	Clearly, if $F\in M^k_2(\QQ)$, then $\phi_m\in J_{k,m}(\QQ)$, and if $F\in S^k_2$, then $\phi_m\in J_{k,m}^\text{cusp}$.
	Furthermore, if~$T=\big(\begin{smallmatrix}
	n & r/2\\
	r/2 & m \end{smallmatrix}\big)$, then the $T$-th Fourier coefficient of $F$ coincides with one of the Fourier coefficients of $\phi_m$, more precisely~$c_T(F)=c_{n,r}(\phi_m)$.
	
	\subsection{Siegel Eisenstein series}\label{sec;siegeis}
	
	This section is a focus on the Siegel Eisenstein series~$E^k_2$ of genus~$2$ and even weight~$k\ge 4$. We deal with the Fourier coefficients $a^k_2(T)$ of $E^k_2$ associated to positive definite matrices $T$ and certain ratios of the form $a^k_2\big(\begin{smallmatrix}
	n & r/2t\\
	r/2t & m/t^2
	\end{smallmatrix}\big)\big/a^k_2\big(\begin{smallmatrix}
	n & r/2\\
	r/2 & m
	\end{smallmatrix}\big)$, for some positive $t$.
	The possible limits of these ratios, where $t$ is fixed and with respect to sequences of matrices of increasing determinant, are essential to classify the accumulation rays of the cone generated by the coefficient extraction functionals indexed over $\halfint_2^+$, and are extensively used in Section \ref{subsec;accrays}.
	
	\begin{defi}
	Let $P_0$ be the Siegel parabolic subgroup of $\Sp_4(\ZZ)$. The (normalized) \emph{Siegel Eisenstein series} of even weight $k\ge 4$ is defined as
	\bes
	E^k_2:\HH_2\longrightarrow \CC,\quad
	Z \longmapsto \sum_{\big(\begin{smallmatrix}A & B\\ C & D\end{smallmatrix}\big)\in P_0\backslash \Sp_4(\ZZ)}\det(CZ+D)^{-k}.
	\ees
	\end{defi}
	It is well-known that $E^k_2$ is a Siegel modular form of weight $k$. We denote its Fourier expansion by
	\bes
	E^k_2(Z)=\sum_{T\in\halfint_2}a^k_2(T)e^{2\pi i\trace(TZ)}.
	\ees
	We reserve the special notation $a^k_2(T)$ for the Fourier coefficients of $E^k_2$, instead of~$c_T(E^k_2)$, since they play a key role in the whole theory.
	
	To state the Coefficient Formula of $a^k_2(T)$, we need to recall some definitions.
	An integer~$D$ is said to be a \emph{fundamental discriminant} if either~$D\equiv 1$ mod~$4$ and squarefree, or~$D\equiv 4s$ for some squarefree integer~$s\equiv 2$ or~$3$ mod~$4$.
	Its associated Dirichlet character~$\chi_D$ is the one given by the Kronecker symbol~$\left(\frac{D}{\boldsymbol{\cdot}}\right)$.
	\begin{defi}[See {\cite[Section 2]{Co75}}]\label{defi;cohenH}
	Let $r$ and $N$ be non-negative integers, with $r$ positive. The Cohen~$H$-function~$H(r,N)$ is defined as follows. If~$N>0$ and~$(-1)^rN\equiv 0$ or~$1$ mod~$4$, we decompose~$N=-Dc^2$ with~$D$ a fundamental discriminant. In this case we set
	\bes
	H(r,N)=L(1-r,\chi_D)\sum_{d|c}\mu(d)\chi_D(d)d^{r-1}\sigma_{2r-1}(c/d).
	\ees
	If $N=0$, then $H(r,0)=\zeta(1-2r)$.
	\end{defi}
	
	\begin{lemma}[Coefficient Formula, see {\cite[p.\ 80]{eiza;theoryjacobi}}]\label{lemma;coeffourmeis}
	The Fourier coefficients of the Siegel Eisenstein series $E^k_2$ are rational and given by
	\be\label{eq:coefffor}
	a^k_2(T)=\begin{cases}
	\frac{2}{\zeta(1-k)\zeta(3-2k)}\sum_{d|(n,r,m)}d^{k-1}H\left(k-1,\frac{4\det T}{d^2}\right), & \text{if $T\neq 0$,}\\
	1, & \text{if $T=0$,}
	\end{cases}
	\ee
	for any $T=\big(\begin{smallmatrix}
	n & r/2\\
	r/2 & m
	\end{smallmatrix}\big)\in\halfint_2$.
	\end{lemma}
	The value $a^k_2\left(\begin{smallmatrix}
	n & 0\\ 0 & 0
	\end{smallmatrix}\right)$ coincide with the $n$-th Fourier coefficient of the (normalized) elliptic Eisenstein series \eqref{eq;normeisserelliptic}.
	The following lemma summarizes well-known properties of~$a^k_2(T)$.
	\begin{lemma}\label{classsiegeis}
	The Fourier coefficients of the Siegel Eisenstein series $E^k_2$ satisfy the following properties.
	\begin{enumerate}[label=(\roman*)]
	\item Suppose that $k\equiv 2$ mod $4$ and $T\in\halfint_2\setminus\{0\}$. If $\det T>0$, resp.\ $\det T=0$, then~$a^k_2(T)$ is a positive, resp.\ negative, rational number. \label{lem;pt1class}
	\item Suppose that $k\equiv 0$ mod $4$ and $T\in\halfint_2\setminus\{0\}$. If $\det T>0$, resp.\ $\det T=0$, then~$a^k_2(T)$ is a negative, resp.\ positive, rational number. \label{lem;pt2class}
	\item There exist positive constants $c_1$ and $c_2$ such that \label{classsiegeis;magn}
	\bes
	c_1\det(T)^{k-3/2}<|a^k_2(T)|<c_2\det(T)^{k-3/2},\quad \text{for every $T>0$}.
	\ees
	\end{enumerate}
	\end{lemma}
	We will usually refer to Lemma \ref{classsiegeis} \ref{classsiegeis;magn} saying that~$a^k_2(T)$ has \emph{the same order of magnitude} of~$\det(T)^{k-3/2}$, usually abbreviated as $a^k_2(T)\asymp\det(T)^{k-3/2}$.
	\begin{proof}
	The proof of first two points is a simple check using the Coefficient Formula.
	The idea is to show that all values of the $H$-function appearing in Formula \eqref{eq:coefffor} have the same sign if $\det T>0$ (resp.\ $\det T=0$).
	This can be proved by induction on the number of prime factors of $4\det(T)/d^2$, or via the equivalent definition of the $H$-function given in~\cite[Section 2]{Co75}. We follow the latter argument.
	\begin{enumerate}
	\item[\emph{\ref{lem;pt1class}}] Suppose that $\det T>0$, then $4\det T\equiv 0$ or $-1$ mod $4$.
	Decompose the $H$-function in $h$-functions as in \cite[Section 2]{Co75}, that is
	\bes
	H(k-1,4\det T)=\sum_{d^2|4\det T}h(k-1,4\det T/d^2).
	\ees
	Under the hypothesis that~$k\equiv 2$ mod~$4$, the $h$-functions are defined as
	\bes
	\qquad h(k-1,4\det T)=(k-2)! 2^{2-k}\pi^{1-k}(4\det T)^{k-3/2}L(k-1,\chi_{-4\det T}),
	\ees
	for every $T\in \halfint^+_2$.
	Clearly, the sign of $h(k-1,4\det T)$ depends on the sign of the last factor, which is positive since
	\bas
	\quad\qquad L(k-1,\chi_{-4\det T})\coloneqq\sum_{n=1}^\infty \frac{\chi_{-4\det T}(n)}{n^{k-1}}=\prod_p\frac{1}{1-\frac{\chi_{-4\det T}(p)}{p^{k-1}}}\ge\prod_p\frac{1}{1+\frac{1}{p^{k-1}}}>0.
	\eas
	Suppose now $\det T=0$, then $H(k-1,0)=\zeta(3-2k)$ and
	\[
	a^k_2(T)=\frac{2}{\zeta(1-k)}\sigma_{k-1}\big(\gcd(n,r,m)\big).
	\]
	Since $\zeta(1-k)=(-1)^{k-1}B_k/k$ and $k-1\equiv 1$ mod $4$, where $B_k$ is the $k$-th Bernoulli number, the coefficient $a^k_2(T)$ is negative.
	\item[\emph{\ref{lem;pt2class}}] It is analogous to the previous one. If~$k\equiv 0$ mod~$4$, then the decomposition in~$h$-functions is as above but with a factor of~$-1$, changing the sign of~$H(k-1,4\det T)$, for every.
	\item[\emph{\ref{classsiegeis;magn}}] This is well-known; see e.g.\ {\cite[Remark 2.2]{das;survey}}.\qedhere
	\end{enumerate}
	\end{proof}
	\begin{rem}\label{rem;siegcuspgrow}
	Let $k\ge 4$ be an even integer and let $F\in S^k_2$ be a Siegel cusp form.
	Suppose that $\seq{T_j}{j\in\NN}$ is a sequence of matrices in $\halfint_2$ of increasing determinant, that is, such that~$\det T_j\to+\infty$ when $j\to+\infty$.
	As explained e.g.\ in~\cite[Section 1.1.1]{das;survey}, the growth of the Fourier coefficients $\seq{c_{T_j}(F)}{j\in\NN}$ is estimated by the \emph{Hecke bound} as
	\bes
	c_{T_j}(F)=O_F\big(\det(T_j)^{k/2}\big).
	\ees
	By Lemma \ref{classsiegeis} \ref{classsiegeis;magn}, we deduce that the Fourier coefficients~$a^k_2(T_j)$ of the Siegel Eisenstein series~$E^k_2$ grow faster than~$c_{T_j}(F)$ when $j\to\infty$, for every cusp form~$F\in S^k_2$ and for every sequence of matrices~$\seq{T_j}{j}$ of increasing determinant.
	\end{rem}
	
	\subsection{Siegel series and ratios of Fourier coefficients}\label{subsec;ratFCeis}
	
	The aim of this section is to provide a classification of certain quotients of coefficients of Siegel Eisenstein series, and their limits over sequences of matrices with increasing determinant.
	The idea is to simplify the explicit formulas of these ratios using the so-called \emph{Siegel series}.
	These results will play a key role in Section \ref{sec;accrays} and Section \ref{sec;propC}, namely to classify the accumulation rays of the modular cone $\cone_k$ and to translate the polyhedrality of $\cone_k$ in terms of weight $k$ Jacobi cusp forms.
	We suggest the reader to skip this rather technical section during the first reading.
	
	We begin with an introduction on Siegel series. If $a$ is a non-zero integer, we denote by~$\nu_p(a)$ the maximal power of $p$ dividing $a$.
	\begin{defi}
	Let $T=\big(\begin{smallmatrix}
	n & r/2\\
	r/2 & m
	\end{smallmatrix}\big)\in\halfint^+_2$ and let $D$ be the fundamental discriminant such that $4\det T=-Dc^2$.
	For every prime $p$, we define~$\alpha_1(T,p)=\nu_p\big(\gcd(n,r,m)\big)$ and~$\alpha(T,p)=\nu_p(-4\det T/D)/2=\nu_p(c)$. The \emph{local Siegel series $F_p(T,s)$} is defined as
	\bes
	F_p(T,s)=\sum_{\ell=0}^{\alpha_1(T,p)}p^{\ell(2-s)}\bigg(\sum_{w=0}^{\alpha(T,p)-\ell}p^{w(3-2s)}-\chi_D(p)p^{1-s}\sum_{w=0}^{\alpha(T,p)-\ell-1}p^{w(3-2s)}\bigg),
	\ees
	where $s\in\CC$ and $\chi_D(n)=\left(\frac{D}{n}\right)$ is the Dirichlet character associated to the Kronecker symbol $\left(\frac{D}{\boldsymbol{\cdot}}\right)$.
	\end{defi}
	Conventionally, any sum from zero to a negative number is zero.
	We remark that if~$p$ does not divide~$-4\det T$, then~$F_p(T,s)=1$.
	
	Sometimes, in the literature, the definition of the local Siegel series differs from ours by a factor, more precisely it is defined as
	\[
	b_p(T,s)=\gamma_p(T,s)F_p(T,s),\quad\text{where}\quad
	\gamma_p(T,s)=\frac{(1-p^{-s})(1-p^{2-s})}{1-\chi_D(p)p^{1-s}};
	\]
	see \cite{Kat;expsiegser} and \cite[p.\ 473, Hilfssatz 10]{kauf}. For our purposes, the factor $\gamma_p(T,s)$ plays no role.
	\begin{defi}
		Let $T\in\halfint^+_2$. The \emph{Siegel series $F_T(s)$} is the product of local Siegel series
		\[
		F_T(s)=\prod_{p|4\det T}F_p(T,s).
		\]
	\end{defi}
	Using Siegel series, we may rewrite some of the Fourier coefficients $a^k_2(T)$ of Siegel Eisenstein series, as stated in the following result; see \cite{kauf}. 
	\begin{prop}\label{prop;newcoefform}
	Let $T\in\halfint^+_2$ and let $k\ge 4$ be an even integer. We may rewrite the Fourier coefficient of the Siegel Eisenstein series~$E^k_2$ associated to the matrix~$T$ as
	\[
	a^k_2(T)=\frac{2 L(2-k,\chi_D)}{\zeta(1-k)\zeta(3-2k)}\cdot F_T(3-k),
	\]
	where $D$ is the fundamental discriminant such that $4\det T=-Dc^2$.
	\end{prop}
	\textcolor{\myblack}{In the remaining part of this section we provide} some results on quotients of certain Fourier coefficients of~$E^k_2$ and their possible limits, as previously announced.
To simplify the explanation, for every positive integer $t$ and every
	\be\label{not;matTdiv}
	T=\left(\begin{smallmatrix}
	n & r/2\\
	r/2 & m
	\end{smallmatrix}\right)\in\halfint^+_2,\quad \text{we define}\quad \Tdiv{T}{t}=\left(\begin{smallmatrix}
	n & r/2t\\
	r/2t & m/t^2
	\end{smallmatrix}\right). 
	\ee
	\begin{lemma}\label{lem;lemratnew}
	Let $T=\big(\begin{smallmatrix}
	n & r/2\\
	r/2 & m
	\end{smallmatrix}\big)$ be a matrix in $\halfint^+_2$ and let $k\ge 4$ be an even integer.
	If~$t$ is a positive integer such that $t\neq 1$, $t|r$ and $t^2|m$, then
	\begin{equation}\label{eq;lemratnew}
	\frac{a^k_2(\Tdiv{T}{t})}{a^k_2(T)}=\prod_{p|t}\frac{F_p(\Tdiv{T}{p^{\nu_p(t)}},3-k)}{F_p(T,3-k)}.
	\end{equation}
Moreover $0< a^k_2(\Tdiv{T}{t})/a^k_2(T)< 1$.
	\end{lemma}
	\begin{proof}
	Let $D$ be the fundamental discriminant such that $4\det \Tdiv{T}{t}=-Dc^2$, then $4\det T=4 t^2 \det \Tdiv{T}{t}=-D (tc)^2$, hence the fundamental discriminants associated to~$T$ and~$\Tdiv{T}{t}$ are equal.
	We use Proposition~\ref{prop;newcoefform} to deduce
	\begin{align}\label{eq;prquotrat}
	\frac{a^k_2(\Tdiv{T}{t})}{a^k_2(T)}=\frac{F_{\Tdiv{T}{t}}(3-k)}{F_T(3-k)}=
	\prod_{p|t}\frac{F_p(\Tdiv{T}{t},3-k)}{F_p(T,3-k)}\cdot\prod_{p\nmid t}\frac{F_p(\Tdiv{T}{t},3-k)}{F_p(T,3-k)}.
	\end{align}
	
	Let $p$ be a prime such that $p$ does not divide $t$, then
	\[
	\alpha_1(T,p)=\nu_p(\gcd(n,r,m))=\nu_p(\gcd(n,r/t,m/t^2)=\alpha_1(\Tdiv{T}{t},p).
	\]
	Analogously, we deduce $\alpha(T,p)=\alpha(\Tdiv{T}{t},p)$. This implies that~$F_p(\Tdiv{T}{t},3-k)=F_p(T,3-k)$ for every $p$ which does not divide $t$, hence the last factor in \eqref{eq;prquotrat} simplifies to $1$.
	
	Suppose now that $p$ divides $t$.
	Since $\alpha(\Tdiv{T}{t},p)=\alpha(\Tdiv{T}{p^{\nu_p(t)}},p)$ and
	\bes
	\alpha_1(\Tdiv{T}{t},p)=\nu_p\big(\gcd(n,r/t,m/t^2)\big)=\nu_p\big(\gcd(n,r/p^{\nu_p(t)},m/p^{2\nu_p(t)})\big)=\alpha_1(\Tdiv{T}{p^{\nu_p(t)}},p),
	\ees
	we deduce that \eqref{eq;prquotrat} simplifies to \eqref{eq;lemratnew}.	
	Furthermore, since $k>4$, the value $F_p(T,3-k)$ is \emph{positive} for every~$T\in\halfint^+_2$.
	Moreover, since $\alpha_1(\Tdiv{T}{t},p)\le\alpha_1(T,p)$ and $\alpha(\Tdiv{T}{t},p)=\nu_p(c)$ is less than $\alpha(T,p)=\nu_p(ct)$, then $F_p(\Tdiv{T}{t},3-k)<F_p(T,3-k)$.
	This concludes the proof.
	\end{proof}
	We want to classify all possible limits of ratios of the form~\eqref{eq;lemratnew}, indexed over a sequence of matrices $\seq{T_j}{j\in\NN}$ in~$\halfint^+_2$, with increasing determinant and fixed bottom-right entry.
	To do so, we need to define certain \emph{special limits} associated to such families.
	For the purposes of this paper, we may consider only \emph{reduced} matrices.
	\begin{prop}\label{prop;newspeclim}
	Let $k\ge 4$ be even and let $m$ be a positive integer. Consider a sequence of reduced matrices $\seqBig{T_j=\big(\begin{smallmatrix}
	n_j & r_j/2\\
	r_j/2 & m
	\end{smallmatrix}\big)}{j\in\NN}$ in $\halfint^+_2$, of increasing determinant.
	Suppose that a prime $p$ is chosen such that $p^s|r_j$ and~$p^{2s}|m$ for some positive integer~$s$.
	If the sequence of ratios~$a^k_2(\Tdiv{T_j}{p^s})/a^k_2(T_j)$ converges to a value~$\lambda_{p^s}$ and~$\alpha(T_j,p)$ diverges when~${j\to\infty}$, then the sequence $\seqBig{\big(\alpha_1(\Tdiv{T_j}{p^s},p),\alpha_1(T_j,p)\big)}{j\in\NN}$ is eventually constant and
	\begin{equation}\label{eq;2newlim}
	\lambda_{p^s}=p^{s(3-2k)}\cdot\frac{1-p^{(2-k)(\alpha_1(\Tdiv{T_j}{p^s},p)+1)}}{1-p^{(2-k)(\alpha_1(T_j,p)+1)}},\qquad\text{for $j$ large enough.}
	\end{equation}
	\end{prop}
	We remark that for different values of $\alpha_1(T_j,p)$ and $\alpha_1(\Tdiv{T_j}{p^s},p)$, the ratio \eqref{eq;2newlim} assumes different values.  
	\begin{defi}\label{def;setspeclim}
	Let $k\ge 4$ be even and let $m$ be a positive integer. For all positive integers $s$ and all primes $p$ such that $p^{2s}$ divides $m$, the \emph{special limits (of weight $k$ and index~$m$) associated to $p^s$} are the limits of ratios arising as in Proposition~\ref{prop;newspeclim}.
	We denote by~$\limitmsp{m}{k}{p^s}$ the set of these special limits.
	\end{defi}
	As we are going to see with Proposition~\ref{prop;infmanylp}, the elements of~$\limitmsp{m}{k}{p^s}$ are those limits of ratios which can be obtained only asymptotically, since they are not ratios of Fourier coefficients of~$E^k_2$ arising from any matrix in~$\halfint^+_2$.
	For this reason, we call them ``special''.
	\begin{rem}\label{rem;finnuminlimitmspp}
	Let $k\ge 4$ be even and let $m$ be a positive integer.
	Since $\alpha_1(T_j,p)$ and~$\alpha_1(\Tdiv{T_j}{p^s},p)$ can assume only a finite number of values in \eqref{eq;2newlim}, the set~$\limitmsp{m}{k}{p^s}$ is \emph{finite} for every positive integer~$s$ and every prime~$p$ such that $p^{2s}$ divides $m$.
	\end{rem}
	\begin{proof}[Proof of Proposition \ref{prop;newspeclim}]
	The local Siegel series evaluated in $s=3-k$ is
	\[
	F_p(T,3-k)=\sum_{\ell=0}^{\alpha_1(T,p)}p^{\ell(k-1)}\Bigg(
	\underbrace{\sum_{w=0}^{\alpha(T,p)-\ell}p^{w(2k-3)}}_{(*)}-
	\chi_{D}(p)p^{k-2}\underbrace{\sum_{w=0}^{\alpha(T,p)-\ell-1}p^{w(2k-3)}}_{(**)}
	\Bigg).
	\]
	We remark that $(*)$ and $(**)$ are two different truncates of a geometric series. Since the truncate of a geometric series can be computed as $\sum_{i=0}^n r^i=\frac{1-r^{n+1}}{1-r}$ for every $r\neq 1$, then
	\begin{align*}
	F_p(T,&3-k)=\sum_{\ell=0}^{\alpha_1(T,p)}p^{\ell(k-1)}\Bigg(
	\frac{1-\big(p^{2k-3}\big)^{\alpha(T,p)-\ell+1}}{1-p^{2k-3}}-\chi_{D}(p)p^{k-2}\frac{1-\big(p^{2k-3}\big)^{\alpha(T,p)-\ell}}{1-p^{2k-3}}
	\Bigg)=\\
	=&\frac{1}{1-p^{2k-3}}\sum_{\ell=0}^{\alpha_1(T,p)}p^{\ell(k-1)}\Big(
	1-\chi_D(p)p^{k-2}+\big(\chi_D(p)p^{k-2}-p^{2k-3}\big)p^{(\alpha(T,p)-\ell)(2k-3)}
	\Big)=\\
	=&\frac{1}{1-p^{2k-3}}\Big(
	(1-\chi_D(p)p^{k-2})\underbrace{\sum_{\ell=0}^{\alpha_1(T,p)}p^{\ell(k-1)}}_{(\star)}+\\
	&\qquad\qquad\qquad\qquad\qquad+(\chi_D(p)p^{k-2}-p^{2k-3})p^{\alpha(T,p)(2k-3)}\underbrace{\sum_{\ell=0}^{\alpha_1(T,p)}p^{\ell(2-k)}}_{(\star\star)}
	\Big).
	\end{align*}
The terms $(\star)$ and $(\star\star)$ are truncates of two different geometric series. Computing their values, we deduce
\ba\label{eq;splitbigproof}
\begin{split}
F_p(T,3-k)&=\frac{1}{1-p^{2k-3}}\bigg(
\underbrace{\frac{\big(1-\chi_D(p)p^{k-2}\big)\big(1-p^{(k-1)(\alpha_1(T,p)+1)}\big)}{1-p^{k-1}}}_{(\clubsuit)}-\\
&-\underbrace{p^{\alpha(T,p)(2k-3)}}_{(\spadesuit)}
\underbrace{\big(
p^{2k-3}-\chi_D(p)p^{k-2}
\big)
\frac{1-p^{(2-k)(\alpha_1(T,p)+1)}}{1-p^{2-k}}}_{(\clubsuit\clubsuit)}
\bigg).
\end{split}
\ea

Let $\seq{T_j}{j\in\NN}$ be a sequence of reduced matrices in $\halfint^+_2$ with bottom-right entry fixed to~$m$ and increasing determinant, such that $\alpha(T_j,p)\to\infty$ when $j\to\infty$.
We want to study the asymptotic behavior of $F_p(T_j,3-k)$ with respect to $j\to\infty$ via~\eqref{eq;splitbigproof}.
The terms~$\mathsmaller{(\clubsuit)}$ and~$\mathsmaller{(\clubsuit\clubsuit)}$ are independent from $\alpha(T_j,p)$, and they remain bounded since~$\alpha_1(T_j,p)$ and~$\chi_{D_j}(p)$ assume only a finite number of values.
In contrast, the value of $\mathsmaller{(\spadesuit)}$ diverges if~$j\to\infty$, since $k>4$ by hypothesis.
This implies that
\begin{equation}\label{eq;splitbigproof2}
F_p(T_j,3-k)\sim p^{\alpha(T_j,p)(2k-3)}\cdot\frac{\big(
p^{2k-3}-\chi_{D_j}(p)p^{k-2}
\big)\big(1-p^{(2-k)(\alpha_1(T_j,p)+1)}\big)}{\big(p^{2k-3}-1\big)\big(1-p^{2-k}\big)},
\end{equation}
if $j\to\infty$.

We conclude the proof studying the asymptotic behavior of the ratios $a^k_2(\Tdiv{T_j}{p^s})/a^k_2(T_j)$. We compute these ratios via the local Siegel series and~\eqref{eq;splitbigproof2}, deducing
\begin{equation*}
\frac{a^k_2(\Tdiv{T_j}{p^s})}{a^k_2(T_j)}=\frac{F_p(\Tdiv{T_j}{p^s},3-k)}{F_p(T_j,3-k)}\sim p^{(3-2k)(\alpha(T_j,p)-\alpha(\Tdiv{T_j}{p^s},p)}\frac{1-p^{(\alpha_1(\Tdiv{T_j}{p^s},p)+1)(2-k)}}{1-p^{(\alpha_1(T_j,p)+1)(2-k)}},
\end{equation*}
when $j\to\infty$.
Since $\alpha(T_j,p)=\alpha(\Tdiv{T_j}{p^s},p)+s$, the claim follows.
\end{proof}
	\begin{cor}\label{cor;ratwithsiegser}
	Let $k\ge 4$ be an even integer and let $\seq{T_j}{j\in\NN}$ be a sequence of reduced matrices in $\halfint^+_2$ of increasing determinant, of the form~$T_j=\big(\begin{smallmatrix}
	n_j & r_j/2\\
	r_j/2 & m
	\end{smallmatrix}\big)$, where $m$ is a fixed positive integer.
	Suppose that a prime~$p$ is chosen such that~$p^{2s}|m$ for some positive integer~$s$. 
	There exists a positive constant~$C_{p^s}$ such that if
	\begin{equation}\label{eq;ratwithsiegser}
	\frac{a^k_2(\Tdiv{T_j}{p^s})}{a^k_2(T_j)}\xrightarrow[j\to\infty]{}\lambda_{p^s}
	\end{equation}
	for some $\lambda_{p^s}$, then either $\lambda_{p^s}=0$, and this happens only when the entries $r_j$ are eventually not divisible by $p^s$, or $C_{p^s}<\lambda_{p^s}<1$.
	Furthermore, if $\lambda_{p^s}$ is not a special limit in $\limitmsp{m}{k}{p^s}$, the sequence of ratios~$a^k_2(\Tdiv{T_j}{p^s})/a^k_2(T_j)$ is eventually constant equal to $\lambda_{p^s}$.
	\end{cor}
	\begin{proof}
	By Lemma \ref{classsiegeis}, the limit $\lambda_{p^s}$ is non-negative.
	If eventually $p^s\nmid r_j$, then the numerators of the ratios in \eqref{eq;ratwithsiegser} are eventually zero, and $\lambda_{p^s}=0$.
	From now on, we suppose that eventually $p^s$ divides $r_j$.
	
	The value of $F_p(T_j,3-k)$ depends only on $\alpha_1(T_j,p)$, $\alpha(T_j,p)$ and the fundamental discriminant $D_j$ such that $-4\det T_j=D_jc_j^2$.
	The value of $D_j$ influences $F_p(T_j,3-k)$ only via $\chi_{D_j}(p)$, which can assume only three values.
	Also the values of $\alpha_1(T_j,p)$ are finite, because $\alpha_1(T_j,p)=\nu_p(\gcd(n_j,r_j,m))$ with $m$ fixed. Only $\alpha(T_j,p)$ can diverge if~$j\to\infty$.
	If~$\alpha(T_j,p)$ does not diverge, then clearly there are only finitely many values that~${F_p(\Tdiv{T_j}{p^s},3-k)/F_p(T_j,3-k)}$ can assume, and they are strictly positive; see Lemma~\ref{lem;lemratnew}.
	In this case, the sequence of ratios $\seq{a^k_2(\Tdiv{T_j}{p^s}/a^k_2(T_j)}{j\in\NN}$ is eventually constant.
	If~$\alpha(T_j,p)$ diverges, then the limit~$\lambda_{p^s}$ is a \emph{special limit} in~$\limitmsp{m}{k}{p^s}$ by Proposition~\ref{prop;newspeclim}.
	\end{proof}
	\begin{defi}\label{defi;setoflimfortmfixp}
	Let $k$ and $m$ be positive integers, with $k\ge4$ even.
	For every positive integer $s$ and every prime $p$ such that $p^{2s}$ divides $m$, we denote by $\limitm{m}{k}{p^s}$ the set of all limits of ratios
	\bes
	\frac{a^k_2(\Tdiv{T_j}{p^s})}{a^k_2(T_j)}\xrightarrow[j\to\infty]{}\lambda_{p^s},
	\ees
	arising as in Corollary \ref{cor;ratwithsiegser}.
	\end{defi}
	We remark that $\limitmsp{m}{k}{p^s}\subseteq\limitm{m}{k}{p^s}\subset[0,1)\cap\QQ$ and that $0\in\limitm{m}{k}{p^s}$.
	The following result clarifies the structure of $\limitm{m}{k}{p^s}$.
	\begin{prop}\label{prop;infmanylp}
	Let $k$, $s$ and $m$ be positive integers, with $k\ge4$ even.
	Let~$p$ be a prime such that~$p^{2s}$ divides $m$. The set $\limitm{m}{k}{p^s}$ is \emph{infinite}, and splits into a \emph{disjoint union} as
	\be\label{eq;propinfmanylp}
	\limitm{m}{k}{p^s}=\left\{\frac{a^k_2(\Tdiv{T}{p^s})}{a^k_2(T)}\,:\,\text{$T\in\halfint^+_2$ reduced with bottom-right entry $m$}\right\}\coprod \limitmsp{m}{k}{p^s}.
	\ee
	In particular, the special limits in $\limitmsp{m}{k}{p^s}$ are not the values of ratios $a^k_2(\Tdiv{T}{p^s})/a^k_2(T)$ in~$\limitm{m}{k}{p^s}$ associated to \emph{reduced} matrices~$T\in\halfint^+_2$ with bottom-right entry~$m$.
	\end{prop}
	\begin{proof}
	The proof is divided in two steps.
	With the former, we prove that $\limitm{m}{k}{p^s}$ is infinite, and that the special limits in $\limitmsp{m}{k}{p^s}$ are never the value of a ratio $a^k_2(\Tdiv{T}{p^s})/a^k_2(T)$ for any \emph{reduced} matrix $T$ in $\halfint^+_2$ with bottom-right entry $m$.
	With the latter step, we prove that the values of such ratios are the elements of $\limitm{m}{k}{p^s}\setminus\limitmsp{m}{k}{p^s}$.
	
	\textbf{First step.} The idea is to find a sequence of sequences of matrices
	\be\label{eq;prseqseq}
	\seqBig{\seq{T_{j,0}}{j\in\NN},\seq{T_{j,1}}{j\in\NN},\dots,\seq{T_{j,x}}{j\in\NN},\dots}{x\in\NN},
	\ee
	where the $T_{j,x}$ are pairwise different reduced matrices of $\halfint^+_2$, such that for any fixed $x$ the sequence $\seq{T_{j,x}}{j\in\NN}$ is of increasing determinant, with
	\bes
	\alpha(T_{i,x},p)=\alpha(T_{j,x},p)\quad\text{and}\quad\alpha(T_{j,x},p)\neq\alpha(T_{j,y},p),
	\ees
	for every $i,j,x,y\in\NN$ with $x\neq y$.

	There exist infinitely many reduced matrices $M_x$ in $\halfint^+_2$ of increasing determinant and with pairwise different values of $\alpha(M_x,p)$.
	In fact, we may choose $M_x=\big(\begin{smallmatrix}
	p^{2x} & 0\\ 0 & m
	\end{smallmatrix}\big)$ with~${x\ge x_0}$, for some $x_0$ such that $p^{2x_0}\ge m$, for which we have
	\bes
	-4\det M_x=-4mp^{2x}=D(cp^x)^2,
	\ees
	where we decompose $-4m=Dc^2$ with $D$ a fundamental discriminant.
	It is clear that~$\alpha(M_x,p)=\nu_p(c)+x$ assumes different values for different choices of $x\ge x_0$.
	From any such $M_x$, we construct the family of reduced matrices
	\bes
	T_{j,x}=\begin{pmatrix}
	(p+1)^{2j}p^{2x} & 0\\ 0 & m
	\end{pmatrix},\qquad\text{where $j\in\NN$.}
	\ees
	Since $T_{0,x}=M_x$ for every $x\ge x_0$, we deduce that
	\ba\label{eq;proofmanyeqseqseq}
	-4\det T_{j,x}&=-4\det M_x\cdot(p+1)^{2j}=D(cp^x(p+1)^j)^2\\
	\alpha(T_{j,x},p)&=\nu_p\big(cp^x(p+1)^j\big)=\nu_p(cp^x)=\alpha(M_x,p)\\
	\alpha_1(T_{j,x},p)&=\nu_p\big(\gcd(p^{2x}(p+1)^{2j},m)\big)=\nu_p\big(\gcd(p^{2x},m)\big)=\alpha_1(M_x,p),
	\ea
	for every $j\in\NN$ and for every~$x\ge x_0$.
	Analogous equalities are satisfied with~$\Tdiv{T_{j,x}}{p^s}$ and~$\Tdiv{M_x}{p^s}$ in place of~$T_{j,x}$ and~$M_x$, respectively.
	By Lemma~\ref{lem;lemratnew}, the equalities \eqref{eq;proofmanyeqseqseq} imply that the sequence of ratios
	\be\label{eq;prconvalseq}
	\seqBig{\frac{a^k_2(\Tdiv{T_{j,x}}{p^s})}{a^k_2(T_{j,x})}}{j\in\NN}=\seqBig{\frac{F_p(\Tdiv{T_{j,x}}{p^s},3-k)}{F_p(T_{j,x},3-k)}}{j\in\NN}
	\ee
	is \emph{constant} for every $x\ge x_0$.
	This implies that the ratio~$a^k_2(\Tdiv{M_x}{p^s})/a^k_2(M_x)$ is an element of~$\limitm{m}{k}{p^s}$ for every $x\ge x_0$.
	
	Since $\alpha(M_x,p)\to\infty$ when $x\to\infty$, then by Proposition \ref{prop;infmanylp} we deduce that
	\be\label{eq;prconvasy}
	\frac{a^k_2(\Tdiv{M_x}{p^s})}{a^k_2(M_x)}\longrightarrow \lambda_{p^s}\in\limitmsp{m}{k}{p^s},\qquad\text{if $x\to\infty$},
	\ee
	that is, the value $\lambda_{p^s}$ is a \emph{special limit}.
	
	We are ready to prove that $\limitm{m}{k}{p^s}$ is infinite.
	Suppose that it is not.
	Then the number of values assumed by the constant sequences \eqref{eq;prconvalseq} with $x\ge x_0$ is finite.
	We deduce from~\eqref{eq;prconvasy} that there exist $\tilde{x}\ge x_0$ and a special limit $\lambda_{p^s}\in\limitmsp{m}{k}{p^s}$ such that~${F_p(\Tdiv{M_x}{p^s},3-k)/F_p(M_x,3-k)= \lambda_{p^s}}$ for every $x\ge \tilde{x}$.
	The rational number~$\lambda_{p^s}$, as fraction in lowest terms, has denominator always divisible by~$p$. In fact, we may rewrite~\eqref{eq;2newlim} as
	\bes
	\lambda_{p^s}=\frac{p^{(k-2)\big(\alpha_1(M_x,p)-\alpha_1(\Tdiv{M_x}{p},p)\big)}}{p^{s(2k-3)}}\cdot\frac{p^{(k-2)\big(\alpha_1(\Tdiv{M_x}{p^s},p)+1\big)}-1}{p^{(k-2)(\alpha_1(M_x,p)+1)}-1}.
	\ees
	This fraction, if reduced in lowest terms, has denominator divisible by $p$, since~$k\ge4$ and~$\alpha_1(M_x,p)-\alpha_1(\Tdiv{M_x}{p^s},p)\le 2s$.
	
	Since both $F_p(M_x,3-k)$ and $F_p(\Tdiv{M_x}{p^s},3-k)$ are integers, the power of $p$ dividing the denominator of $\lambda_{p^s}$, as fraction in lowest term, must eventually divide $F_p(M_x,3-k)$, for every $x\ge\tilde{x}$.
	This is not possible, since $F_p(M_x,3-k)$ is not divisible by $p$.
	In fact, under the hypothesis $k>4$, we deduce \textcolor{\myblack}{with} simple congruences modulo $p$ that
	\be\label{eq;prcongFp}
	F_p(M_x,3-k)\equiv 1-\chi_D(p)p^{k-2}(1-\delta_{0,\alpha_1(M_x,p)})\equiv 1\quad\text{mod $p$},
	\ee
	for every $x\ge\tilde{x}$.
	Hence $\limitm{m}{k}{p^s}$ must be infinite.
	
	Since \eqref{eq;prcongFp} is satisfied for every $T\in\halfint^+_2$ in place of $M_x$, we deduce that the special limits in $\limitmsp{m}{k}{p^s}$ can not be obtained as ratios~$a^k_2(T_{p^s})/a^k_2(T)$ for any~$T\in\halfint^+_2$ reduced with bottom-right entry $m$.
	
	\textbf{Second step.} Let $T=\big(\begin{smallmatrix}
	n & r/2\\
	r/2 & m
	\end{smallmatrix}\big)$ be a reduced matrix in $\halfint^+_2$.
	Consider the sequence of reduced matrices $\seq{T_j}{j\in\NN}$, where $T_j$ is defined as
	\bes
	T_j=\begin{pmatrix}
	n-j(r^2-4nm) & r/2\\
	r/2 & m
	\end{pmatrix}.
	\ees
	We remark that $T_0=T$, and that $\det T_j\to\infty$ when $j\to\infty$.
	We decompose $-4\det T=Dc^2$, where~$D$ is a fundamental discriminant and deduce that
	\be\label{eq;propinfmanylpdecj}
	-4\det T_j=r^2-4m\big(n-j(r^2-4nm)\big)=(r^2-4nm)(4mj+1)=Dc^2(4mj+1).
	\ee
	Let $\seq{T_x}{x}$ be the sub-sequence of $\seq{T_j}{j\in\NN}$ such that $4mx+1$ is a perfect square. We denote the latter by $c_x^2$, with $c_x$ positive.
	There are infinitely many natural numbers~$x$ satisfying this condition. In fact, we may choose~$x=y(my+1)$, where~$y$ is a positive integer, since in this case
	\bes
	4mx+1=4my(my+1)+1=(2my+1)^2.
	\ees
	We deduce from \eqref{eq;propinfmanylpdecj} that the matrices of the sequence $\seq{T_x}{x}$ satisfy
	\bes
	-4\det T_x=Dc^2(4mx+1)=D(c\cdot c_x)^2,
	\ees
	therefore
	\bes
	\alpha(T_x,p)=\nu_p(c\cdot c_x)=\nu_p(c)=\alpha(T_0,p),
	\ees
	since $c_x^2=4mx+1$ and $p$ divides $m$.
	
	We claim that $\alpha_1(T_x,p)=\alpha_1(T,p)$. To prove it, we firstly remark that
	\bes
	\alpha_1(T_x,p)=\nu_p\big(\gcd(n-xDc^2,r,m)\big)=\min\{\nu_p(n-xDc^2),\nu_p(r),\nu_p(m)\}.
	\ees
	Clearly $p^{\alpha_1(T,p)}|Dc^2$.
	If~${\alpha_1(T,p)=\nu_p(n)}$, then also $\alpha_1(T,p)=\nu_p(n-xDc^2)$ for every index~$x$ of the sub-sequence $\seq{T_x}{x}$.
	If~$\alpha_1(T,p)\neq\nu_p(n)$, then~$p^{\alpha_1(T,p)}|(n-xDc^2)$ for every~$x$, hence $\alpha_1(T,p)=\min\{\nu_p(n-xDc^2),\nu_p(r),\nu_p(m)\}$.
	These imply what we claimed above.
	
	Since $\alpha_1(\Tdiv{T_x}{p^s},p)=\alpha_1(\Tdiv{T}{p^s},p)$ and $\alpha(\Tdiv{T_x}{p^s},p)=\alpha(\Tdiv{T}{p^s},p)$ for every $x$, the sequence of ratios
	\bes
	\seqBig{\frac{a^k_2(\Tdiv{T_x}{p^s})}{a^k_2(T_x)}}{x}=\seqBig{\frac{F_p(\Tdiv{T_x}{p^s},3-k)}{F_p(T_x,3-k)}}{x}
	\ees
	is \emph{constant}.
	This implies that the value of the ratios $a^k_2(\Tdiv{T}{p^s})/a^k_2(T)$ is an element of~$\limitm{m}{k}{p^s}$.
	Since $T$ was chosen arbitrarily among the reduced matrices in $\halfint^+_2$ with~$m$ as bottom-right entry, the proof is concluded.
	\end{proof}
	The remaining part of this section aims to generalize the previous results, replacing the limits of ratios $\lambda_{p^s}$ by \emph{tuples} of limits of ratios, indexed over the positive integers $t$ such that $t^2$ divides $m$.
	\begin{cor}\label{cor;ratwithsiegsernonprime}
	Let $k\ge4$ be an even integer and let $\seqbig{T_j=\big(\begin{smallmatrix}
	n_j & r_j/2\\
	r_j/2 & m
	\end{smallmatrix}\big)}{j\in\NN}$ be a sequence of reduced matrices in $\halfint^+_2$ of increasing determinant, where the bottom-right entries are fixed to a positive integer $m$.
	Let $t$ be a positive integer such that~$t^2|m$ and that
	\bes
	\frac{a^k_2(\Tdiv{T_j}{t})}{a^k_2(T_j)}\xrightarrow[j\to\infty]{}\lambda_t
	\ees
	for some $\lambda_t$. There exists a positive constant $C_t$, depending on $t$, such that either~$\lambda_t=0$, and this happens only when the entries $r_j$ are eventually not divisible by $t$, or~$C_t<\lambda_t<1$. There exist also a sub-sequence $\seq{T_i}{i}$ of $\seq{T_j}{j\in\NN}$ and $\lambda_{p^{\nu_p(t)}}\in\limitm{m}{k}{p^{\textcolor{\myblack}{\nu_p(t)}}}$ for every prime divisor $p$ of $t$, such that
	\be\label{eq;statratwithsiegsernonprime}
	\frac{a^k_2\big(\Tdiv{T_i}{p^{\nu_p(t)}}\big)}{a^k_2(T_i)}\xrightarrow[i\to\infty]{}\lambda_{p^{\nu_p(t)}}\quad\text{and}\quad\lambda_t=\prod_{p|t}\lambda_{p^{\nu_p(t)}}.
	\ee
	Furthermore, if $\lambda_{p^{\nu_p(t)}}$ is a \emph{non-special limit} for every $p$, then the sequence $\seq{a^k_2(\Tdiv{T_i}{t})/a^k_2(T_i)}{i}$ is eventually equal to $\lambda_t$.
	\end{cor}
	\begin{proof}
	It is a consequence of Lemma \ref{lem;lemratnew}. The result follows as in Corollary \ref{cor;ratwithsiegser}, working on each factor appearing on the right-hand side of \eqref{eq;lemratnew} applied with $T_j$ in place of $T$.
	\end{proof}
	\begin{defi}\label{def;tupleoflim}
	Let $m$ be a positive integer and let $t_0=1<t_1<\dots<t_d$ be the divisors of $m$ such that $t_i^2|m$ for all $i$. We denote by $\limitmm{m}{k}$ the set of tuples of rational numbers $(\lambda_{t_1},\dots,\lambda_{t_d})$ for which there exists a sequence of reduced matrices $\seq{T_j}{j\in\NN}$ in $\halfint^+_2$, with increasing determinant and bottom-right entry $m$, such that
	\bes
	\frac{a^k_2(\Tdiv{T_j}{t_i})}{a^k_2(T_j)}\xrightarrow[j\to\infty]{}\lambda_{t_i},\qquad\text{for every $i=1,\dots,d$.}
	\ees
	We say that a tuple $(\lambda_{t_1},\dots,\lambda_{t_d})$ in $\limitmm{m}{k}$ is a \emph{special tuple of limits (of weight~$k$ and index~$m$)} if there exists a $i$ such that $t_i=p^s$ for some prime $p$ and some positive integer $s$, and such that $\lambda_{t_i}$ is a special limit in $\limitmsp{m}{k}{p^s}$.
	We denote by $\limitmmsp{m}{k}$ the set of these special tuples of limits.
	\end{defi}
	Let $(\lambda_{t_1},\dots,\lambda_{t_d})$ be a tuple of limits in $\limitmm{m}{k}$.
	If $t$ is a divisor of $m$ such that~$t^2|m$ and with prime decomposition~$t=p_1^{s_1}\cdots p_x^{s_x}$, then both~$t$ and the powers of primes of the prime decomposition of~$t$ appear among the~$t_i$'s. Moreover~$\lambda_{t}=\prod_{j=1}^x\lambda_{p_j^{s_j}}$ by Corollary \ref{cor;ratwithsiegsernonprime}.
	
	The tuples in $\limitmm{m}{k}$ will be used in Section \ref{sec;accrays} to index the accumulation rays of the modular cone $\cone_k$ associated to sequences of matrices with bottom-right entries fixed to $m$.
	
	We conclude this section with the following generalization of Proposition~\ref{prop;infmanylp}, which shows that the tuples of ratios of Siegel Eisenstein series associated to the same matrix in~$\halfint^+_2$ lie in some~$\limitmm{m}{k}$.
	\begin{cor}\label{cor;infmanytup}
	Let $k$ and $m$ be positive integers, with $k\ge4$ even.
	We denote by ${t_0=1<t_1<\dots<t_d}$ the divisors of $m$ such that $t_i^2|m$ for all~$i$.
	If~$d\ge 1$, i.e.~$m$ is non-squarefree, then the set~$\limitmm{m}{k}$ is \emph{infinite}, and splits into a disjoint union as
	\bes
	\limitmm{m}{k}\!=\!\bigg\{\bigg(\frac{a^k_2(\Tdiv{T}{t_1})}{a^k_2(T)},\dots,\frac{a^k_2(\Tdiv{T}{t_d})}{a^k_2(T)}\bigg):\text{$T\in\halfint^+_2$ reduced with bottom-right entry $m$}\bigg\}\coprod \limitmmsp{m}{k}.
	\ees
	Furthermore, if $m$ is divisible by the squares of two different primes, then also $\limitmmsp{m}{k}$ is infinite.
	\end{cor}
	\begin{proof}
	In the \emph{second step} of the proof of Proposition \ref{prop;infmanylp}, we proved that for every ${T=\big(\begin{smallmatrix}
	n & r/2\\
	r/2 & m
	\end{smallmatrix}\big)}$ in $\halfint^+_2$, the sequence of reduced matrices with increasing determinant
	\bes
	\seqBig{T_j=\big(\begin{smallmatrix}
	n-j(r^2-4nm) & r/2\\
	r/2 & m
	\end{smallmatrix}\big)}{j\in\NN}\subset\halfint^+_2
	\ees
	contains a sub-sequence $\seq{T_x}{x}$ such that $\seq{a^k_2\big(\Tdiv{T_x}{p^s}\big)/a^k_2(T_x)}{x}$ is a constant sequence.
	This implied that the value $a^k_2(\Tdiv{T}{p^s})/a^k_2(T)$ lies in $\limitm{m}{k}{p^s}$ for every reduced matrix $T$ in $\halfint^+_2$.
	
	We remark that the definition of the matrices $T_j$ does not depend on the chosen power of prime $p^s$.
	We recall that the matrices $T_x$ were chosen to be the ones in $\seq{T_j}{j\in\NN}$ such that~$4mj+1$ is a perfect square.
	Also this choice does not depend on the power of prime~$p^s$.
	This means that the sequence~$\seq{a^k_2\big(\Tdiv{T_x}{p^s}\big)/a^k_2(T_x)}{x}$ is constant \emph{for every}~$p^s$. By Corollary~\ref{cor;ratwithsiegsernonprime}, we deduce that the sequence of tuples
	\bes
	\seqbigg{\Big(\frac{a^k_2\big(\Tdiv{T_x}{t_1}\big)}{a^k_2(T_x)},\dots,\frac{a^k_2\big(\Tdiv{T_x}{t_d}\big)}{a^k_2(T_x)}\Big)}{x}
	\ees
	is \emph{constant}. This implies that
	\be\label{eq;proofcor319tuples}
	\bigg\{\bigg(\frac{a^k_2(\Tdiv{T}{t_1})}{a^k_2(T)},\dots,\frac{a^k_2(\Tdiv{T}{t_d})}{a^k_2(T)}\bigg)\,:\,\text{$T\in\halfint^+_2$ reduced with bottom-right entry $m$}\bigg\}\subseteq\limitmm{m}{k}.
	\ee
	
	Since the number of values of the ratios $a^k_2(\Tdiv{T}{t_1})/a^k_2(T)$, with $T$ in $\halfint^+_2$ reduced and with bottom-right entry $m$, is infinite by Proposition \ref{prop;infmanylp}, then also $\limitmm{m}{k}$ is infinite.
	
	Any tuple in $\limitmmsp{m}{k}$ has an entry which is a special limit $\lambda_{p^s}$ in $\limitm{m}{k}{p^s}$ associated to a power of a prime $p^s$ such that $p^{2s}|m$.
	By Proposition \ref{prop;infmanylp}, the limit $\lambda_{p^s}$ is not the value of a ratio $a^k_2(\Tdiv{T}{p^s})/a^k_2(T)$ for any reduced matrix $T$ in $\halfint^+_2$ with bottom-right entry $m$.
	This implies that the subset of $\limitmm{m}{k}$ appearing in~\eqref{eq;proofcor319tuples} is \emph{disjoint} with~$\limitmmsp{m}{k}$.
	
	We conclude the proof showing that if $m$ is divisible by the squares of two different primes, then also $\limitmmsp{m}{k}$ is infinite.
	We suppose without loss of generality that~$t_1$ and~$t_2$ are two different primes.
	We follow the same idea of the \emph{first step} in the proof of Proposition~\ref{prop;infmanylp}.
	For every $j$ and $x$ in $\NN$, let $T_{j,x}$ be the reduced matrix in $\halfint^+_2$ defined as
	\bes
	T_{j,x}=\begin{pmatrix}
	t_2^{2j}\cdot t_1^{2x} & 0\\
	0 & m
	\end{pmatrix}.
	\ees
	The tuple
	\be
	\bigg(
	\frac{a^k_2(\Tdiv{T_{j,x}}{t_1})}{a^k_2(T_{j,x})},\frac{a^k_2(\Tdiv{T_{j,x}}{t_2})}{a^k_2(T_{j,x})},\dots,\frac{a^k_2(\Tdiv{T_{j,x}}{t_d})}{a^k_2(T_{j,x})}
	\bigg)
	\ee
	is an element of $\limitmm{m}{k}$ for every $j,x\in\NN$, as we showed at the beginning of this proof.
	For every choice of $j$, the limit
	\bes
	\lambda_{t_1,j}=\lim_{x\to\infty}\frac{a^k_2(\Tdiv{T_{j,x}}{t_1})}{a^k_2(T_{j,x})}
	\ees
	is a special limit in $\limitmsp{m}{k}{t_1}$.
	We recall that $a^k_2(\Tdiv{T_{j,x}}{t_2})/a^k_2(T_{j,x})=a^k_2(\Tdiv{T_{j,\tilde{x}}}{t_2})/a^k_2(T_{j,\tilde{x}})$ for every~$x,\tilde{x}$ large enough.
	This was actually proven using \eqref{eq;proofmanyeqseqseq} in the proof of Proposition~\ref{prop;infmanylp}.
	Hence, there exist tuples in $\limitmmsp{m}{k}$ of the form
	\be\label{eq;proofsecentinf}
	\bigg(
	\lambda_{t_1,j},\frac{a^k_2(\Tdiv{T_{j,x}}{t_2})}{a^k_2(T_{j,x})},\dots
	\bigg)
	\ee
	for some $x$ large enough.
	By Proposition \ref{prop;infmanylp}, the number of values assumed by the second entry of \eqref{eq;proofsecentinf}, with $j\to\infty$, is infinite. In fact $\lim_{j\to\infty}a^k_2(\Tdiv{T_{j,x}}{t_2})/a^k_2(T_{j,x})$ is a special limit in $\limitmsp{m}{k}{t_2}$.
	This implies that there are infinitely many special tuples of limits of the form~\eqref{eq;proofsecentinf} in~$\limitmmsp{m}{k}$.
	\end{proof}
	
	\subsection{Klingen Eisenstein series}\label{sec;klingeisser}
	
	Any elliptic modular form $f\in M^k_1$ of even weight~$k$ can be written in a unique way as the sum of a multiple of the Eisenstein series~$E^k_1$, and a cusp form~$g\in S^k_1$, that is, there exists a complex numbers~$a$ such that~${f=aE^k_1+bg}$.
	The same decomposition holds for the Fourier coefficients of~$f$.
	Namely, we can decompose~${c_n(f)=a\cdot c_n(E^k_1)+ c_n(g)}$ for every natural number~$n$.
	It is well-known that the coefficients of $E^k_1$ grow faster than the coefficients of any cusp form, with respect to~$n\to\infty$.
	This means that if~$f$ is not a cusp form, then the relevant part for the growth of $c_n(f)$ is given by its Eisenstein part.
	Such a clean decomposition is characteristic of elliptic modular forms, and does not hold for Siegel modular forms.
	The main obstacles are the so-called \emph{Klingen Eisenstein series}, whose coefficient growths behave sometimes as for the Siegel Eisenstein series $E^k_2$ and sometimes as for cusp forms, depending on the chosen sequences of matrices of increasing determinant.
	In the recent paper \cite{boda;petnorm}, B\"{o}cherer and Das have proposed an extensive study of the growth of these coefficients.
	The aim of this section is to clarify the previous issue and to recall from \cite{boda;petnorm} the necessary results for the purposes of this paper.\\
	
	We denote by $C_{2,1}$ the Klingen parabolic subgroup of $\Sp_4(\ZZ)$, defined as
	\bes
	C_{2,1}=\big\{\gamma\in\Sp_4(\ZZ)\,:\,\gamma=\left(\begin{smallmatrix}
	* & *\\
	0_{1,3} & *
	\end{smallmatrix}\right)\big\}.
	\ees
	\begin{defi}
	Let $k> 4$ be an even integer. Given an elliptic cusp form $f\in S^k_1$, the \emph{Klingen Eisenstein series of weight $k$} attached to $f$ is defined as
	\bes
	E^k_{2,1}(f,Z)=\sum_{\gamma\in C_{2,1}\backslash\Sp_4(\ZZ)}\det(CZ+D)^{-k}f\big((\gamma\cdot Z)^*\big),
	\ees
	where we denote by $Z^*$ the upper-left entry of $Z\in\HH_2$, and where~$\gamma=\left(\begin{smallmatrix}
	A & B\\ C & D
	\end{smallmatrix}\right)$.
	\end{defi}
	It is well-known that Klingen Eisenstein series are Siegel modular forms.
	We use the special notation $a^k_2(f,T)$ for the Fourier coefficient of $E^k_{2,1}(f)$ associated to the matrix~$T$ in~$\halfint_2$.
	The subspace of Klingen Eisenstein series is denoted by $N^k_2$.
	This subspace has complex dimension equal to the one of $S^k_1$, and any basis is of the form $E^k_{2,1}(f_1),\dots,E^k_{2,1}(f_\ell)$ for some basis $f_1,\dots,f_\ell$ of $S^k_1$.
	Moreover, if $f\in S^k_1(\QQ)$, then also $E^k_{2,1}(f)$ has rational Fourier coefficients.
	\begin{rem}\label{rem;coefklingcuspeq}
	Any Fourier coefficient of $E^k_{2,1}(f)$ associated to a \emph{singular} matrix in $\halfint_2$ is equal to a coefficient of the elliptic cusp form $f$, as we briefly recall.
	If~$T\in\halfint_2$ is singular, then there exist $u\in \GL_2(\ZZ)$ and $n\in\NN$ such that
	\bes
	u^t\cdot T\cdot u = \left(\begin{smallmatrix}
	n & 0\\ 0& 0
	\end{smallmatrix}\right).
	\ees
	By Remark \ref{rem:reducedmat}, we deduce that
	\be\label{eq;inrefsingSmodform}
	a^k_2(f,T)=a^k_2\left(f,\left(\begin{smallmatrix}
	n & 0\\ 0& 0
	\end{smallmatrix}\right)\right).
	\ee
	It is well-known that the coefficient appearing on the right-hand side of~\eqref{eq;inrefsingSmodform} equals~$c_n(f)$; see e.g.\ \cite[Section 5, Proposition 5]{kli;intro}.
	\end{rem}
	The Structure Theorem for Siegel modular forms \cite[Theorem 2, p.\ 73]{kli;intro} allows us to decompose the space of Siegel modular forms in
	\be\label{eq;structthm}
	M^k_2=\langle E^k_2\rangle_\CC\oplus N^k_2\oplus S^k_2,
	\ee
	with analogous decomposition of~$M^k_2(\QQ)$ over~$\QQ$.
	We highlighted in Remark \ref{rem;siegcuspgrow} some bounds for the growth of the Fourier coefficients of $E^k_2$ and the cusp forms in $S^k_2$. We provide now the missing bounds for the Klingen Eisenstein series in $N^k_2$.
	The following result \textcolor{\myblack}{of Kitaoka~\cite{kit;modforms}} is a first attempt in this direction.
	\begin{prop}[\textcolor{\myblack}{See \cite[Theorem p.\ 113 and Corollary p.\ 120]{kit;modforms}}]\label{prop;growthKEvsSE}
	Let $k>4$ be an even integer and let $\seq{T_j}{j\in\NN}$ be a sequence of matrices in~$\halfint^+_2$ of increasing determinant. For every elliptic cusp form~${f\in S^k_1}$, the Fourier coefficients of the associated Klingen Eisenstein series~$E^k_{2,1}(f)$ satisfy the bound
	\bes
	a^k_2(f,T_j)=O\big(\det(T_j)^{k-3/2}\big),\quad\text{for $j\to\infty$.}
	\ees
	\end{prop}
	Proposition \ref{prop;growthKEvsSE}, jointly with \textcolor{\myblack}{Lemma~\ref{classsiegeis}~\ref{classsiegeis;magn}}, ensures that the Fourier coefficients of any Klingen Eisenstein series of weight $k>4$ do not grow faster than the coefficients of the Siegel Eisenstein series of the same weight.
	This is not enough for our purposes.
	In fact, we need to know with respect to which sequences~$\seq{T_j}{j\in\NN}$ in~$\halfint^+_2$ the coefficients~$a^k_2(f,T_j)$ grow with the same order of magnitude of~$a^k_2(T_j)$.
	We illustrate here a solution of this problem following the wording of~\cite{boda;petnorm}.
	
	Let $k>4$ be an even integer and let $f\in S^k_1$. We write the Fourier--Jacobi expansion of~$E^k_{2,1}(f)$ as $E^k_{2,1}(f,Z)=\sum_m\phi_m(\tau,z)e^{2\pi im\tau'}$, where $Z=\left(\begin{smallmatrix}
	\tau & z\\
	z & \tau'
	\end{smallmatrix}\right)\in\HH_2$.
	For every~$m$, the Fourier--Jacobi coefficient $\phi_m$ decomposes as a sum of its Eisenstein and cuspidal parts, respectively~$\mathcal{E}_{k,m}\in J^{\text{Eis}}_{k,m}$ and~$\phi^0_m\in J^{\text{cusp}}_{k,m}$, that is,
	\begin{equation}\label{JacobiEiscusp}
	\phi_m=\mathcal{E}_{k,m}+\phi^0_m.
	\end{equation}
	This implies the decomposition of Fourier coefficients
	\be\label{eq;decfourinjac}
	a^k_2(f,T)=c_{n,r}(\mathcal{E}_{k,m})+c_{n,r}(\phi^0_m),\quad\text{for every $T=\left(\begin{smallmatrix}
	n & r/2\\
	r/2 & m
	\end{smallmatrix}\right)\in\halfint_2$.}
	\ee
	
	The idea is to deduce the growth of $a^k_2(f,T)$ from the growth of the two members appearing on the right-hand side of~\eqref{eq;decfourinjac}.
	The next result connects the growth of the Eisenstein part $c_{n,r}(\mathcal{E}_{k,m})$ with the one of the coefficients of the Siegel Eisenstein series.
	\begin{prop}[See {\cite[Theorem 6.8]{boda;petnorm}}]\label{thm;bodapetn}
	The Eisenstein part of $a^k_2(f,T)$ appearing in~\eqref{eq;decfourinjac} can be further decomposed as
	\[
	c_{n,r}(\mathcal{E}_{k,m})=\frac{\zeta(1-k)}{2}\sum_{t^2|m} \alpha_m(t,f)\cdot a^k_2 \begin{pmatrix}n & r/2t\\ r/2t & m/t^2\end{pmatrix},
	\]
	where we use the usual convention that $a^k_2 \left(\begin{smallmatrix}
	n & r/2t\\
	r/2t & m/t^2
	\end{smallmatrix}\right)=0$ whenever $t$ does not divide $r$, and
	\be\label{eq;alphamaux}
	\alpha_m(t,f)=\sum_{\ell\mid t}\mu(t/\ell)\frac{\gf{f}{m/\ell^2}}{g_k(m/\ell^2)},
	\ee
	for which we defined the auxiliary functions
	\ba\label{eq;gkauxfunct}
	\gf{f}{m} &=\sum_{y^2\mid m}\mu(y)c_{m/y^2}(f),\\
	g_k(m) &=\sum_{y^2\mid m}\mu(y)\sigma_{k-1}(m/y^2)=m^{k-1}\prod_{p\mid m}(1+p^{-k+1}).
	\ea
	\end{prop}
	We conclude this section with a bound for the cuspidal part $c_{n,r}(\phi^0_m)$.
	\begin{prop}[See {\cite[Corollary 6.5]{boda;petnorm}}]\label{prop;Ocusppart}
	For all sequences $\seqbig{T_j=\big(\begin{smallmatrix}
	n_j & r_j/2\\
	r_j/2 & m_j
	\end{smallmatrix}\big)}{j\in\NN}$ of \emph{reduced} matrices in~$\halfint^+_2$ of increasing determinant, the cuspidal part of~$a^k_2(f,T_j)$ appearing in~\eqref{eq;decfourinjac} satisfies the bound
	\bes
	c_{n_j,r_j}(\phi^0_{m_j})=O\big(\det(T_j)^{k/2+1/4+\varepsilon}\big),
	\ees
	for every $\varepsilon>0$.
	\end{prop}
	
	\section{Background on cones}\label{sec;backgrcones}
	
	In this section we introduce the cones of special cycles of codimension two on orthogonal Shimura varieties associated to \emph{unimodular} lattices, and the cone of coefficient extraction functionals of Siegel modular forms.
	Eventually, we explain how to deduce geometric properties of the former via the ones of the latter.
	To fix the notation, we briefly recall the needed background on cones.\\
	
	Let $V$ be a non-trivial finite-dimensional vector space over $\QQ$, and let~$\gencone$ be a non-empty subset of~$V$.
	The \emph{(convex) cone generated by~$\gencone$} is the smallest subset of~$V$ that contains~$\gencone$ and is closed under linear combinations with non-negative coefficients.
	We denote it either by~$C_\QQ(\gencone)$, or by~$\langle \gencone\rangle_{\QQ_{\ge 0}}$.
	If there exists a finite subset~$\gencone'\subseteq \gencone$ such that~$C_\QQ(\gencone')=C_\QQ(\gencone)$, we say that the cone~$C_\QQ(\gencone)$ is \emph{polyhedral} (or \emph{finitely generated}).
	A cone is said to be \emph{pointed} if it contains no lines.
	
	The \emph{convex hull of~$\gencone$} is the smallest convex subset of~$V$ containing~$\gencone$. It is denoted by~$\conv_{\QQ}(\gencone)$, and coincides with the set of all convex combinations of elements of~$\gencone$, namely
	\bes
	\conv_{\QQ}(\gencone)=\bigg\{\sum_{g\in J} x_g \cdot g\,:\,\text{$J\subseteq\gencone$ is finite, $\sum_{g\in J} x_g=1$ and $x_g\in\QQ_{\ge0}$}\bigg\}.
	\ees
	Analogous definitions holds over $\RR$.
	
	For simplicity, from now on we suppose that $C_\QQ(\gencone)$ is full-dimensional in~$V$.
	The~\emph{$\RR$-closure of $C_\QQ(\gencone)$} is the topological closure~$\overline{C_\QQ(\gencone)}=\overline{C_\QQ(\gencone)\otimes_\QQ\RR}$ of~$C_\QQ(\gencone)$ in the vector space~$V\otimes\RR$ endowed with the Euclidean topology.
	The \emph{boundary rays of~$C_\QQ(\gencone)$} are the rays of $\overline{C_\QQ(\gencone)}$ lying on its boundary.
	An \emph{extremal ray} of~$C_\QQ(\gencone)$ is a boundary ray of $C_\QQ(\gencone)$ that does not lie in the interior of any subcone of~$\overline{C_\QQ(\gencone)}$ of dimension higher than one.
	We say that~$C_\QQ(\gencone)$ is a \emph{rational cone} if all its extremal rays can be generated by vectors of~$V$.
	
	A ray $r$ of~$V\otimes\RR$ is said to be an \emph{accumulation ray of~$C_\QQ(\gencone)$ with respect to the set of generators~$\gencone$} if there exists a sequence of pairwise different generators $\seq{g_j}{j\in\NN}$ in~$\gencone$ such that~$\RR_{\ge0}\cdot g_j \to r$ when~$j\to\infty$.	
	Clearly, all accumulation rays lie in~$\overline{C_\QQ(\gencone)}$.
	The \emph{accumulation cone of $C_\QQ(\gencone)$ with respect to~$\gencone$} is defined as the subcone of~$\overline{C_\QQ(\gencone)}$ generated by the accumulation rays of~$C_\QQ(\gencone)$ with respect to~$\gencone$.
	If there is no accumulation ray, it is defined as the trivial cone~$\{0\}$.
	
	Clearly, all previous definitions extend also to cones defined on real vector spaces.
	\begin{ex}
	Consider the subset of~$\QQ^2$ defined as
	\bas
	\gencone_1&=\{(1,a)\,:\,a\in[0,1]\cap\QQ\},\\
	\gencone_2&=\{(1,a)\,:\,a\in[0,1)\cap\QQ\},\\
	\gencone_3&=\{(1,a)\,:\,a\in[0,\pi)\cap\QQ\}.
	\eas
	The cone $C_\QQ(\gencone_1)$ is rational and polyhedral, with extremal rays $\RR_{\ge0}\cdot(1,0)$ and $\RR_{\ge0}\cdot (1,1)$.
	The cone~$C_\QQ(\gencone_2)$ is rational but non-polyhedral, and its $\RR$-closure~$\overline{C_\QQ(\gencone_2)}$ is rational and polyhedral.
	The cone $C_\QQ(\gencone_3)$ is neither rational nor polyhedral, while its $\RR$-closure is polyhedral but non-rational.
	\end{ex}
	In Section \ref{sec;propC}, we will study how infinitely many extremal rays of a cone could converge towards an accumulation cone.
	Along this process, it is important to keep in mind that even if a sequence of extremal rays converges, the boundary ray obtained as a limit does not have to be extremal.
	This behavior, which can happen only in dimension higher than~$3$, is illustrated in the following example.
	\begin{ex}\label{ex;acextrnotextr}
	Let $c$ be the semicircle in $\RR^3$ defined as $c=\{(\cos\theta,\sin\theta,0)\,:\,\theta\in[0,\pi]\}$, and let $A=(1,0,1)$ and $B=(1,0,-1)$.
	We define the convex hull $\gencone=\conv_\RR(c\cup\{A,B\})$ and the inclusion $\iota\colon\RR^3\to\RR^4$ by $(x,y,z)\mapsto(x,y,z,1)$.
	The cone $C_\RR(\iota(\gencone))$ is a full-dimensional pointed cone in $\RR^4$, with extremal rays
	\bes
	\text{$\RR_{\ge 0}\cdot\iota(A)$,$\quad\RR_{\ge 0}\cdot\iota(B)\quad$and$\quad\RR_{\ge 0}\cdot\iota(P)$ for all $P\in c-\{(1,0,0)\}$.}
	\ees
	In fact, every vector of the boundary ray $\RR_{\ge 0}\cdot\iota(1,0,0)$ is a non-negative combination of some vectors lying on the two extremal rays given by $\iota(A)$ and $\iota(B)$.
	Let $\seq{\theta_j}{j\in\NN}$ be a sequence of pairwise different elements in the interval $(0,\pi)$ converging to $0$.
	The sequence of extremal rays $\RR_{\ge 0}\cdot\iota(\cos\theta_j,\sin\theta_j,0)$ converges to the boundary ray $\RR_{\ge 0}\cdot\iota(1,0,0)$, which is non-extremal.
	\end{ex}
	
	Let~$\psi\colon V\to W$ be a linear map of~$\QQ$-vector spaces of finite dimensions.
	If a cone~$\cone\subset V$ is rational, resp.\ polyhedral, then also the cone~$\psi(\cone)\subset W$ is rational, resp.\ polyhedral.
	Nevertheless, there are properties of~$\cone$ that may not be preserved by~$\psi$.
	In fact, as shown in the following example, there are linear maps~$\psi$ mapping a pointed cone~$\cone$ to a cone that contains a line, and mapping the accumulation cone of~$\cone$ with respect to a set of generators~$\gencone$, to a cone that is not the accumulation cone of~$\psi(\cone)$ with respect to the set of generators~$\psi(\gencone)$.
	
	\begin{ex}
	Let~$P_j=(1,1/t,0)\in\RR^3$, for every positive integer~$j$, and let
	\bes
	A=(1,0,0),\qquad B=(0,1,0),\qquad C=(0,0,1).
	\ees
	We define~$\cone$ as the cone in~$\RR^3$ generated by the set
	\bes
	\gencone=\{A,B,C\}\cup\{P_t : \text{$t\in\ZZ_{>0}$}\}.
	\ees
	Let~$\pi\colon \RR^3\to\RR^2$, $(x,y,z)\mapsto(x,z)$, and let~$\Pi\colon\RR^2\to\RR$ be the projection to the line generated by the vector~$(1,-1)$ in~$\RR^2$.
	We define the linear map~$\psi\colon\RR^3\to\RR$ as the composition~$\Pi\circ\pi$.
	The cone~$\cone$ is pointed, but its image~$\psi(\cone)$ is not.
	The accumulation cone of~$\cone$ with respect to~$\gencone$ is given by the ray~$\RR_{\ge0}\cdot A$, which maps to a non-trivial ray \textcolor{\myblack}{under}~$\psi$.
	However, since the set~$\psi(\gencone)$ is finite, the accumulation cone of~$\psi(\cone)$ with respect to~$\psi(\gencone)$ is trivial.
	\end{ex}
	The following result provides a sufficient condition for the contraction of an accumulation ray via a linear map.
	It will be used in Section~\ref{sec:concoefextrfun} to show that many of the properties of the cone of special cycles are inherited from the cones of coefficients of Siegel modular forms; see Corollary~\ref{cor;prespoint&acc}.
	\begin{lemma}\label{lem;tecnondeduceaccones}
	Let~$\psi\colon V \to W$ be a linear map of Euclidean vector spaces of finite dimensions, and let~$\seq{v_j}{j\in\NN}$ be a sequence of pairwise different vectors in~$V$.
	Suppose that~$v_j=r_j e +\widetilde{v_j}$, for some~$r_j\in\RR_{>0}$ and~$e,\widetilde{v_j}\in V$, such that
	\bes
	e\perp\widetilde{v_j},\qquad r_j\rightarrow\infty,\qquad\text{and}\qquad \frac{\widetilde{v_j}}{r_j}\rightarrow 0.
	\ees
	If~$\seqbig{\psi(v_j)}{j\in\NN}$ is a constant sequence in~$W$, then~$\psi(e)=0$.
	In particular, the accumulation ray~$\RR_{\ge0}\cdot e$ arising from the sequence of vectors~$\seq{v_j}{j}$ is contracted by~$\psi$.
	\end{lemma}
	\begin{proof}
	Since~$\psi(v_j-v_0)=0$ for every~$j\in\NN$, we may divide both terms of such equality by~$r_j$, and deduce that
	\be\label{eq;proofoftecnondeduceaccones}
	0=\psi\Big(\frac{v_j-v_0}{r_j}\Big)=\psi\Big(e + \frac{\widetilde{v_j}}{r_j} + \frac{v_0}{r_j}\Big).
	\ee
	Since both~$\widetilde{v_j}/r_j$ and~$v_0/r_j$ tends to zero when~$j\to\infty$ by hypothesis, the right-hand side of~\eqref{eq;proofoftecnondeduceaccones} tends to~$\psi(e)$, hence~$\psi(e)=0$.
	\end{proof}
	\subsection{Cones of special cycles of codimension~$\boldsymbol{2}$}\label{sec;conesspeccycles}
	
	In this section, we define the cones of special cycles associated to orthogonal Shimura varieties.
	We restrict the illustration to cycles of codimension two on varieties associated to \emph{unimodular} lattices.
	The relationship with Siegel modular forms is provided in Section \ref{sec:concoefextrfun}.\\
	
	Let $X$ be a normal irreducible complex space of dimension~$b$.
	A \emph{cycle $Z$ of codimension~$g$} in $X$ is a locally finite formal linear combination
	\bes
	Z=\sum n_Y Y,\qquad\text{$n_Y\in\ZZ$,}
	\ees
	of distinct closed irreducible analytic subsets~$Y$ of codimension~$g$ in~$X$.
	The \emph{support} of the cycle~$Z$ is the closed analytic subset $\supp(Z)=\bigcup_{n_Y\neq 0} Y$ of pure codimension~$g$.
	The integer~$n_Y$ is the \emph{multiplicity} of the irreducible component~$Y$ of~$\supp(Z)$ in the cycle~$Z$.
	
	If~$X$ is a manifold, and~$\Gamma$ is a group of biholomorphic transformations of~$X$ acting properly discontinuously, we may consider the preimage~$\pi^*(Z)$ of a cycle~$Z$ of codimension~$g$ on~$X/\Gamma$ under the canonical projection $\pi\colon X\to X/\Gamma$.
	For any irreducible component~$Y$ of~$\pi^{-1}\big(\supp(Z)\big)$, the multiplicity~$n_Y$ of~$Y$ with respect to~$\pi^*(Z)$ equals the multiplicity of~$\pi(Y)$ with respect to~$Z$.
	This implies that~$\pi^*(Z)$ is a \emph{$\Gamma$-invariant cycle}, meaning that if~$\pi^*(Z)=\sum n_Y Y$, then
	\bes
	\gamma\big(\pi^*(Z)\big)\coloneqq \sum n_Y \gamma(Y)\qquad\text{equals $\pi^*(Z)$, for every $\gamma\in \Gamma$.}
	\ees
	Note that we do not take account of possible ramifications of the cover $\pi$.
	
	We now focus on orthogonal Shimura varieties associated to \emph{unimodular lattices}.
	Let $L$ be an even non-degenerate unimodular lattice of signature $(b,2)$.
	We denote by $(\cdot{,}\cdot)$ the bilinear form of $L$, and by $q$ the quadratic form defined as $q(\lambda)=(\lambda,\lambda)/2$, for every $\lambda\in L$.
	The $b$-dimensional complex manifold
	\bes
\mathcal{D}_b=\{ z\in L\otimes\CC\setminus\{0\}\, :\,\text{$( z,z)=0$ and $( z,\bar{z})<0$}\}/\CC^* \subset \mathbb{P}(L\otimes\CC)
	\ees
	has two connected components.
	The action of the group of the isometries of~$L$, denoted by~$\bigO(L)$, extends to an action on~$\mathcal{D}_b$.
	We choose a connected component of~$\mathcal{D}_b$ and denote it by~$\projmod$.
	We define~$\bigO^+(L)$ as the subgroup of~$\bigO(L)$ containing all isometries which preserve~$\projmod$.
	
	Let~$\Gamma$ a subgroup of finite index in~$\bigO^+(L)$.
	The \emph{orthogonal Shimura variety} associated to~$\Gamma$ is
	\bes
	X_\Gamma=\Gamma\backslash\projmod.
	\ees
	By the Theorem of Baily and Borel, the analytic space $X_\Gamma$ admits a unique algebraic structure, which makes it a \emph{quasi-projective algebraic variety}. Each of these varieties inherits a line bundle from the restriction of the tautological line bundle $\mathcal{O}(-1)$ on $\PP(L\otimes\CC)$ to $\projmod$. This is the so-called \emph{Hodge bundle}, which we denote by $\omega$.
	
	An attractive feature of this kind of varieties is that they have many algebraic cycles.
	We recall here the construction of the so-called \emph{special cycles}; see~\cite{ku;algcycle} for further information.
	They are a generalization of the Heegner divisors in higher codimension; see~\cite[Section~$5$]{br;borchp} for a description of such divisors in a setting analogous to the one of this paper.
	
	Recall that $\halfint_2$, resp.~$\halfint_2^+$, is the set of symmetric half-integral positive semi-definite, resp.\ positive definite, $2\times 2$-matrices.
	If $\vect{\lambda}=(\lambda_1,\lambda_2)\in L^2$, the \emph{moment matrix} of~$\vect{\lambda}$ is defined as~${q(\vect{\lambda})\coloneqq\frac{1}{2}\big( (\lambda_i,\lambda_j)\big)_{i,j}}$, where~$\big( (\lambda_i,\lambda_j)\big)_{i,j}$ is the matrix given by the inner products of the entries of~$\vect{\lambda}$, while its orthogonal complement in~$\projmod$ is~$\vect{\lambda}^\perp = \lambda_1^\perp\cap\lambda_2^\perp$.
	For every~$T\in\halfint^+_2$, the codimension~$2$ cycle
	\ba\label{eq;defSpeccydimmcones}
	\sum_{\substack{\vect{\lambda}\in L^2\\ q(\vect{\lambda})=T}} \vect{\lambda}^\perp
	\ea
	is $\Gamma$-invariant in $\projmod$.
	Since the componentwise action of~$\Gamma$ on the vectors $\vect{\lambda}\in L^2$ of fixed moment matrix $T\in\halfint^+_2$ has finitely many orbits, the cycle~\eqref{eq;defSpeccydimmcones} descends to a cycle of codimension $2$ on~$X_\Gamma$, which we denote by $Z(T)$ and call \emph{the special cycle associated to~$T$}.
	\textcolor{\myblack}{Special cycles} are preserved \textcolor{\myblack}{by} pullbacks of quotient maps~$\pi\colon X_{\Gamma'}\to X_\Gamma$, for every subgroup~$\Gamma'$ of finite index in~$\Gamma$. This is the reason why we usually drop~$\Gamma$ from the notation, writing only~$Z(T)$ instead of~$Z(T)_\Gamma$.
	
	\begin{rem}\label{rem;Heegnerspeccy}
	An analogous construction works for matrices~$T\in\halfint_2$, where the associated special cycles have codimension~$\rank(T)$.
	The divisors $Z\left(\begin{smallmatrix}
	n & 0\\ 0 & 0
	\end{smallmatrix}\right)$ of $X_\Gamma$, where $n$ is a positive integer, are the so-called \emph{Heegner divisors}, usually denoted by $H_n$.
	These, together with~$\omega^*$, are the special cycles of codimension one of $X_\Gamma$, and their classes generate the whole~$\Pic(X_\Gamma)$, as proved in~\cite[Corollary~$3.8$]{blmm;conj}.
	\end{rem}
	
	If~$Z$ is a cycle of codimension~$r$ in~$X_\Gamma$, we denote by~$\{Z\}$ its rational class in the Chow group~$\CH^r(X_\Gamma)$, and by~$[Z]$ its cohomology class in~$H^{r,r}(X_\Gamma)$.
	By Poincaré duality, we may consider~$[Z]$ as a linear functional on the cohomology space of compactly supported closed~$(r,r)$-forms on~$X_\Gamma$; see e.g.~\cite[Section~$8.1$]{blmm;conj}.
	
	Eventually, we define the cones of special cycles we treat in this paper.
	\begin{defi}\label{defi;conesspeccycl}
	Let $X_\Gamma$ be an orthogonal Shimura variety associated to a non-degenerate even unimodular lattice of signature $(b,2)$, with $b>2$.
	The \emph{cone of special cycles (of codimension~$2$) on $X_\Gamma$} is the cone in $\CH^2(X_\Gamma)\otimes\QQ$ defined as
	\bes
	\conecy=\langle \{Z(T)\}\,:\, T\in\halfint_2^+\rangle_{\QQ_{\ge 0}},
	\ees
	while the \emph{cone of rank one special cycles (of codimension~$2$) on $X_\Gamma$} is
	\bes
	\conecyone=\langle \{Z(T)\}\cdot\{\omega^*\}\,:\, \text{$\T\in\halfint_2$ and $\rank(T)=1$} \rangle_{\QQ_{\ge 0}}.
	\ees
	\end{defi}
	Whenever we refer to the~\emph{accumulation cones} of~$\conecy$ and~$\conecyone$, we implicitly consider them with respect to the set of generators of~$\conecy$ and~$\conecyone$ used in Definition~\ref{defi;conesspeccycl}.
	
	Although it is still unclear whether $\CH^2(X_\Gamma)\otimes\QQ$ is finite-dimensional, it is known that the span over~$\QQ$ of the special cycles of codimension two is of finite dimension; see \cite[Corollary 6.3]{brra;modconj}.
	In particular, both $\conecy$ and $\conecyone$ are of finite dimension.
	
	The cone $\cone_{X_\Gamma}'$ is \emph{pointed, rational, and polyhedral}. We provide a proof based on the main result of \cite{brmo} at the end of this section.
	In Section \ref{sec:concoefextrfun}, we will explain how to deduce these properties using the Fourier coefficients of Siegel modular forms associated to singular matrices.
	The main property of $\conecyone$ is that it has only one accumulation ray, which is generated by an internal point of the $\RR$-closure of~$\conecyone$.
	
	The geometry of $\conecy$ is more interesting, although more complicated.
	We prove in Section~\ref{sec:concoefextrfun} that the accumulation cone of $\conecy$ is \emph{pointed, rational, and polyhedral}, deducing the rationality of $\conecy$.
	The explicit classification of all accumulation rays of~$\conecy$ is provided in Section~\ref{sec;acrayschow}.
	
	The rational class $\{\omega^*\}^2\in\CH^2(X_\Gamma)$ does not appear neither in the set of generators of~$\conecy$ nor in the one of~$\conecyone$.
	It will be clear at the end of Section \ref{sec;accrays} that it is contained in the interior of~$\conecy$.
	The properties of $\conecy$ and $\conecyone$ stated above are summarized in the following result.
	\begin{thm}\label{thm;nammainintro}
	Let $X$ be an orthogonal Shimura variety associated to a non-degenerate even \emph{unimodular} lattice of signature $(b,2)$, with $b>2$.
	\begin{enumerate}[label=(\roman{*})]
	\item The cone of rank one special cycles ${\cone'}_{\!\! X}$ is pointed, rational, polyhedral, and of dimension~$\dim {{M_1}^{\!1+b/2}}$. (Bruinier--Möller)\label{nammainintrofirstp}
	\item The accumulation cone of the cone of special cycles $\cone_{X}$ is pointed, rational, polyhedral, and of the same dimension as ${\cone'}_{\!\! X}$.
	\item The cone $\cone_{X}$ is rational and of maximal dimension in the subspace of~$\CH^2(X)\otimes\QQ$ generated by the special cycles of codimension~$2$.
	\item The cones~$\cone_{X}$ and~${\cone'}_{\!\! X}$ intersect only at the origin.
	Moreover, if the accumulation cone of~$\cone_{X}$ is enlarged with a non-zero element of~${\cone'}_{\!\! X}$, the resulting cone is non-pointed.
	\textcolor{\myblack}{In particular~$-{\cone'}_{\!\! X}\subseteq\overline{\mathcal{C}_X}$.}
	\end{enumerate}
	\end{thm}
	
	As we explain in Section~\ref{sec:concoefextrfun}, the properties of the cones of special cycles appearing in Theorem~\ref{thm;nammainintro} are strictly connected with the analogous properties of certain cones of coefficient extraction functionals of Siegel modular forms.
	While working with rational classes of cycles on a variety is notoriously hard, the coefficient extraction functionals of Siegel modular forms can be computed \emph{explicitly} over a basis of~$M^k_2$.
	In this paper, we use the arithmetic properties of such functionals to prove Theorem~\ref{thm;nammainintro}.
	We will see also how the polyhedrality problem of~$\cone_{X_\Gamma}$ can be studied with Siegel modular forms via Conjecture~\ref{conj;our}.
	
	We conclude this section with the proof of the first part of our main result.
	\begin{rem}\label{rem;onHLforXgammareferee}
	\textcolor{\myblack}{The map given by wedging elements of~$H^r(X_\Gamma,\CC)$ with the Hodge class~$\omega$ is injective for all~$r<b-2$, although the variety~$X_\Gamma$ is only quasi-projective.
	This can be deduced in terms of the analogue of the Hard Lefschetz Theorem~\cite[Corollary~$9.2.3$]{interscoho} for the intersection cohomology of the Baily--Borel compactification~$\bbcomp$ of~$X_\Gamma$.
	In fact, the intersection cohomology group~$IH^r(\bbcomp,\CC)$ is isomorphic to~$H^r(X_\Gamma,\CC)$ for every~$r<b-1$, as proved in~\cite{looijgengaL2} \cite{SapSt;L2coho}, and the Kähler class~$\omega$ of~$X_\Gamma$ is identified with the Chern class of an ample line bundle in~$\bbcomp$; see~\cite[Sections~$2.4$ and~$2.5$]{blmm;conj} for further information.}
	\end{rem}
	\begin{proof}[Proof of Theorem \ref{thm;nammainintro} \textit{\ref{nammainintrofirstp}}.]
	The cone
	\bes
	\widetilde{\cone}=\langle \{Z(T)\}\,:\,\text{$T\in\halfint_2$ and $\rank(T)=1$}\rangle_{\QQ_{\ge 0}}
	\ees
	is the cone in $\Pic(X_\Gamma)\otimes\QQ$ generated by the Heegner divisors $\{H_n\}$.
	In fact, since
	\bes
	Z(T)=Z(u^t\cdot T\cdot u)\quad\text{for every $u\in\GL_2(\ZZ)$},
	\ees
	we deduce that for every~$T\in\halfint_2$ of rank one there exists a positive integer~$n$ such that~$\{Z(T)\}$ is equal to~$\left\{Z\left(\begin{smallmatrix}
	n & 0 \\ 0 & 0
	\end{smallmatrix}\right)\right\}$.
	The latter is the Heegner divisor~$H_n$; see Remark~\ref{rem;Heegnerspeccy}.
	The intersection map
	\bes
	\rho\colon \Pic(X_\Gamma)\otimes\QQ\longrightarrow \CH^2(X_\Gamma)\otimes\QQ,\quad \{H_n\}\longmapsto \{H_n\}\cdot\{\omega^*\}=\{Z\left(\begin{smallmatrix}
	n & 0\\
	0 & 0
	\end{smallmatrix}\right)\}\cdot\{\omega^*\},
	\ees
	is linear and maps $\widetilde{\cone}$ to $\conecyone$.
	By \cite[Theorem 3.4]{brmo} the former cone is rational, polyhedral, and of dimension $\dim M^k_1$.
	Since $\rho$ is linear, also $\cone_{X_\Gamma}'$ is rational, polyhedral, and of dimension at most~$\dim M^k_1$.
	
	We conclude the proof showing that the dimension of~$\widetilde{\cone}$ and its pointedness are preserved \textcolor{\myblack}{by}~$\rho$.
	To do so, it is enough to show that~$\rho$ is injective.
	Consider the commutative diagram
	\begin{center}
	\begin{tikzcd}
		\Pic(X_\Gamma)\otimes\QQ \arrow[r, "\rho"] \arrow[d, "cl_1"]
		& \CH^2(X_\Gamma)\otimes\QQ \arrow[d, "cl_2"] \\
		H^{1,1}(X_\Gamma,\QQ) \arrow[r, "\sigma"]
		& H^{2,2}(X_\Gamma,\QQ)
	\end{tikzcd}
	\end{center}
	where the vertical arrows are the cycle maps, and~$\sigma$ is the map induced by the exterior product with~$-\omega$.
	By~\cite[Corollary~$3.8$]{blmm;conj}, the map~$cl_1$ is an isomorphism, and by \textcolor{\myblack}{Remark~\ref{rem;onHLforXgammareferee}}, the map~$\sigma$ is injective, hence~$\rho$ is injective as well.
	\end{proof}
	\subsection{Cones of coefficient extraction functionals and modularity}\label{sec:concoefextrfun}
	Let $k>4$ be an even integer.
	Recall that we denote by $M^k_2(\QQ)$ the space of weight~$k$ Siegel modular forms of genus~$2$ with rational Fourier coefficients, and by $c_T$ the coefficient extraction functional associated to a matrix $T\in\halfint_2$; see Section \ref{sec;siegmodforms} for further information.
	\begin{defi}\label{def;modcones}
	The \emph{modular cone of weight $k$} is the cone in the dual space $M^k_2(\QQ)^*$ defined as
	\bes
	\conemod=\langle c_T\,:\, T\in\halfint^+_2\rangle_{\QQ_{\ge0}},
	\ees
	while the \emph{rank one modular cone of weight $k$} is
	\bes
	\conemodone=\langle c_T\,:\,\text{$T\in\halfint_2$ and $\rank(T)=1$}\rangle_{\QQ_{\ge0}}.
	\ees
	\end{defi}
	Whenever we refer to the~\emph{accumulation cones} of~$\conemod$ and~$\conemodone$, we implicitly consider the ones with respect to the set of generators of~$\conemod$ and~$\conemodone$ appearing in Definition~\ref{def;modcones}.
	
	The following proposition is the key result to relate the cones of functionals with the cones of special cycles; see also \cite[Corollary~1.8]{ra;forfou}.
	\begin{prop}\label{prop;functpsi}
	Let $X_\Gamma$ be an orthogonal Shimura variety associated to an even unimodular lattice of signature~$(b,2)$, with~${b>2}$. The map
	\bes
	\psi_\Gamma\colon M^{1+b/2}_2(\QQ)^*\longrightarrow \CH^2(X_\Gamma)\otimes\QQ,\quad c_T\longmapsto\{Z(T)\}\cdot\{\omega^*\}^{2-\rank(T)},
	\ees
	is \emph{well-defined} and \emph{linear}.
	\end{prop}
	\begin{proof}
	The function over $\HH_2$ defined as
	\bes
	\Theta_\Gamma(Z)=\sum_{T\in\halfint_2}\{Z(T)\}\cdot\{\omega^*\}^{2-\rank(T)}e^{2\pi i \trace(TZ)}
	\ees
	is a Siegel modular form of weight $1+b/2$ with values in $\CH^2(X_\Gamma)\otimes\CC$.
	This follows from the so-called Kudla's Modularity Conjecture, proved for the case of genus~$2$ in~\textcolor{\myblack}{\cite{br;vvformal} and}~\cite{ra;forfou}, and for general genus in~\cite{brra;modconj}.
	The previous compact formulation \textcolor{\myblack}{implies} the following one; see \cite[Corollary 6.2]{ra;forfou}.
	For every linear functional~${f\in (\CH^2(X_\Gamma)\otimes\CC)^*}$, the formal Fourier expansion
	\bes
	\Theta_{\Gamma,f}(Z):=\sum_{T\in\halfint_2}f\big(\{Z(T)\}\cdot\{\omega^*\}^{2-\rank(T)}\big)e^{2\pi i \trace(TZ)}
	\ees
	is a Siegel modular form of weight $1+b/2$.
	
	Let $\{T_j\}_{j=1}^s$ be a finite set of matrices in $\halfint_2$.
	Suppose that there exist complex numbers~$\lambda_j$ such that~$\sum_j \lambda_j c_{T_j}=0$ in $\big(M^{1+b/2}_2\big)^*$, or equivalently that~${\sum_j\lambda_j c_{T_j}(F)=0}$, for every~${F\in M^{1+b/2}_2}$.
	We deduce that
	\bes
	\sum_{j=1}^s\lambda_j c_{T_j}(\Theta_{\Gamma,f})=
	\sum_{j=1}^s\lambda_j f\Big(\{Z(T_j)\}\cdot\{\omega^*\}^{2-\rank(T_j)}\Big)
	=f\Big(\sum_{j=1}^s\lambda_j \{Z(T_j)\}\cdot\{\omega^*\}^{2-\rank(T_j)}\Big)
	=0,
	\ees
	for every functional $f$.
	This implies that the complex extension of~$\psi_\Gamma$ is a homomorphism. Since the complex space $M^{1+b/2}_2$ admits a basis of Siegel modular forms with rational Fourier coefficients, the restriction $\psi_\Gamma$ over $\QQ$ is well-defined.
	\end{proof}
	\begin{thm}\label{thm;mainthm2}
	Let $k>4$ be an even integer such that $k\equiv 2$ mod $4$.
	\begin{enumerate}[label=(\roman{*})]
	\item The rank one modular cone $\conemodone$ is pointed, rational, polyhedral, and of the same dimension as $M^k_1$.\label{pt;conek'}
	\item The accumulation cone of the modular cone $\conemod$ is pointed, rational, polyhedral, and of the same dimension as $M^k_1$.\label{pt2;conek'}
	\item The cone $\conemod$ is pointed, rational, and of the same dimension as $M^k_2$.\label{pt3;conek'}
	\item The cones $\conemod$ and $\conemodone$ intersect only at the origin.\label{pt4;conek'}
	Moreover, if the cone $\conemod$ is enlarged with a non-zero element of \textcolor{\myblack}{$\conemodone$}, the resulting cone is non-pointed.
	\textcolor{\myblack}{In particular~$-\conemodone\subseteq\overline{\conemod}$}
	\end{enumerate}
	\end{thm}
	Let~$\psi_\Gamma$ be as in Proposition~\ref{prop;functpsi}.
	The images \textcolor{\myblack}{under}~$\psi_\Gamma$ of the cones~$\cone_{1+b/2}$ and~${\cone'}_{\!\!1+b/2}$ are~$\conecy$ and~$\conecyone$, respectively.
	Note that since the variety $X_\Gamma$ is associated to a \emph{unimodular} lattice of signature~$(b,2)$, the weight $1+b/2$ is an even integer congruent to $2$ mod $4$.
	
	By means of Lemma~\ref{lem;tecnondeduceaccones}, we may deduce the following non-trivial properties of~$\conecy$ via the ones of~$\conemod$.
	\begin{cor}\label{cor;prespoint&acc}
	The accumulation cone of~$\conecy$ is pointed, and every accumulation ray of~$\conemod$ maps \textcolor{\myblack}{under}~$\psi_\Gamma$ to an accumulation ray of~$\conecy$.
	Moreover, the accumulation cones of~$\conemod$ and~$\conecy$ have the same dimension.
	\end{cor}
	\begin{proof}
	Let~$\iota\colon M^k_1(\RR)^*\to M^k_2(\RR)^*$, be the embedding defined as~$c_n\mapsto c_{\big(\begin{smallmatrix}n & 0\\ 0 & 0\end{smallmatrix}\big)}$, for every~${n\in\NN}$.
	Consider the commutative diagram
	\begin{center}
		\begin{tikzcd}
		M^k_1(\RR)^* \arrow[r, "\iota"] \arrow[d, "\psi'_\Gamma"]
		& M^k_2(\RR)^* \arrow[d, "\psi_\Gamma"] \\
		\Pic(X_\Gamma)\otimes\RR \arrow[r, "\rho"]
		& \CH^{\textcolor{\myblack}{2}}(X_\Gamma)\otimes\RR
		\end{tikzcd}
	\end{center}
	where~$\rho$ is the map given by the intersection with~$\{\omega^*\}$, and~$\psi'_\Gamma$ is the \textcolor{\myblack}{analogue} of~$\psi_\Gamma$ for Heegner divisors in~$X_\Gamma$, namely it maps~$c_n\mapsto\{H_n\}$, for every positive integer~$n$, and~$c_0\mapsto\{\omega^*\}$.
	As explained in~\cite[Section~$4$]{brmo}, the map~$\psi'_\Gamma$ is an isomorphism.
	In fact, also the composition of~$\psi'_\Gamma$ with the cycle map is so.
	Since the map induced by~$\rho$ in cohomology is injective by the Hard Lefschetz Theorem, we deduce that~$\rho$ is injective as well\textcolor{\myblack}{; see Remark~\ref{rem;onHLforXgammareferee}}.
	
	We will see in Section~\ref{sec;accumulcone} that the accumulation cone~$\accone_k$ of~$\conemod$ is pointed and contained \textcolor{\myblack}{in} the image of the embedding~$\iota$.
	Since the diagram above is commutative, and~$\rho\circ\psi'_\Gamma$ is injective, we deduce that~$\psi_\Gamma$ embeds~$\accone_k$ into~$\CH^2(X_\Gamma)\otimes\RR$, therefore~$\psi_\Gamma(\accone_k)$ is pointed and of dimension~$\dim M^k_1$.
	
	We conclude the proof by showing that every accumulation ray of~$\conemod$ maps to an accumulation ray of~$\conecy$ \textcolor{\myblack}{under}~$\psi_\Gamma$.
	Suppose this is not the case, namely there exists an accumulation ray~$r$ of~$\conemod$ such that~$\psi_\Gamma(r)$ is not an accumulation ray of~$\cone_{X_\Gamma}$.
	This means that there exists a sequence~$\seq{T_j}{j\in\NN}$ of reduced matrices in~$\halfint^+_2$ of increasing determinant, such that the functionals~$c_{T_j}$ are pairwise different and~$\RR_{\ge0}\cdot c_{T_j}\to r$, but such that the sequence of cycles~$\seqbig{\{Z(T_j)\}}{j\in\NN}$ is constant.
	Let~$e$ be a generator of the accumulation ray~$r$.
	We decompose~$c_{T_j}=r_j e+\widetilde{v_j}$, for some~$r_j\in\RR$ and some~$\widetilde{v_j}$ orthogonal to~$e$.
	Since~$\RR_{\ge0}\cdot c_{T_j}\to\RR_{\ge0}\cdot e$, we deduce that~$r_j$ is eventually positive, and that~$\widetilde{v_j}/r_j\to 0$ when~$j\to+\infty$.
	Moreover, since~$c_{T_j}(E^k_2)\to\infty$ by Lemma~\ref{classsiegeis}, we deduce that also~$r_j$ diverges.
	By Lemma~\ref{lem;tecnondeduceaccones}, the map~$\psi_\Gamma$ contracts the ray~$r$.
	But this is not possible, since~$\psi_\Gamma$ is injective on~$\accone_k$, as proved at the beginning of this proof.
	\end{proof}
	\begin{rem}\label{rem;nonpointandinjpsigamma}
	The problem of the pointedness of the whole cone~$\conecy$ is more subtle.
	As shown in Theorem~\ref{thm;mainthm2}, the modular cone~$\conemod$ is pointed.
	However, the map~$\psi_\Gamma$ might contract some of the rays of~$\conemod$, making~$\conecy$ non-pointed.
	This is not the case if e.g.~$\psi_\Gamma$ is injective.
	Such injectivity is a non-trivial open problem.
	It seems reasonable it may be tackled proving the injectivity of the Kudla--Millson lift of genus~$2$, as explained in~\cite{br;borchp} for the counterpart of~$\psi_\Gamma$ for elliptic modular forms.
	We will return to such interesting problem in a future work.
	\end{rem}
	Since the rationality and the polyhedrality are geometric properties of cones which are preserved by linear maps between vector spaces over $\QQ$, Theorem~\ref{thm;nammainintro} follows from Theorem~\ref{thm;mainthm2} and Corollary~\ref{cor;prespoint&acc}.
	We remark that the polyhedrality of the cone of special cycles~$\cone_{X_\Gamma}$ is implied by Conjecture~\ref{conj;our}.
	
	We conclude this section with the proof of Theorem~\ref{thm;mainthm2}~\ref{pt;conek'}.
	The remaining points of Theorem~\ref{thm;mainthm2} are proven in the following sections.
	\begin{proof}[Proof of Theorem~\ref{thm;mainthm2}~\textit{\ref{pt;conek'}}]
	If the matrix $T\in\halfint_2$ has rank one, then there exists~${u\in \GL_2(\ZZ)}$ such that ${u^t\cdot T\cdot u=\left(\begin{smallmatrix}
	n & 0\\
	0 & 0
	\end{smallmatrix}\right)}$, for some $n\in\NN$.
	We denote the latter matrix by $M(n)$, \textcolor{\myblack}{for} simplicity.
	Since $c_T=c_{u^t T u}$ for every $u\in \GL_2(\ZZ)$, we deduce that
	\bes
	\cone_k'=\langle \{c_{M(n)}\,:\,\text{$n\in\ZZ_{>0}$}\}\rangle_{\QQ_{\ge0}}.
	\ees
	As basis of $M^k_2(\QQ)$, we choose
	\bes
	E^k_2,E^k_{2,1}(f_1),\dots,E^k_{2,1}(f_\ell),F_1,\dots,F_{\ell'},
	\ees
	where $f_1,\dots,f_\ell$ is a basis of $S^k_1(\QQ)$ and $F_1,\dots,F_{\ell'}$ is a basis of $S^k_2(\QQ)$, on which we rewrite the functionals $c_{M(n)}$ as
	\be\label{eq;proofptck'}
	c_{M(n)}=\big(a^k_1(n), c_n(f_1), \dots, c_n(f_\ell), 0, \dots, 0\big)^t\in\QQ^{\dim(M^k_2)}.
	\ee
	Here we used the well-known fact that the Fourier coefficient of the Siegel Eisenstein series $E^k_2$ associated to $M(n)$ is the $n$-th coefficient of the elliptic Eisenstein series $E^k_1$.
	Analogously, the coefficients of the Klingen Eisenstein series $E^k_{2,1}(f)$ associated to the matrix $M(n)$ coincide with the coefficients $c_n(f)$, for all elliptic cusp forms $f$.
	In fact, the images of~$E^k_2$ and $E^k_{2,1}(f)$ \textcolor{\myblack}{under} the Siegel~$\Phi$-operator are respectively~$E^k_1$ and~$f$; see e.g.~\cite[Section 5]{kli;intro}.
	
	Let $\widetilde{\cone}_k$ be the cone of coefficient extraction functionals of \emph{elliptic} modular forms defined~as
	\bes
	\widetilde{\cone}_k=\langle c_n\,:\,n\in\ZZ_{>0} \rangle_{\QQ_{\ge0}}\subset M^k_1(\QQ)^*.
	\ees
	It is clear from \eqref{eq;proofptck'} that the cone $\cone_k'$ is the embedding in $\QQ^{\dim M^k_2}$ of the cone $\widetilde{\cone}_k$ written over the basis $E^k_1,f_1,\dots,f_\ell$.
	The latter is pointed, rational, polyhedral, and of dimension~$\dim M^k_1$ by \cite[Theorem 3.4]{brmo}. Hence, also $\cone_k'$ satisfies the same properties.
	\end{proof}
	
	\section{The accumulation rays of the modular cone}\label{sec;accrays}
	
	We fix, once and for all, \textbf{a weight $\mathbf{k>4}$ such that $\mathbf{k\equiv 2}$ mod $\mathbf{4}$}.
	The purpose of this section is to classify the accumulation rays of the modular cone $\cone_k$.
	For simplicity, we represent the functionals $c_T$ over a chosen basis of $M^k_2(\QQ)$ of the form
	\be\label{eq;basisMk2}
	E^k_2,E^k_{2,1}(f_1),\dots,E^k_{2,1}(f_\ell),F_1,\dots,F_{\ell'},
	\ee	
	where the Klingen Eisenstein series $E^k_{2,1}(f_j)$ are associated to a basis~$f_1,\dots,f_\ell$ of elliptic cusp forms of~$S^k_1(\QQ)$, and $F_1,\dots,F_{\ell'}$ is a basis of Siegel cusp forms of~$S^k_2(\QQ)$.
	With respect to the basis \eqref{eq;basisMk2}, we may rewrite the functional~$c_T$ as column vectors
	\bes
	c_{T}=\big(a^k_2(T), a^k_2(f_1,T), \dots, a^k_2(f_\ell,T), c_T(F_1), \dots, c_T(F_{\ell'})\big)^t\in\QQ^{\dim M^k_2}.
	\ees
	Recall that we denote by $a^k_2(T)$ and $a^k_2(f_j,T)$ the $T$-th Fourier coefficient associated to $E^k_2$ and $E^k_{2,1}(f_j)$ respectively, in contrast with the coefficients $c_T(F_i)$ of cusp forms.
	
	By Proposition \ref{thm;bodapetn}, if $T=\big(\begin{smallmatrix}
	n & r/2\\
	r/2 & m
	\end{smallmatrix}\big)$,
	the coefficients $a^k_2(f_j,T)$ can be decomposed in Eisenstein and cuspidal parts as
	\bes
	a^k_2(f_j,T)=\frac{\zeta(1-k)}{2}\cdot\sum_{t^2|m}\alpha_m(t,f_j)a^k_2(\Tdiv{T}{t})+c_{n,r}\big((\phi_m^{f_j})^0\big).
	\ees
	The notation used in this decomposition is the same of Section \ref{sec;siegmodforms}.
	In particular, the auxiliary function $\alpha_m$ is defined as in \eqref{eq;alphamaux}, while we denote by $(\phi_m^{f_j})^0$ the cuspidal part of the $m$-th Fourier--Jacobi coefficient associated to $E^k_{2,1}(f_j)$.
	We recall that the matrix denoted by~$\Tdiv{T}{t}$ is constructed from $T$ as defined in \eqref{not;matTdiv}.
	
	Since $k\equiv 2$ mod $4$, the first entry of $c_T$ is \emph{positive}, namely~$a^k_2(T)>0$, by Lemma~\ref{classsiegeis}.
	This implies we can rewrite the ray $\RR_{\ge0}\cdot c_T$, dividing the generator $c_T$ by $a^k_2(T)$, as
	\be\label{eq;decompcef}
	\RR_{\ge0}\cdot c_T=
	\RR_{\ge0}\cdot\begin{pmatrix}
	1 \\
	\const\cdot\sum_{t^2|m}\alpha_m(t,f_1)\frac{a^k_2(\Tdiv{T}{t})}{a^k_2(T)}+\frac{c_{n,r}\big((\phi_m^{f_1})^0\big)}{a^k_2(T)} \\ \vdots \\
	\const\cdot\sum_{t^2|m}\alpha_m(t,f_\ell)\frac{a^k_2(\Tdiv{T}{t})}{a^k_2(T)}+\frac{c_{n,r}\big((\phi_m^{f_\ell})^0\big)}{a^k_2(T)}
	\\ \frac{c_{T}(F_1)}{a^k_2(T)} \\ \vdots \\ \frac{c_{T}(F_{\ell'})}{a^k_2(T)}\end{pmatrix},
	\ee
	where we simply write $\const$ instead of the negative constant $\frac{\zeta(1-k)}{2}$.
	\begin{defi}
	We denote by $\mathcal{S}_k$ the section of the modular cone $\cone_k$ obtained by intersecting it with the hyperplane of points with first coordinate $1$.	
	Equivalently, it is the convex subset in $\cone_k$ of functionals with value $1$ on $E^k_2$. 
	\end{defi}
	We present some basic properties of $\mathcal{S}_k$ and $\cone_k$ in the following result.
	\begin{prop}\label{lemma;compactness}
	The section $\mathcal{S}_k$ is bounded and the modular cone $\cone_k$ is pointed of maximal dimension.
	\end{prop}
	\begin{proof}
	If~$\mathcal{S}_k$ is unbounded, then there exists a sequence of matrices $\seq{T_j}{j\in\NN}$ in~$\halfint^+_2$ such that one of the entries of the point $\mathcal{S}_k\cap\QQ_{\ge0}\cdot c_{T_j}$ diverges when $j\to\infty$.
	This means that either~${|a^k_2(f,T_j)/a^k_2(T_j)|\to\infty}$ or~$c_{T_j}(F)/a^k_2(T_j)\to\infty$. Both cases are impossible, the former by Proposition \ref{prop;growthKEvsSE} and Lemma \ref{classsiegeis} \ref{classsiegeis;magn}, the latter by Remark~\ref{rem;siegcuspgrow}.
	
	The rays in $\overline{\conemod}$ associated to the generators~$c_T$ intersect~$\overline{\mathcal{S}_k}$ in exactly one point; see~\eqref{eq;decompcef}. Since~$\overline{\mathcal{S}_k}$ is compact, \emph{all} rays of~$\overline{\cone_k}$ intersect~$\overline{S_k}$ in one point. These observations imply that~$\overline{\conemod}$ (hence~$\conemod$) is pointed.
	
	We prove now that $\dim\cone_k=\dim M^k_2$. It is enough to show that the functionals $c_T$ associated to matrices $T\in \halfint^+_2$ generate $M^k_2$ over $\CC$. Suppose that this is false. Then, there exists a non-zero $F\in M^k_2$ such that $c_T(F)=0$ for every $T\in\halfint^+_2$. Such Siegel modular forms are called \emph{singular}. It is well-known that, for $k>4$ even, there are no non-zero singular modular forms; see e.g.\ \cite[Section 8, Theorem 2]{kli;intro}. This implies the claim.
	\end{proof}
	We want to classify all possible accumulation rays of the modular cone $\cone_k$.
	We recall that a ray $r$ in $\overline{\conemod}$ is an accumulation ray of $\conemod$ (with respect to the generators appearing in Definition~\ref{def;modcones}) if there exists a family of matrices~$\seq{T_j}{j\in\NN}$ in~$\halfint^+_2$ such that the functionals~$c_{T_j}$ are pairwise different, and the sequence of rays $\seq{\RR_{\ge0}\cdot c_{T_j}}{j\in\NN}$ converges to~$r$.
	
	To classify the accumulation rays of~$\conemod$, we proceed as follows.
	Let~$\seq{c_{T_j}}{j\in\NN}$ be a sequence of pairwise different functionals associated to positive definite matrices~$T_j\in\halfint^+_2$.
	Since~$c_{T_j}=c_{u^t\cdot T_j\cdot u}$ for every $u\in \GL_2(\ZZ)$, we may suppose without loss of generality that the matrices $T_j=\big(\begin{smallmatrix}
	n_j & r_j/2\\
	r_j/2 & m_j
	\end{smallmatrix}\big)$ are \emph{reduced}, i.e.\ the entries satisfy $0\le r_j\le m_j\le n_j$ for every~$j$; see Remark~\ref{rem:reducedmat}.
	For every fixed determinant $d$, there are finitely many reduced matrices~$T$ in~$\halfint^+_2$ with $\det T=d$.
	Since the functionals $c_{T_j}$ are	pairwise different, the matrices $T_j$ have \emph{increasing determinant}, i.e.\ $\det T_j\to\infty$ when $j\to\infty$.
	Suppose that the sequence of rays~$\seq{\RR_{\ge0}\cdot c_{T_j}}{j}$ converges. We classify the accumulation rays arising from such sequences with respect to the chosen family of reduced matrices~$\seq{T_j}{j\in\NN}$.
	In Section~\ref{subsec;accrays}, we treat the cases where the entries~$m_j$ are \textcolor{\myblack}{bounded}.
	In Section~\ref{sec;mnonfixed}, we treat the \textcolor{\myblack}{complementary} cases where the entries~$m_j$ \textcolor{\myblack}{diverge}.
	
	Along the way, we illustrate also some properties that the accumulation rays satisfy.
	These are translated into properties of the points of intersection of~$\overline{\mathcal{S}_k}$ with the accumulation rays. In fact, by Proposition~\ref{lemma;compactness}, also the accumulation rays intersect~$\overline{\mathcal{S}_k}$ in one point.
	
	\subsection{The case of~\textcolor{\myblack}{$\boldsymbol{m_j}$ bounded}}\label{subsec;accrays}
	\textcolor{\myblack}{We firstly consider the case of~$m_j$ eventually constant.}
	
	We fix, once and for all, a positive integer~$m$. Let~$\seq{T_j}{j\in\NN}$ be a sequence of reduced matrices~$T_j=\big(\begin{smallmatrix}
	n_j & r_j/2\\
	r_j/2 & m
	\end{smallmatrix}\big)$ in~$\halfint^+_2$, of increasing determinant. Suppose that the sequence of rays $\seq{\RR_{\ge0}\cdot c_{T_j}}{j}$ is convergent. We rewrite these rays as in \eqref{eq;decompcef}, over the chosen basis \eqref{eq;basisMk2}.
	We already observed in Remark \ref{rem;siegcuspgrow} that
	\bes
	\frac{c_{T_j}(F_s)}{a^k_2(T_j)}\xrightarrow[\,j\to\infty\,]{}0,
	\ees
	for every $s=1,\dots,\ell'$.
	Since the matrices $T_j$ are reduced, by Lemma~\ref{classsiegeis}~\ref{classsiegeis;magn} and Proposition~\ref{prop;Ocusppart} we deduce that analogously
	\bes
	\frac{c_{n_j,r_j}\big((\phi_m^{f_s})^0\big)}{a^k_2(T_j)}\xrightarrow[\,j\to\infty\,]{}0,
	\ees
	for every $s=1,\dots,\ell$.
	Up to considering a sub-sequence of~$\seq{T_j}{j\in\NN}$, we may suppose that the ratios~$a^k_2(\Tdiv{T_j}{t})\big/a^k_2(T_j)$ converge for every square-divisor~$t$ of~$m$.
	In fact, these ratios are bounded between~$0$ and~$1$ by Lemma~\ref{lem;lemratnew}.
	We denote by~$\lambda_t$ the associated limits of ratios.
	These observations imply that
	\be\label{eq;accbarSR}
	\RR_{\ge0}\cdot c_{T_j}\xrightarrow[\,j\to\infty\,]{}
	\RR_{\ge0}\cdot\underbrace{\left(\begin{smallmatrix}
	1 \\
	\const\cdot\sum_{t^2|m}\lambda_t\alpha_m(t,f_1)\\ \vdots \\
	\const\cdot\sum_{t^2|m}\lambda_t\alpha_m(t,f_\ell)\\
	0\\ \vdots \\ 0
	\end{smallmatrix}\right)}_{\in\overline{\mathcal{S}_k}}.
	\ee
	\begin{defi}\label{defi;Qmlambda}
	Let $1=t_0<t_1<\dots<t_\sst$ be the positive integers whose squares divide~$m$. We denote by $\Qm$ the point of intersection of $\overline{\mathcal{S}_k}$ and the accumulation ray obtained in \eqref{eq;accbarSR}.
	If~$m$ is squarefree, we simply write~$\Qmsqf$.
	\end{defi}
	In the notation $\Qm$, there is no need to keep track neither of $\lambda_{t_0}=\lambda_1$, since it is always equal to $1$, nor of the chosen sequence of matrices $\seq{T_j}{j\in\NN}$.
	Note that~$(\lambda_{t_1},\dots,\lambda_{t_d})$ is a \emph{tuple of limits} in~$\limitmm{m}{k}$, as studied in Section~\ref{subsec;ratFCeis}; see Definition~\ref{def;tupleoflim} for more details.
	
	In the remainder of this section, we explain the geometric properties of the accumulation rays~$\RR_{\ge0}\cdot \Qm$ in $\overline{\cone_k}$ \textcolor{\myblack}{by means of} the ones of the points $\Qm$ on~$\overline{\mathcal{S}_k}$.
	We firstly introduce a piece of notation.
	\begin{defi}\label{defi;Vs}
	For every positive integer $s$, we define the point $V_s\in\QQ^{\dim M^k_2}$ as
	\bes
	V_s=\Big(1, \const\cdot\alpha_s(1,f_1),\dots, \const\cdot\alpha_s(1,f_\ell), 0, \dots, 0\Big)^t.
	\ees
	\end{defi}

	The points~$V_s$ are contained in~$\overline{\mathcal{S}_k}$. In fact, consider a sequence of reduced matrices~$\seq{T_j}{j\in\NN}$ in~$\halfint^+_2$ with increasing determinant, such that the bottom-right entry is fixed to~$m$ as above.
	If the entry $r_j$ of $T_j$ is eventually non-divisible by any square-divisor of $m$ different from $1$, then the sequence of rays $\seq{\RR_{\ge0}\cdot c_{T_j}}{j\in\NN}$ converges to the accumulation ray~${\RR_{\ge0}\cdot Q_m(0,\dots,0)}$.
	The point~$V_m$ coincides exactly with~$Q_m(0,\dots,0)$.
	Hence~$\RR_{\ge0}\cdot V_m$ is always an accumulation ray of the modular cone~$\cone_k$.
	
\begin{figure}[h]
	\centering
	\includegraphics[scale=0.6]{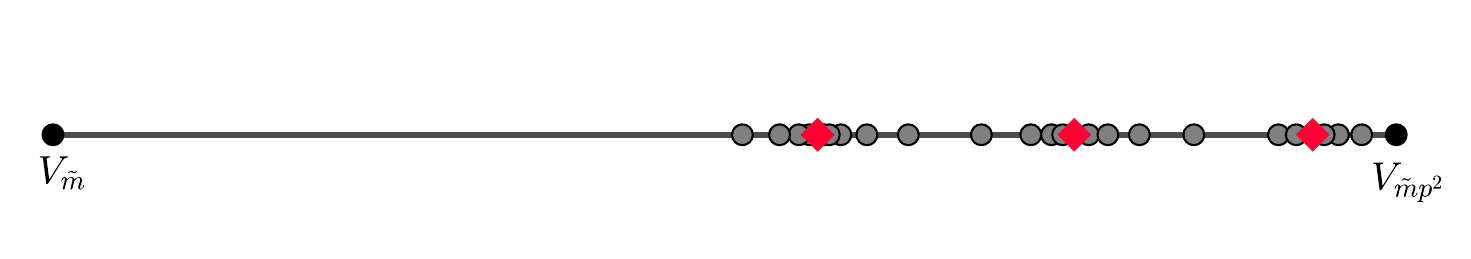}
	\caption{\emph{An idea of the convex hull generated by $V_{\tilde{m}}$ and $V_{\tilde{m}p^2}$, where $\tilde{m}$ is a positive squarefree integer and $p$ is a prime.
	The grey points represent the infinitely many points~$Q_{\tilde{m}p^2}(\lambda_p)$.
	These points accumulate towards some~$Q_{\tilde{m}p^2}(\lambda_p')$, in red, where~$\lambda_p'$ is a special limit in~$\limitmsp{m}{k}{p}$. These are in finite number by Remark \ref{rem;finnuminlimitmspp}; see Section~\ref{subsec;examplesmfixed} for further information.}}
	\label{fig;segment_with_Qm}
	\end{figure}

\begin{figure}[h]
	\centering
	\includegraphics[scale=0.5]{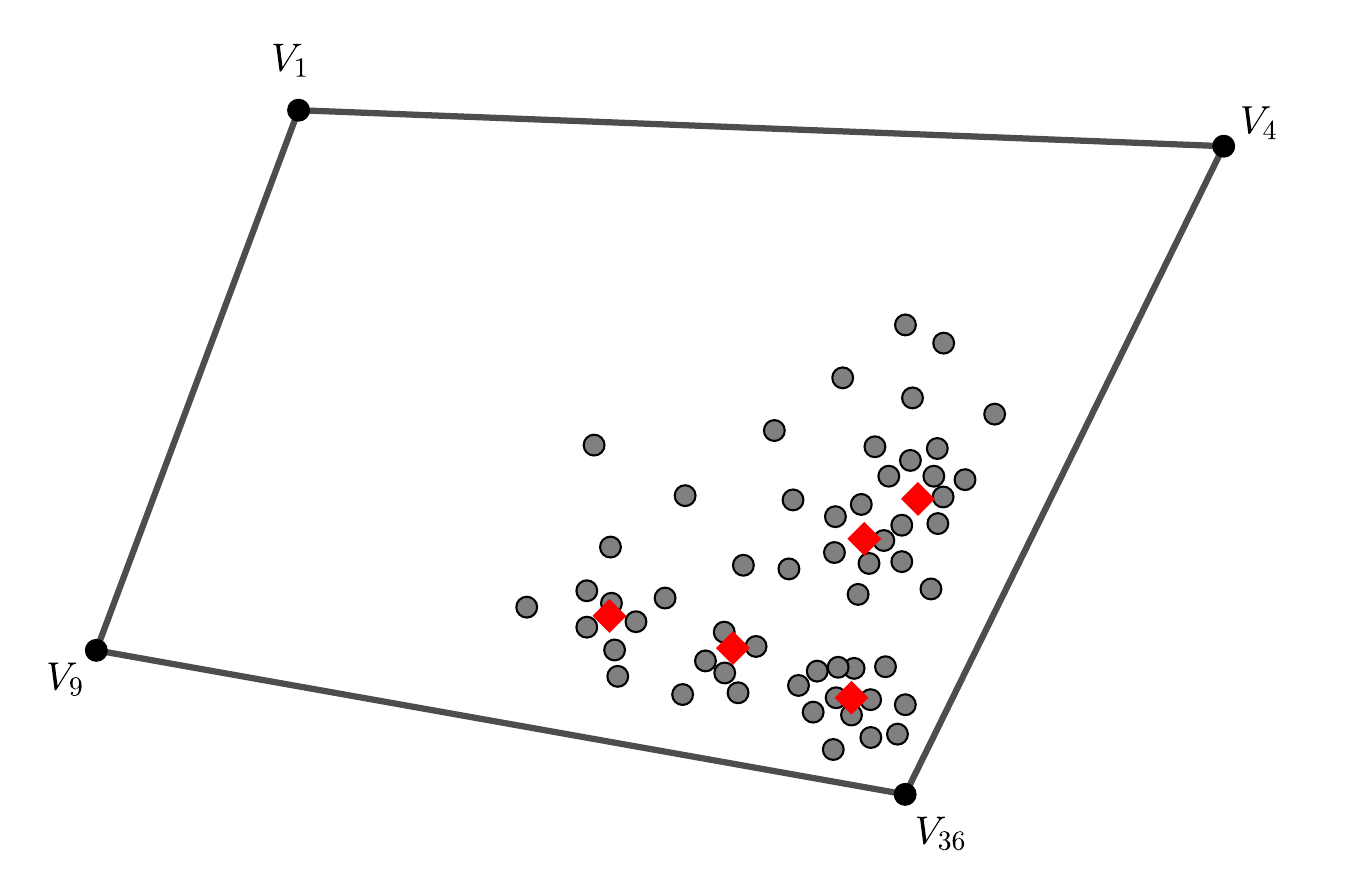}
	\caption{\emph{An idea of the convex hull generated by $V_1$, $V_4$, $V_9$ and $V_{36}$.
	The grey points are some of the infinitely many points $Q_4(\lambda_2,\lambda_3,\lambda_6)$.
	The red points are some of the points towards which the $Q_4(\lambda_2,\lambda_3,\lambda_6)$ accumulate. The number of these can be infinite, depending on the arrangement of the vertexes of the convex hull; see Section~\ref{subsec;examplesmfixed} for further information.}}
	\label{fig;conv_hull_m_36}
	\end{figure}

	We remark that if $m$ is non-squarefree, there are \emph{infinitely many} $\lambda_{t}$ arising as limits of ratios as above; see Proposition \ref{prop;infmanylp} and Corollary \ref{cor;infmanytup}. 
	We are going to prove that, nevertheless, for every $m$, the points $\Qm$ are always contained in the convex hull of \emph{finitely many} $V_s$ for some $s\le m$; see Theorem \ref{thm;convhull}.
	This is essential to prove that the accumulation cone of~$\cone_k$ is rational polyhedral.	
	\begin{lemma}\label{lemma;decQmlambdas}
	Let~$\vect{\lambda}=\tupleoflim\in\limitmm{k}{m}$.
	The point~$\Qmsqf(\vect{\lambda})$ may be written as
	\be\label{eq;sumQmintermsofVm}
	\Qmsqf(\vect{\lambda})
	=
	\sum_{j=0}^\sst \bigg(\sum_{\{t_i\,:\,t_j|t_i\}}\mu\Big(\frac{t_i}{t_j}\Big)\lambda_{t_i}\bigg)\cdot V_{m/t_j^2},
	\ee
	where $\mu$ is the Möbius function.
	\textcolor{\myblack}{In particular, the sum of the coefficients multiplying the points~$V_{m/t_j^2}$ on the right-hand side of~\eqref{eq;sumQmintermsofVm} equals~$1$}. 
	\end{lemma}
	\begin{proof}
	For every $f\in S^k_1$, we may rewrite the defining sum~\eqref{eq;alphamaux} of the auxiliary function~$\alpha_m$ to deduce that
	\begin{align*}
	\sum_{j=0}^\sst\lambda_{t_j}\alpha_m(t_j,f)
	=
	\sum_{j=0}^\sst\sum_{\ell_j|t_j}\mu\Big(\frac{t_j}{\ell_j}\Big)\lambda_{t_j}\alpha_{m/\ell_j^2}(1,f)=
	\sum_{j=0}^\sst\bigg(\sum_{\{t_i\,:\,t_j|t_i\}}\mu\Big(\frac{t_i}{t_j}\Big)\lambda_{t_i}\bigg)\alpha_{m/t_j^2}(1,f).
	\end{align*}
	If evaluated in~$f=f_i$, the left-hand side of the previous formula gives the~$i+1$ entry of the vector~$\Qmsqf(\vect{\lambda})$ up to the factor~$\zeta$.
	Since the value~$\zeta\cdot\alpha_{m/t_j^2}(1,f_i)$ is the~$i+1$ entry of~$V_{m/t_j^2}$, it remains to show that \textcolor{\myblack}{the first entry of the tuple on the left-hand side of~\eqref{eq;sumQmintermsofVm}, which equals~$1$, is the same as the first entry of the right-hand side of~\eqref{eq;sumQmintermsofVm}.
	This is equivalent of showing that} the sum of the coefficients multiplying the~$V_{m/t_j^2}$'s on the right-hand side of~\eqref{eq;sumQmintermsofVm} equals~$1$.
	This is an easy check, since
	\begin{align}\label{eq;inlemmadecQmlambdas}
	\sum_{j=0}^\sst\sum_{\{t_i\,:\,t_j|t_i\}}\mu\left(\frac{t_i}{t_j}\right)\lambda_{t_i}
	=
	\sum_{i=0}^\sst\lambda_{t_i}\sum_{\ell|t_i}\mu\left(\frac{t_i}{\ell}\right)
	=
	1+\sum_{i=1}^\sst\lambda_{t_i}\sum_{\ell|t_i}\mu(\ell)=1.
	\end{align}
	Here we used that if~$\ell$ divides~$t_i$, then~$\ell=t_j$ for some $j\le i$, together with the well-known formula~$\sum_{a|b}\mu(a)=\delta_{b,1}$.
	\end{proof}
	\begin{thm}\label{thm;convhull}
	Let~$\vect{\lambda}=\tupleoflim\in\limitmm{k}{m}$.
	The points $\Qmsqf(\vect{\lambda})$ lie in the convex hull over~$\RR$ generated by the points~$V_{m/t_j^2}$ for $j=0,\dots,\sst$.
	\end{thm}
	To make the previous result as clear as possible, in Section~\ref{subsec;examplesmfixed} we compute explicitly the convex hull in $\overline{\cone_k}$ generated by the points $\Qmsqf(\vect{\lambda})$, for a few $m$;
	see Figures \ref{fig;segment_with_Qm} and \ref{fig;conv_hull_m_36} as examples of such convex hulls.
	\begin{proof}	
	By Lemma~\ref{lemma;decQmlambdas}, the points~$\Qmsqf(\vect{\lambda})$ are linear combinations of the points~$V_{m/t_j^2}$ for~${j=0,\dots,\sst}$.
	We check that the coefficients of these combinations fulfill the definition of convex hull, i.e.\ their sum is one and they are non-negative; see the introduction of Section~\ref{sec;backgrcones}.
	The fact that their sum equals~$1$ has already been checked \textcolor{\myblack}{with Lemma~\ref{lemma;decQmlambdas}}.
	We now check the non-negativity.
	Decompose $m=v^2\widetilde{m}$, where $\widetilde{m}$ is squarefree. Let $v=p_1^{a_1}\cdots p_b^{a_b}$ be the prime decomposition of $v$. Choose a positive integer $t$ such that $t^2|m$. This implies that $t|v$ and $t=p_1^{s_1}\cdots p_b^{s_b}$ for some $0\le s_j\le a_j$, where $j=1,\dots,b$.
	We want to show that
	\begin{equation}\label{eq;thmconvhull}
	\sum_{j=1}^b\sum_{x_j=s_j}^{a_j}\mu\left(p_1^{x_1-s_1} \cdots p_b^{x_b-s_b}\right)\lambda_{p_1^{x_1}\cdots p_b^{x_b}}\ge0.
	\end{equation}
	\textcolor{\myblack}{
	As described in Corollary~\ref{cor;ratwithsiegsernonprime}, there exist~$\lambda_{{p_j}^{x_j}}\in\limitm{m}{k}{{p_j}^{x_j}}$ such that~$\lambda_{p_1^{x_1}\cdots p_b^{x_b}}=\prod_{j}\lambda_{{p_j}^{x_j}}$.
	Since also the Möbius function~$\mu$ is multiplicative, to prove~\eqref{eq;thmconvhull} it is enough to show that
	\begin{equation}\label{eq;thmconvhullbis}
	\sum_{x_j=s_j}^{a_j}\mu\big({p_j}^{x_j-s_j}\big)\lambda_{{p_j}^{x_j}}\ge0,
	\end{equation}
	for every~$j=1,\dots,b$.}
	To verify this property, we prove the analogous inequality where
	\bes
	\textcolor{\myblack}{\lambda_{{p_j}^{x_j}}\quad\text{is replaced with}\quad a^k_2\left(\begin{smallmatrix}
	n & r/2{p_j}^{x_j}\\
	r/2{p_j}^{x_j} & m/{p_j}^{2x_j}
	\end{smallmatrix}\right)\Big/ a^k_2\left(\begin{smallmatrix}
	n & r/2\\
	r/2 & m
	\end{smallmatrix}\right),}
	\ees
	where $\big(\begin{smallmatrix}
	n & r/2\\
	r/2 & m
	\end{smallmatrix}\big)\in\halfint^+_2$. In fact, the former is the limit of a sequence of ratios of Fourier coefficients as the latter, which are \textcolor{\myblack}{\emph{non-negative} by Lemma \ref{classsiegeis}}.
	
	\textcolor{\myblack}{Recall the notation $\Tdiv{T}{x}=\big(\begin{smallmatrix}
	n & r/2x\\
	r/2x & m/x^2
	\end{smallmatrix}\big)$, for every $T=\big(\begin{smallmatrix}
	n & r/2\\
	r/2 & m
	\end{smallmatrix}\big)\in\halfint_2^+$.}
	Without loss of generality, we may assume that~$s_j<a_j$.
	\textcolor{\myblack}{Since~$\mu\big({p_j}^{x_j-s_j}\big)=0$} whenever $x_j-s_j\ge2$, \textcolor{\myblack}{to prove~\eqref{eq;thmconvhullbis} it is enough to show that}
	\ba\label{eq;thmconvhulleasyrev}
	\textcolor{\myblack}{a^k_2({T}^{[{p_j}^{s_j}]}) - a^k_2(T^{[{p_j}^{s_j+1}]}) \ge0.}
	\ea
	Since $k\equiv 2$ mod $4$ by hypothesis, the coefficients of the Siegel Eisenstein series appearing in \eqref{eq;thmconvhulleasyrev}, which are always evaluated on positive definite matrices, are either positive or zero.
	The latter case happens only when the upper-right entries of the \textcolor{\myblack}{matrices~$\Tdiv{T}{{p_j}^{s_j}}$ and~$\Tdiv{T}{{p_j}^{s_j+1}}$ are} not half-integral.
	
	Consider the first summand~\textcolor{\myblack}{$a^k_2(\Tdiv{T}{{p_j}^{s_j}})$} of~\eqref{eq;thmconvhulleasyrev}. If the entry $r$ of $T$ is not divisible by~\textcolor{\myblack}{${p_j}^{s_j}$}, then~\textcolor{\myblack}{$a^k_2(\Tdiv{T}{{p_j}^{s_j}})=0$, and so is also all the remaining summand on the left-hand side of~\eqref{eq;thmconvhulleasyrev}.}
	If instead~$\textcolor{\myblack}{{p_j}^{s_j}}|r$, then~\textcolor{\myblack}{$a^k_2(\Tdiv{T}{{p_j}^{s_j}})>0$}.
	
	Consider the second summand~\textcolor{\myblack}{$a^k_2(T^{[{p_j}^{s_j+1}]})$ on} the left-hand side of~\eqref{eq;thmconvhulleasyrev}.
	If $\textcolor{\myblack}{{p_j}^{s_j+1}}\nmid r$, then this term is zero.
	By Lemma~\ref{lem;lemratnew}, if~$\textcolor{\myblack}{{p_j}^{s_j+1}}| r$, then~\textcolor{\myblack}{${a^k_2(\Tdiv{T}{{p_j}^{s_j+1}})<a^k_2(\Tdiv{T}{{p_j}^{s_j}})}$}.
	In fact~\textcolor{\myblack}{$F_p(\Tdiv{T}{{p_j}^{s_j+1}},3-k)<F_p(\Tdiv{T}{{p_j}^{s_j}},3-k)$} for every prime $p$, as shown in Section~\ref{subsec;ratFCeis}.
	This implies~\eqref{eq;thmconvhulleasyrev}.
	\end{proof}
	\begin{cor}\label{cor;polconacc}
	Let $\seq{\RR_{\ge0}\cdot c_{T_j}}{j\in\NN}$ be a convergent sequence of rays, where $T_j\in\halfint^+_2$ are reduced, of increasing determinant and with the bottom-right entries eventually equal to some positive $m$. The accumulation ray of the modular cone $\conemod$ obtained as limit of such sequence is contained in the subcone $\langle V_{m/t^2}\,:\,t^2|m\rangle_{\RR_{\ge0}}$ of $\overline{\mathcal{C}_k}$, which is \emph{rational polyhedral}.
	\end{cor}
	\begin{proof}
	The limit of $\RR_{\ge0}\cdot c_{T_j}$ is as in \eqref{eq;accbarSR}, that is, it is generated by~$\Qmsqf(\vect{\lambda})$, for some~${\vect{\lambda}\in\limitmm{k}{m}}$.
	By Theorem~\ref{thm;convhull}, this point is contained in the convex hull generated by the $V_{m/t^2}$, where~$t$ runs among the positive integers whose squares divide $m$. The points $V_s$ have rational entries for every $s$, because so are the values of $\alpha_s(1,f)$ for every $f\in S^k_1(\QQ)$. The polyhedrality of the cone generated by the $V_{m/t^2}$ is trivial, since these points are in finite number.
	\end{proof}
	
		\textcolor{\myblack}{Suppose now that the bottom-right entries $m_j$ of $T_j$ oscillate among a finite set of positive integers, and that the sequence of rays $\seq{\RR_{\ge0}\cdot c_{T_j}}{j\in\NN}$ converges.
	Then, the accumulation ray obtained as a limit for $j\to\infty$ must be $\RR_{\ge0}\cdot Q_{\widetilde{m}}(\lambda_{t_1},\dots,\lambda_{t_d})$ for some $\widetilde{m},\lambda_{t_1},\dots,\lambda_{t_d}$.
	In fact, consider the sub-sequence~$\seq{T_i}{i}$ of~$\seq{T_j}{j\in\NN}$ where the matrices~$T_i$ have the entry~$m_i$ fixed to one of the values appearing infinitely many time as bottom-right entry of~$T_j$, say~$\widetilde{m}$.
	We have seen above that the limit of $\RR_{\ge0}\cdot c_{T_i}$, for $i\to\infty$, must be generated by~$Q_{\widetilde{m}}(\vect{\lambda})$ for some tuple of limits $\vect{\lambda}\in\limitmm{\widetilde{m}}{k}$.}
	
	
	\subsection{The case of \textcolor{\myblack}{$\boldsymbol{m_j}$ divergent}}\label{sec;mnonfixed}
	The aim of this section is to describe the geometric properties of the accumulation rays of the modular cone $\cone_k$ arising as limits of sequences~$\seq{\RR_{\ge0}\cdot c_{T_j}}{j\in\NN}$, where $T_j$ are reduced matrices in $\halfint^+_2$ of increasing determinant, such that the bottom-right \textcolor{\myblack}{entries~$m_j$ diverge}.
	For this reason, this section may be considered as the complementary of Section \ref{subsec;accrays}, where the bottom-right entries were \textcolor{\myblack}{bounded}.
	\begin{defi}
	We define the point $P_\infty\in\QQ^{\dim M^k_2}$ as
	\[
	P_\infty=(1,0,\dots,0)^t.
	\]
	\end{defi}
	The point $P_\infty$ lies in $\overline{\mathcal{S}_k}$, as follows from the next result.
	\begin{lemma}\label{lem;VstoPinf}
	The points $V_s\in\overline{\mathcal{S}_k}$ converge to $\Pinf$, when $s\to\infty$.
	\end{lemma}
	\begin{proof}
	It is sufficient to prove that if $s\to\infty$, then $\alpha_s(1,f)\to 0$ for every elliptic cusp form~${f\in S^k_1}$. It is straightforward to check that
	\begin{align*}
	|\alpha_s(1,f)|&=\left|\frac{\gf{f}{s}}{g_k(s)}\right|=
	\frac{\left|\sum_{d^2|s}\mu(d)c_{s/d^2}(f)\right|}{s^{k-1}\prod_{p|s}(1+p^{-k+1})}\le\frac{\sum_{d^2|s}|c_{s/d^2}(f)|}{s^{k-1}}\le\\
	&\le\frac{\sigma_0(s)\cdot\max_{1\le y \le s}|c_y(f)|}{s^{k-1}}=O_f\left(s^{\frac{2+\varepsilon-k}{2}}\right),
	\end{align*}
	for all $\varepsilon>0$.
	The last equality is deduced using the classical Hecke-bound for Fourier coefficients of elliptic cusp forms and the well-known property $\sigma_0(s)=o(s^\varepsilon)$ for all $\varepsilon>0$.
	Since~$k>4$, the claim follows.
	\end{proof}
	\begin{prop}\label{prop;infmanyVs}
	If $k\ge 18$, then the modular cone $\cone_k$ has \emph{infinitely many} accumulation rays.
	\end{prop}
	\begin{proof}
	Since this is a counting problem of rays of~$\overline{\cone_k}$, we may work with~$M^k_2(\RR)$ in place of~$M^k_2(\QQ)$, namely with Siegel modular with real Fourier coefficients.
	In this way, we may choose the basis~\eqref{eq;basisMk2} such that the elliptic cusp forms~$f_1,\dots,f_\ell$ are (normalized) Hecke forms.
	
	Since $k\ge 18$, then~$\dim S^k_1(\RR)>0$.
	We firstly note that the point~$V_1$ on $\overline{\mathcal{S}_k}$ is different from~$\Pinf$.
	In fact, suppose that it is not, then $c_1(f_j)=0$ for every $j=1,\dots,\ell$.
	This means that $c_1(f_j)=0$, but this is not possible since the Hecke forms $f_j$ are normalized with $c_1(f_j)=1$.
	Suppose that there is only a finite number of accumulation rays of $\cone_k$. Since $\RR_{\ge0}\cdot V_s$ is an accumulation ray for every positive integer~$s$, also the number of points~$V_s$ must be finite.
	By Lemma~\ref{lem;VstoPinf}, the points~$V_s$ converge to~$\Pinf$ when~$s\to\infty$.
	This implies that there exists a positive integer~$s_0$ such that~$V_s=\Pinf$ for every~$s\ge s_0$.
	Suppose that~$s\ge s_0$ is \emph{squarefree}.
	We deduce from the equality~$V_s=\Pinf$ that~$c_s(f_j)=0$ for every~$j$, in particular
	\be\label{eq;prinfmanyVs}
	c_s(f_j)=0\qquad \text{for every $s\ge s_0$ squarefree}.
	\ee
	It is known that the coefficients of normalized Hecke eigenforms satisfy
	\bes
	c_{p^{\nu+1}}(f_j)=c_p(f_j)\cdot c_{p^\nu}(f_j)-p^{k-1}c_{p^{\nu-1}}(f_j),
	\ees
	for every prime $p$ and every $\nu\ge1$; see e.g.\ \cite[Part~I, Section~4.2]{1-2-3}.
	Let $p$ be a prime number greater than $s_0$. Since~$V_{p^2}$ coincides with $\Pinf$, then $c_{p^2}(f_j)-c_1(f_j)=0$ and
	\bes
	0=c_{p^2}(f_j)-1=\big(c_p(f_j)\big)^2-p^{k-1}-1.
	\ees
	We deduce that $c_p(f_j)$ is non-zero.
	Hence the relation \eqref{eq;prinfmanyVs} can not be satisfied by~$c_p(f_j)$, for any prime~$p\ge s_0$.
	This implies that there are infinitely many~$V_s$.
	\end{proof}
	The following result concludes the classification of \emph{all} possible accumulation rays in $\overline{\mathcal{C}_k}$.
	\begin{prop}\label{prop;acraysminf}
	Let $\seq{T_j}{j\in\NN}$ be a sequence of reduced matrices in $\halfint^+_2$ of increasing determinant, such that the bottom-right entries $m_j$ diverge when $j\to\infty$. The sequence of rays $\seq{\RR_{\ge0}\cdot c_{T_j}}{j\in\NN}$ converges to $\RR_{\ge0}\cdot \Pinf$.
	\end{prop}
	\begin{proof}
	As usual, we consider every functional $c_T$ as a point in $\QQ^{\dim M^k_2}$, writing it with respect to the basis~\eqref{eq;basisMk2}.
	It is enough to prove that $c_{T_j}/a^k_2(T_j)\to \Pinf$, when $j\to \infty$.
	Since the cuspidal parts of the entries of $c_{T_j}$ grow slower than $a^k_2(T_j)$ when $j\to\infty$, the accumulation ray obtained as limit of $\RR_{\ge0}\cdot c_{T_j}$ depends only on the Eisenstein parts of the entries of $c_{T_j}$; see Remark~\ref{rem;siegcuspgrow} and Proposition \ref{prop;Ocusppart}. 
	Analogously to the proof of Lemma~\ref{lem;VstoPinf}, we may compute
	\begin{align*}
	\bigg|\sum_{t^2|m}\alpha_m(t,f)\frac{a^k_2(\Tdiv{T}{t})}{a^k_2(T)}\bigg|\le
	\sum_{t^2|m}\left|\alpha_m(t,f)\right|\le
	\sigma_0(m)\cdot\max_{t^2|m}|\alpha_m(t,f)|\le \hspace{2cm}\\
	\le\sigma_0(m)\cdot \max_{t^2|m}\sum_{s|t}\left|\frac{\gf{f}{m/s^2}}{g_k(m/s^2)}\right|\le
	\sigma_0(m)\cdot\max_{t^2|m}\left(\sigma_0(t)\cdot\max_{s|t}\left|\frac{\gf{f}{m/s^2}}{g_k(m/s^2)}\right|\right)=\\
	=
	O_f\left(m^{\frac{2+\varepsilon-k}{2}}\right),
	\end{align*}
	for every $f\in S^k_1$ and every $\varepsilon>0$, when $\det T\to\infty$.
	Here we used the well-known property~${\sigma_0(s)=o(s^\varepsilon)}$, for all $\varepsilon>0$, the Hecke bound for elliptic cusp forms, and the inequality 
	\bes
	0\le a^k_2(\Tdiv{T}{t})/a^k_2(T)\le1,
	\ees
	for all positive integers $t$ whose squares divide $m$.
	Since~$k>4$, the claim follows.
	\end{proof}
	
	\section{The accumulation cone of the modular cone is rational and polyhedral}\label{sec;accumulcone}
	
	We recall that a ray of the $\RR$-closure $\overline{\cone_k}$ is an \emph{accumulation ray} of the modular cone $\conemod$ (with respect to the set of generators appearing in Definition~\ref{def;modcones}), if it is the limit of some sequence of rays $\seq{\RR_{\ge0}\cdot c_{T_j}}{j\in\NN}$, where $T_j\in\halfint^+_2$ are reduced and of increasing determinant; see Section \ref{sec;backgrcones}.
	The \emph{accumulation cone} of $\cone_k$ is the cone generated by the accumulation rays of $\cone_k$.
	We denote it by $\accone_k$.
	
	By the classification of the accumulation rays of~$\conemod$ given in Section~\ref{sec;accrays}, in particular by Corollary~\ref{cor;polconacc} and Proposition~\ref{prop;acraysminf}, the cone~$\accone_k$ may be generated as
	\be\label{eq;genaccone}
	\accone_k=\langle\Pinf, V_s\,:\, s\ge1\rangle_{\RR_{\ge0}}.
	\ee 
	The goal of this section is to prove the following result.
	\begin{thm}\label{thm;acconeratpol}
	If $k>4$ and $k\equiv 2$ mod $4$, then the accumulation cone $\accone_k$ of the modular cone $\cone_k$ is rational and polyhedral, of the same dimension as $M^k_1$.
	\end{thm}
	We firstly present some preparatory results.
	As in Section \ref{sec;accrays}, we consider all coefficient extraction functionals~$c_T$ as vectors in~$\QQ^{\dim M^k_2}$, written over a fixed basis of~$M^k_2(\QQ)$ of the form~\eqref{eq;basisMk2}.
	\begin{defi}\label{defi;Ps}
	For every positive integer $s$, we define the point $P_s\in\QQ^{\dim M^k_2}$ as
	\bes
	P_s= \Big(1,\const\cdot\frac{c_s(f_1)}{\ssigma(s)},\dots,\const\cdot\frac{c_s(f_\ell)}{\ssigma(s)},0,\dots,0\Big)^t,
	\ees
	where we write $\const$ instead of the negative constant $\frac{\zeta(1-k)}{2}$.
	\end{defi}
	We remark that whenever~$s$ is squarefree, the point~$P_s$ coincides with the point~$V_s$ defined in Section~\ref{sec;accrays}.
	The points~$P_s$ are contained in the section~$\overline{\mathcal{S}_k}$, as showed by the following result.
	Recall the auxiliary function~$g_k$ from~\eqref{eq;gkauxfunct}.
	\begin{prop}\label{Pmnonsqfr}
	Let $s$ be a positive integer. The point $P_s$ satisfies the relation
	\be\label{eq;relPm}
	P_s=\sum_{t^2|s}\frac{g_k(s/t^2)}{\ssigma(s)}V_{s/t^2}.
	\ee
	In particular, the point $P_s$ lies in the convex hull $\conv_\RR\left(\{V_{s/t^2}:t^2|s\}\right)$.
	\end{prop}
	To make Proposition \ref{Pmnonsqfr} as clear as possible, we propose a direct check of \eqref{eq;relPm} in Section~\ref{subsec;examplesmfixed} for a few choices of $m$.
	\begin{proof}
	\textcolor{\myblack}{Let~$s$ be a positive integer.
	We denote by~$\gamma_{t,s}$ the rational numbers defined recursively as~$\gamma_{1,s}=g_k(s)/\sigma_{k-1}(s)$ and}
	\bes
	\gamma_{t,s}=-\frac{1}{\ssigma(s)}\sum_{1\ne d|t}\mu(d)\ssigma(s/d^2)\gamma_{t/d,s/d^2},
	\ees
	\textcolor{\myblack}{for all square divisors~$t>1$ of~$s$.}
	
	\textcolor{\myblack}{We show that}
	\begin{equation}\label{eq;proofPmintern}
	\textcolor{\myblack}{\frac{c_s(f)}{\ssigma(s)}=\sum_{t^2|s}\gamma_{t,s}\alpha_{s/t^2}(1,f),\qquad\text{for every $f\in S^k_1$,}}
	\end{equation}
	by induction on the number $\sqdiv(s)$ of square-divisors of $s$.
	Suppose that~${\sqdiv(s)=1}$, then~$s$ is squarefree and $c_s(f)/\ssigma(s)=\alpha_s(1,f)$.
	\textcolor{\myblack}{Since}~$\gamma_{1,s}=1$, the desired relation is fulfilled.
	
	Suppose now that $\sqdiv(s)>1$ and that~\eqref{eq;proofPmintern} is satisfied for every positive integer~$\tilde{s}$ such that~$\sqdiv(\tilde{s})<\sqdiv(s)$.
	We \textcolor{\myblack}{rewrite the right-hand side of} \eqref{eq;proofPmintern} as
	\begin{align*}
	\textcolor{\myblack}{\frac{c_s(f)}{\ssigma(s)}}
	+
	\textcolor{\myblack}{\sum_{1\neq t^2|s}\frac{\mu(t)c_{s/t^2}(f)}{\ssigma(s)}}
	+
	\sum_{1\neq t^2|s}\gamma_{t,s}\alpha_{s/t^2}(1,f),
	\end{align*}
	\textcolor{\myblack}{then we apply induction and rewrite it as}
	\bes
	\textcolor{\myblack}{\frac{c_s(f)}{\ssigma(s)}}+\sum_{1\neq t^2|s}\frac{\mu(t)\ssigma(s/t^2)}{\ssigma(s)}\sum_{\tilde{t}^2|\frac{s}{t^2}}\gamma_{\tilde{t},s/t^2}\alpha_{s/t^2\tilde{t}^2}(1,f)+\sum_{1\neq t^2|s}\gamma_{t,s}\alpha_{s/t^2}(1,f).
	\ees
	\textcolor{\myblack}{Using the recursive definition of~$\gamma_{t,s}$, it is easy to check that the sum of second and the third summand of the previous formula is zero.
	This concludes the proof of~\eqref{eq;proofPmintern}.}
	
	\textcolor{\myblack}{We now prove} that $\gamma_{t,s}=g_k(s/t^2)/\ssigma(s)$, by induction on the number of divisors~$\ndiv(t)$ of~$t$.
	If~$\ndiv(t)=1$, then~${t=1}$ and the claim is true by definition, for every $s$. Suppose now that $\ndiv(t)>1$, then
	\begin{align*}
	\gamma_{t,s}&=-\frac{1}{\ssigma(s)}\sum_{1\ne d|t}\mu(d)\ssigma(s/d^2)\frac{g_k(s/t^2)}{\ssigma(s/d^2)}=-\frac{g_k(s/t^2)}{\ssigma(s)}\sum_{1\ne d|t}\mu(d)=\frac{g_k(s/t^2)}{\ssigma(s)},
	\end{align*}
	where we used induction on $\gamma_{t/d,s/d^2}$, since $\ndiv(t/d)<\ndiv(t)$ whenever~${d\ne1}$.
	
	To conclude the proof, we show that the coefficients $\gamma_{t,s}$ satisfy the requirements which make $P_s$ a point of the convex hull~$\conv_\RR\big(\{V_{s/t^2}:t^2|s\}\big)$.
	Firstly, we prove that~${\sum_{t^2|s}\gamma_{t,s}=1}$ by induction on the number of square-divisors~$\sqdiv(s)$ of~$s$.
	This is equivalent to \textcolor{\myblack}{proving} that $\sum_{t^2|s} g_k(s/t^2)=\ssigma(s)$. Suppose that $\sqdiv(s)=1$, then~$s$ is squarefree and~$\sum_{t^2|s}g_k(s/t^2)=g_k(s)=\ssigma(s)$.
	If~$\sqdiv(s)>1$, then
	\begin{align*}
	\sum_{t^2|s}g_k(s/t^2)=\sum_{t^2|s}\sum_{y^2|\frac{s}{t^2}}\mu(y)\ssigma(s/t^2y^2)=\sum_{x^2|s}\ssigma(s/x^2)\sum_{d|x}\mu(d)=\ssigma(s).
	\end{align*}
	Eventually, since $g_k(s)>0$ for every positive integer $s$, so is $\gamma_{t,s}$ for every $t^2|s$.
	
	Since \eqref{eq;proofPmintern} is true for every elliptic cusp form of weight $k$, it is true also for the chosen basis $f_1,\dots,f_\ell$ of $S^k_1(\QQ)$.
	The evaluation of \eqref{eq;proofPmintern} in $f=f_j$ verifies the~$(j+1)$-th entry of the equality \eqref{Pmnonsqfr}. The check for the remaining entries is trivial.
	\end{proof}
	\begin{cor}
	For every positive integer $s$, the ray $\RR_{\ge0}\cdot P_s$ lies in the rational polyhedral subcone~$\langle V_{s/t^2}\,:\,t^2|s\rangle_{\RR_{\ge0}}$ of $\overline{\mathcal{C}_k}$.
	\end{cor}	
	\begin{proof}
	The polyhedrality is a trivial consequence of Proposition~\ref{Pmnonsqfr}.
	Since the basis~\eqref{eq;basisMk2} is made of Siegel modular forms with rational Fourier coefficients, we deduce that the subcone is rational.
	\end{proof}
	\begin{cor}\label{cor;realdimaccone}
	The (real) dimension of $\accone_k$ is equal to the (complex) dimension of $M^k_1$.
	\end{cor}
	\begin{proof}
	The cone $\accone_k\subseteq \RR^{\dim M^k_2}$ is generated by vectors where only the first $1+\ell$ entries can be different from zero, as we can see from~\eqref{eq;genaccone}. Since $1+\ell=\dim M^k_1$, it is clear that $\dim\accone_k\le \dim M^k_1$. Let $\widetilde{\cone}_k$ be the cone generated over $\RR$ by the coefficient extraction functionals of $M^k_1$, that is
	\bes
	\widetilde{\cone}_k=\langle c_s\,:\,s\in\ZZ_{\ge1}\rangle_{\RR_{\ge0}}.
	\ees
	We consider the functionals $c_s$ as vectors in $\RR^{\dim M^k_1}$, represented over the basis $E^k_1,f_1,\dots,f_\ell$, where $E^k_1$ is the normalized elliptic Eisenstein series of weight~$k$, and~$f_1,\dots,f_\ell$ is the basis of~$S^k_1(\QQ)$ chosen in~\eqref{eq;basisMk2}.
	The entries of~$c_s/c_s(E^k_1)$ are the first~$1+\ell$ entries of~$P_s$. This means that the linear map
	\bes
	\iota\colon\widetilde{\cone}_k\longrightarrow\accone_k,\quad c_s/c_s(E^k_1)\longmapsto P_s
	\ees
	is an embedding. Hence, we have also $\dim\accone_k\ge\dim M^k_1$.
	\end{proof}
	\begin{lemma}\label{prop;100isinternalinaccone}
	The point $\Pinf$ is \emph{internal} in $\accone_k$.
	\end{lemma}
	\begin{proof}
	The idea is to rewrite $\Pinf$ as a linear combination with positive coefficients of enough points $P_s$, such that these generate a subcone of $\accone_k$ with maximal dimension. By Lemma~\ref{lemma;propbmc}, there exist a constant $A$ and positive coefficients $\eta_j$ with $1\le j\le A$, such that
	\bes
	\sum_{j=1}^A\eta_j c_j|_{S^k_1(\QQ)}=0\qquad\text{in $S^k_1(\QQ)^*$}.
	\ees
	We recall that~$A$ can be chosen arbitrarily large.
	
	The entries of~$P_s$ associated to the basis~$f_1,\dots,f_\ell$ of~$S^k_1(\QQ)$ are, up to multiplying by the negative constant~$\const/\ssigma(s)$, the values of the functional~$c_s$ on~$f_1,\dots,f_\ell$.
	This implies that
	\bes
	\sum_{j=1}^A\eta_j \ssigma(j) P_j=\Pinf\sum_{j=1}^A\eta_j\ssigma(j).
	\ees
	Since the points $P_s$ are contained in $\accone_k$ by Proposition \ref{Pmnonsqfr}, also $\Pinf$ is contained therein.
	By Corollary \ref{cor;realdimaccone}, we may take $A$ big enough such that the dimension of $\langle P_j\,:\,1\le j\le A\rangle_{\RR_{\ge0}}$ is the same as the one of $\accone_k$.
	In this way, the point $\Pinf$ is internal in $\accone_k$ with respect to the euclidean topology.	
	\end{proof}
	We are ready to illustrate the proof of the main result of this section.
	\begin{proof}[Proof of Theorem \ref{thm;acconeratpol}]
	Suppose that $\accone_k$ is not polyhedral, that is, it has infinitely many extremal rays.
	Since $\accone_k$ is generated by $\Pinf$ and the points $V_s$ with~$s$ positive, and these points accumulate only towards~$\Pinf$ by Lemma~\ref{lem;VstoPinf}, there are infinitely many extremal rays of the form~$\RR_{\ge0}\cdot V_{s'}$, for some~$s'>0$.
	These \textcolor{\myblack}{boundary} rays accumulate towards~$\RR_{\ge0}\cdot \Pinf$, hence \textcolor{\myblack}{also} the latter \textcolor{\myblack}{is} a boundary ray.
	But this is in contrast with Lemma \ref{prop;100isinternalinaccone}.
	Therefore, the cone~$\accone_k$ is polyhedral.
	
	The extremal rays are generated by some of the points $V_s$, which have rational entries. Hence, the cone is rational.
	
	The statement about the dimension of $\accone_k$ is Corollary \ref{cor;realdimaccone}.
	\end{proof}
	
	\section{Additional properties of the modular cone}\label{sec;propC}

	In this section, which is a focus on the geometric properties of the modular cone~$\cone_k$, we generalize some of the results used in Section \ref{sec;accumulcone} to prove that $\accone_k$ is rational polyhedral. The problem of the polyhedrality of $\cone_k$ is more complicated. The issue is to understand how a sequence of rays~$\RR_{\ge0}\cdot c_{T_j}$ converges to an accumulation ray of $\cone_k$, depending on the choice of the family of reduced matrices~$\seq{T_j}{j\in\NN}$ in~$\halfint^+_2$ with increasing determinant.
	We will translate the polyhedrality of~$\conemod$ into a conjecture on Fourier coefficients of Jacobi cusp forms.
	
	\textbf{We fix once and for all a weight $\mathbf{k>4}$ such that $\mathbf{k\equiv 2}$ mod $\mathbf{4}$}, and consider the functionals $c_T$ as vectors in $\QQ^{\dim M^k_2}$ over the basis~\eqref{eq;basisMk2}.
	We begin with the properties of~$\cone_k$ which are a direct consequence of the results in the previous sections.
	We remark that these properties, together with Proposition~\ref{lemma;compactness}, give the previously announced points~\ref{pt3;conek'} and~\ref{pt4;conek'} of Theorem~\ref{thm;mainthm2}.
	\begin{prop}
	The modular cone~$\cone_k$ is rational, and intersects the rank~$1$ modular cone~$\cone_k'$ only at the origin.
	Moreover, if the cone $\cone_k$ is enlarged with a non-zero vector of~$\cone_k'$, the resulting cone is non-pointed.
	\textcolor{\myblack}{In particular~$-\conemodone\subseteq\overline{\conemod}$}
	\end{prop}
	\begin{proof}
	Since the generators $c_T$ are functionals over the space of Siegel modular forms with \emph{rational} Fourier coefficients, the rationality of $\cone_k$ follows trivially by the rationality of its accumulation cone, namely by Theorem~\ref{thm;acconeratpol}.
	
	If we rewrite the functionals with respect to the usual basis~\eqref{eq;basisMk2} of~$M^k_2(\QQ)$, we deduce by Remark~\ref{rem;coefklingcuspeq} that
	\bas
	c_{\left(\begin{smallmatrix}
	s & 0\\ 0 & 0
	\end{smallmatrix}\right)}=&
	\left(a^k_2\left(\begin{smallmatrix}
	s & 0\\ 0 & 0
	\end{smallmatrix}\right),a^k_2\left(f_1,\left(\begin{smallmatrix}
	s & 0\\ 0 & 0
	\end{smallmatrix}\right)\right),\dots, a^k_2\left(f_\ell,\left(\begin{smallmatrix}
	s & 0\\ 0 & 0
	\end{smallmatrix}\right)\right), 0,\dots,0\right)=\\
	=&\left(c_s(E^k_1),c_s(f_1),\dots,c_s(f_\ell),0\dots,0\right)=
	2\sigma_{k-1}(s)/\zeta(1-k)\cdot P_s,
	\eas
	for every positive integer $s$.
	Since $k\equiv 2$ mod $4$, the constant $\zeta(1-k)$ is negative, hence
	\bes
	\RR_{\ge0}\cdot c_{\left(\begin{smallmatrix}
	s & 0\\ 0 & 0
	\end{smallmatrix}\right)}=\RR_{\ge0}\cdot (-P_s).
	\ees
	The ray $\RR_{\ge0}\cdot P_s$ is contained in $\overline{\cone_k}$ by Proposition \ref{Pmnonsqfr}.
	This implies that whenever we enlarge the cone $\overline{\cone_k}$ with one of the generators of $\cone_k'$, which are the functionals $c_T$ associated to non-zero singular matrices, the resulting cone contains also $\RR_{\ge0}\cdot (-P_s)$ for some $s$.
	Since the whole line $\RR\cdot P_s$ is contained in the enlarged cone, the latter is non-pointed.
	This is sufficient to conclude the proof, since a rational cone in a finite-dimensional vector space over $\QQ$ is pointed if and only if its $\RR$-closure is pointed.
	\end{proof}
	In Section \ref{sec;accumulcone} we proved that $\Pinf$ is \emph{internal} in the accumulation cone $\accone_k$ of $\cone_k$; see Lemma~\ref{prop;100isinternalinaccone}.
	This played a key role for the proof of the polyhedrality of $\accone_k$.
	In the following result, we prove that $\Pinf$ lies in the interior of the $\RR$-closure $\overline{\cone_k}$.
	Note that it does not follow from Lemma~\ref{prop;100isinternalinaccone}, and it does not imply it. In fact, the cones~$\accone_k$ and~$\overline{\cone_k}$ may have different dimensions.
	\begin{prop}\label{prop;100isinternalinSk}
	The point $\Pinf$ is \emph{internal} in $\overline{\conemod}$.
	\end{prop}
	We know that the dual space $(M^k_2)^*$ is generated over~$\CC$ by the functionals~$c_T$ with~${T\in\halfint^+_2}$.
	For the proof of Proposition~\ref{prop;100isinternalinSk}, we need to restrict the set of these generators to the ones indexed by an auxiliary subset of~$\halfint^+_2$, as showed by the following result.
	It follows from the noteworthy fact that Siegel modular forms are determined by their Fundamental Fourier coefficients; see~\cite{saha} and~\cite[Section~$7.2$]{boda;petnorm}.
	\begin{lemma}\label{thm;semplific}
	Let $\halfint_2'$ be the subset of $\halfint^+_2$ containing all matrices with squarefree bottom-right entry.
	The dual space~$(M^k_2)^*$ is generated by the functionals~$c_T$ with~$T\in\halfint'_2$.
	\end{lemma}
	\begin{proof}[Proof of Lemma \ref{thm;semplific}]
	We prove the result showing that if~$F\in M^k_2\setminus\{0\}$, then the Fourier coefficient~$c_T(F)$ is non-zero for an infinite number of matrices~$T\in\halfint'_2$.
	We follow closely the proofs of~\cite[Theorem~$1$]{saha} and~\cite[Proposition~$7.7$]{boda;petnorm}.
	
	\textbf{Cuspidal case:}
	Suppose that~$F$ is a Siegel cusp form.
	By~\cite[Proposition~$2.2$]{saha}, there exists an odd prime~$p$ such that the~$p$-th Fourier--Jacobi coefficient~$\phi_p$ of~$F$ is non-zero\textcolor{\myblack}{; see also~\cite[Proof of the Main Theorem]{kriegraum}}.
	In fact, see~\cite[p.~$369$]{saha}, the Jacobi cusp form~$\phi_p$ has an infinite number of non-zero Fourier coefficients.
	More precisely, they are of the form~$c_{(D+\mu^2)/4p,\mu}(\phi_p)$, where~$D$ and~$\mu$ are integers, and~$D$ is odd and squarefree.
	Such coefficients equals the ones of~$F$ corresponding to the matrices~$\big(\begin{smallmatrix}
	(D+\mu^2)/4p & \mu/2\\
	\mu/2 & p
	\end{smallmatrix}\big)$, which are contained in~$\halfint'_2$.
	
	\textbf{Non-cuspidal case:}
	Suppose that~$F\in M^k_2\setminus S^k_2$.
	By Lemma~\ref{classsiegeis}, the property of~$F$ we want to prove is satisfied if~$F=E^k_2$.
	Therefore, without loss of generality, we may suppose that the Siegel Eisenstein part of~$F$ is trivial.
	This means we may rewrite~$F$ as~$F=E^k_{2,1}(F) + G$, for some~$f\in S^k_1\setminus\{0\}$ and~$G\in S^k_2$.
	As illustrated in~\cite[p.~$369$]{boda;petnorm}, it is possible to construct a sequence of matrices in~$\halfint^+_2$ of the form~$T_j=\big(\begin{smallmatrix}
	n_j & 1/2\\
	1/2 & m_j
	\end{smallmatrix}\big)$, for some squarefree~$m_j$ and of increasing determinant, with the property that~$c_{T_j}(F)$ diverges when~$j\to\infty$.
	\end{proof}
	\begin{proof}[Proof of Proposition \ref{prop;100isinternalinSk}]
	The idea of the proof is the following.
	We rewrite~$\Pinf$ as a linear combination with positive coefficients of some $P_s$, as in the proof of Lemma \ref{prop;100isinternalinaccone}.
	Then, we rewrite some of those $P_s$ associated to squarefree \textcolor{\myblack}{indices} as linear combinations with positive coefficients of some functionals $c_T$.
	We will take these combinations in such a way that the subcone generated by those $c_T$ has maximal dimension into the $\RR$-closure $\overline{\cone_k}$.
	
	As we have already shown in the proof of Lemma~\ref{prop;100isinternalinaccone}, by Lemma~\ref{lemma;propbmc} there exist an arbitrarily large constant~$A$, and positive coefficients $\eta_m$, with $1\le m\le A$, such that
	\be\label{eq;dim100int}
	\Pinf=\sum_{m=1}^A\eta_m P_m.
	\ee
	By Lemma \ref{lem:relcoefextjac}, for every positive~$m$ there exist an arbitrarily large constant~$B_m$ and positive coefficients~$\mu_{n,r}^m$ such that
	\be\label{eq;decompsomprinc}
	\sum_{1\le n\le B_m}\sum_{\substack{r\in\ZZ\\ 4nm-r^2>0}}\mu_{n,r}^mc_{n,r}|_{J_{k,m}^{\text{cusp}}(\QQ)}=0.
	\ee
	We recall that $c_{n,r}|_{J_{k,m}^{\text{cusp}}(\QQ)}$ is the functional in~$J^{\text{cusp}}_{k,m}(\QQ)^*$ which extracts the~$(n,r)$-th Fourier coefficient of Jacobi cusp forms in~$J^{\text{cusp}}_{k,m}(\QQ)$.
	Note that the sum appearing in~\eqref{eq;decompsomprinc} is finite.
	
	Suppose that~$m$ is fixed \emph{squarefree}, and write for simplicity~$T^m_{n,r}$ instead of~$\big(\begin{smallmatrix}n & r/2\\ r/2 & m\end{smallmatrix}\big)$.
	\textcolor{\myblack}{Under} the usual decomposition of the entries of~$c_{T^m_{n,r}}$ as explained at the beginning of Section~\ref{sec;accrays}, we deduce from~\eqref{eq;decompsomprinc} that
	\ba\label{eq;explformlarge100isinternalinSk}
	\sum_{1\le n\le B_m}\sum_{\substack{r\in\ZZ\\ 4nm-r^2>0}}\mu_{n,r}^mc_{T^m_{n,r}}
	=&
	\left(\begin{smallmatrix}
	\sum_n\sum_r \mu^m_{n,r} a^k_2(T^m_{n,r}) \\ \const\cdot \frac{c_m(f_1)}{\ssigma(m)}\cdot\sum_n\sum_r \mu^m_{n,r}a^k_2(T^m_{n,r}) \\
	\vdots \\
	\const\cdot \frac{c_m(f_\ell)}{\ssigma(m)}\cdot \sum_n\sum_r \mu^m_{n,r} a^k_2(T^m_{n,r}) \\ 0 \\ \vdots \\ 0
	\end{smallmatrix}\right)
	+
	\underbrace{\left(\begin{smallmatrix}
	0 \\  \sum_n\sum_r \mu^m_{n,r} c_{n,r}\big((\phi_m^{f_1})^0\big)\\ \vdots \\ \sum_n\sum_r 
	\mu^m_{n,r} c_{n,r}\big((\phi_m^{f_\ell})^0\big)\\ \sum_n\sum_r \mu^m_{n,r} c_{T^m_{n,r}}(F_1) \\ \vdots \\ \sum_n\sum_r \mu^m_{n,r} c_{T^m_{n,r}}(F_{\ell'})
	\end{smallmatrix}\right)}_{=0}=\\
	=&P_m\cdot\sum_{1\le n\le B_m}\sum_{\substack{r\in\ZZ\\ 4nm-r^2>0}} \mu^m_{n,r} a^k_2(T^m_{n,r}).
	\ea
	The matrices $T^m_{n,r}$ appearing in the previous equation are contained in $\halfint^+_2$, that is, they are positive definite.
	We define
	\bes
	\tilde{\eta}(m,B_m)=\frac{\eta_m}{\sum_{1\le n\le B_m}\sum_{\substack{r\in\ZZ\\ 4nm-r^2>0}} \mu_{n,r}^m a^k_2(T^m_{n,r})},
	\ees
	for every~$m$ squarefree. The value~$\tilde{\eta}(m,B_m)$ is \emph{positive} by Lemma~\ref{classsiegeis}~\ref{lem;pt1class}.
	We can then further decompose $\Pinf$ from \eqref{eq;dim100int} into
	\bas
	\Pinf=&
	\sum_{1\le m\le A}\eta_m P_m=\\
	=&\sum_{\substack{1\le m\le A\\ \text{$m$ non-squarefree}}}\eta_m P_m+\sum_{\substack{1\le m\le A\\ m \text{ squarefree}}}\tilde{\eta}(m,B_m)\sum_{1\le n\le B_m}\sum_{\substack{r\in\ZZ\\ 4nm-r^2>0}} \mu^m_{n,r} c_{T^m_{n,r}}.
	\eas
	Up to \textcolor{\myblack}{choosing} $A$ and $B_m$ large enough, the functionals $c_{T^m_{n,r}}$ appearing in the previous decomposition of $\Pinf$ generate over $\QQ$ the whole $M^k_2(\QQ)^*$ by Lemma \ref{thm;semplific}.
	This implies that~$\Pinf$ lies in the interior of $\overline{\mathcal{S}_k}$.
	\end{proof}
	We now focus on Conjecture~\ref{conj;our}, namely the problem of the polyhedrality of the~$\RR$-closure~$\overline{\cone_k}$.
	The cone $\overline{\cone_k}$ is polyhedral if and only if it has finitely many extremal rays, or equivalently if its extremal rays do not accumulate anywhere.
	An accumulation ray arising as limit of a sequence of extremal rays is a boundary ray of $\cone_k$, but not necessarily extremal; see Example \ref{ex;acextrnotextr}.
		
\begin{figure}[h]
	\centering
	\includegraphics[scale=0.7]{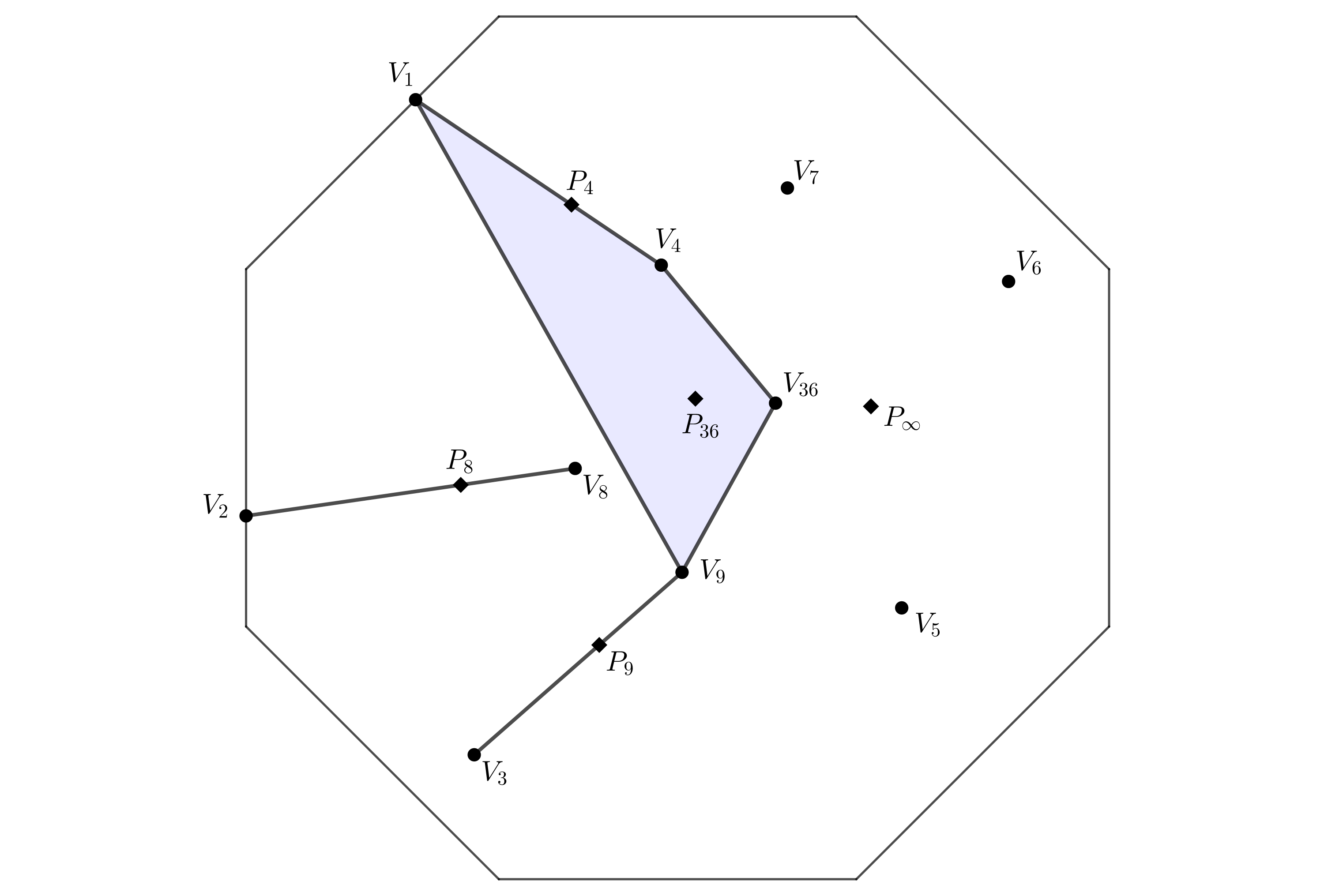}
	\caption{\emph{An idea of the section $\overline{\mathcal{S}_k}$ under the hypothesis of Theorem \ref{thm;polyhedwithhyp}, with highlighted the convex hull associated to~${m=36}$}.}
	\label{fig;conhulinS}
	\end{figure}

	The first, although hopeless, idea to prove that the extremal rays of~$\cone_k$ do not accumulate, is to show that \emph{all} accumulation rays of $\cone_k$ are generated by points lying in the interior of $\overline{\cone_k}$, as we did for the accumulation ray $\RR_{\ge0}\cdot \Pinf$ in Proposition \ref{prop;100isinternalinSk}.
	We checked with SageMath \cite{sage} that, for large weights, this is \emph{false}, since some of the accumulation rays $\RR_{\ge0}\cdot V_s$ may lie in the boundary of $\overline{\cone_k}$.
	The following example collects some of these empirical observations.
	The computation of the coefficients of Siegel modular forms was carried out with the package~\cite{stakemori}.
	\begin{ex}\label{ex;boundarycontacray}
	Suppose that $k>4$ and $k\equiv 2$ mod $4$. We provide in the following table some of the accumulation rays of the modular cone~$\cone_k$ which lie in the boundary of $\overline{\cone_k}$.
	\begin{center}

	\begin{tabular}{ |r||p{8cm}|  }
 	\hline
	 $k$ & Some accumulation rays in the boundary of $\overline{\cone_k}$\\
 	\hline
 	$18$   & $\RR_{\ge0}\cdot V_1$\\
 	$22$, $26$, $30$   & $\RR_{\ge0}\cdot V_1$, $\RR_{\ge0}\cdot V_2$\\
 	$34$, $38$   & $\RR_{\ge0}\cdot V_1$, $\RR_{\ge0}\cdot V_2$, $\RR_{\ge0}\cdot V_3$\\
 	$42$   & $\RR_{\ge0}\cdot V_1$, $\RR_{\ge0}\cdot V_2$, $\RR_{\ge0}\cdot V_3$, $\RR_{\ge0}\cdot V_4$\\
 	\hline
	\end{tabular}
	\end{center}
	\end{ex}
	With Example \ref{ex;boundarycontacray} in mind, we provide an alternative sufficient condition to deduce Conjecture \ref{conj;our}. This is exactly the hypothesis of the following result. The idea is that to deduce the polyhedrality of $\cone_k$, it is enough to show that every accumulation ray of $\cone_k$ is generated by a point which lies in the interior of a subcone of $\overline{\cone_k}$, and that this subcone eventually contains the sequences of rays converging to the chosen accumulation ray.
	\begin{thm}\label{thm;polyhedwithhyp}
	Suppose that for every accumulation ray $\RR_{\ge0}\cdot\Qmsqf(\vect{\lambda})$, with~$\vect{\lambda}\in\limitmm{m}{k}$, there exists a subcone $\subconem(\vect{\lambda})$ of $\overline{\conemod}$ which contains~$\Qmsqf(\vect{\lambda})$ in its interior, and such that any sequence of rays~$\seq{\RR_{\ge0}\cdot c_{T_j}}{j\in\NN}$ converging to $\RR_{\ge0}\cdot\Qmsqf(\vect{\lambda})$, where~$T_j=\big(\begin{smallmatrix}
	n_j & r_j/2\\
	r_j/2 & m
	\end{smallmatrix}\big)$ are reduced in~$\halfint^+_2$ and of increasing determinant, is eventually contained in $\subconem(\vect{\lambda})$.
	Then the modular cone~$\conemod$ is polyhedral.
	\end{thm}
	See Figure \ref{fig;conhulinS} for an idea of the (polyhedral) shape of the section $\mathcal{S}_k$ of $\cone_k$ whenever the hypothesis of Theorem \ref{thm;polyhedwithhyp} are fulfilled.
		
	The rest of this section is devoted to the proof of the previous result and to some remarks on its hypothesis.
	More precisely, we prove the hypothesis of Theorem \ref{thm;polyhedwithhyp} for~$m$ \emph{squarefree}, and we translate that hypothesis for $m$ non-squarefree into a conjecture on Fourier coefficients of Jacobi cusp forms.
	\begin{proof}[Proof of Theorem \ref{thm;polyhedwithhyp}]
	Let~$\{\RR_{\ge0}\cdot a_n\}_{n\in\NN}$ be the set of extremal rays of $\overline{\cone_k}$, where~${a_n\in\overline{\mathcal{S}_k}}$.
	\textcolor{\myblack}{Recall that~$\overline{\mathcal{C}_k}$ is the convex cone defined as the topological closure over~$\RR$ of~$\mathcal{C}_k$.
This implies that~$\overline{\mathcal{C}_k}$ is generated by the generators of~$\mathcal{C}_k$ (i.e.\ the functionals~$c_T$) \emph{together with their accumulation rays}.
We deduce that the rays~$\RR_{\ge0}\cdot a_n$ can be only of the following two types.}
	\begin{enumerate}[label=(\roman*)]
	\item $\RR_{\ge0}\cdot a_n=\RR_{\ge0}\cdot c_{T_n}$ for a suitable coefficient extraction functional $c_{T_n}$, where $T_n$ is a reduced matrix in $\halfint^+_2$.\label{eq;typeproofextrray1}
	\item The ray $\RR_{\ge0}\cdot a_n$ is not of type \ref{eq;typeproofextrray1}, and there exists a sequence~$\seq{T_j}{j\in\NN}$ of reduced matrices in~$\halfint^+_2$ with increasing determinant, such that~$\RR_{\ge0}\cdot c_{T_j}\to\RR_{\ge0}\cdot a_n$ if~$j\to\infty$.\label{eq;typeproofextrray2}
	\end{enumerate}
	The idea is to show that the extremal rays of $\overline{\cone_k}$ are finitely many (hence $\overline{\cone_k}$ is polyhedral) and only of type \ref{eq;typeproofextrray1} (hence $\cone_k$ is polyhedral).
		
	Let $\RR_{\ge0}\cdot a_n$ be an extremal ray of type~\ref{eq;typeproofextrray2}, arising from a sequence of matrices~$\seq{T_j}{j\in\NN}$.
	The bottom-right entries $m_j$ of the matrices $T_j$ can not diverge when $j\to\infty$. In fact, if they diverge, then $a_n=\Pinf$ by Proposition \ref{prop;acraysminf}.
	Since $\Pinf$ is an internal point of $\overline{\mathcal{S}_k}$ by Proposition \ref{prop;100isinternalinSk}, the ray $\RR_{\ge0}\cdot a_n$ is not a boundary ray.
	This is in contrast with the hypothesis that $\RR_{\ge0}\cdot a_n$ is an extremal ray, hence the entries $m_j$ must be bounded.
		
	Let $m$ be one of the values that the entries $m_j$ assume infinitely many times.
	Up to considering a subsequence of~$\seq{T_j}{j\in\NN}$, the ratios~${a^k_2\big(\Tdiv{T_j}{t}\big)\big/a^k_2(T_j)}$ converge to some $\lambda_t$ for every integer $t$ such that $t^2$ divides~$m$.
	Here we use the notation~$\Tdiv{T_j}{t}$ as introduced in~\eqref{not;matTdiv}.
	Denote by~$\vect{\lambda}$ the corresponding tuple of limits in~$\limitmm{m}{k}$.
	Following the same argument of the introduction of Section~\ref{subsec;accrays}, we deduce that
	\bes
	\RR_{\ge0}\cdot c_{T_j}\xrightarrow[j\to\infty]{}\RR_{\ge0}\cdot \Qmsqf(\vect{\lambda}),
	\ees
	in particular $a_n=\Qmsqf(\vect{\lambda})$.
	By hypothesis, there exists a subcone $\subconem(\vect{\lambda})$ of~$\overline{\conemod}$ containing~$a_n$ in its interior.
	If~$\dim\subconem(\vect{\lambda})>1$, then~$\RR_{\ge0}\cdot a_n$ can not be extremal by definition.
	If~$\dim\subconem(\vect{\lambda})=1$, then~$\RR_{\ge0}\cdot a_n$ is of type~\ref{eq;typeproofextrray1}.
	Hence, there are no extremal rays of~$\conemod$ of type~\ref{eq;typeproofextrray2}.
	
	We conclude the proof showing that the extremal rays of type \ref{eq;typeproofextrray1} are finitely many.
	Suppose they are not, that is, there exists a sequence of pairwise different extremal rays~$\seq{\RR_{\ge0}\cdot c_{T_j}}{j\in\NN}$ indexed over a family of reduced matrices~$\seq{T_j}{j\in\NN}$ in~$\halfint^+_2$ of increasing determinant.
	We suppose without loss of generality that this sequence converges to a \emph{boundary ray}~$\RR_{\ge0}\cdot b$ for some~$b\in\overline{\mathcal{S}_k}$.
	The limit does not have to be extremal; see Example~\ref{ex;acextrnotextr}.
	Following the same argument as above, up to considering a subsequence of~$\seq{T_j}{j\in\NN}$, the bottom-right entries of  these matrices are fixed to some positive integer~$m$, and~$b=\Qmsqf(\vect{\lambda})$ for some~$\vect{\lambda}\in\limitmm{m}{k}$.
	By hypothesis, there exists a subcone~$\subconem(\vect{\lambda})$ of~$\overline{\conemod}$ containing~$b$ in its interior, and such that the rays~${\RR_{\ge0}\cdot c_{T_j}}$ eventually lie in~$\subconem(\vect{\lambda})$. Since~${\RR_{\ge0}\cdot c_{T_j}}$ are pairwise different, the dimension of~$\subconem(\vect{\lambda})$ is greater than~$1$. Since these are extremal rays for $\conemod$, they are extremal rays also for~$\subconem(\vect{\lambda})$. But this implies they are boundary rays of~$\subconem(\vect{\lambda})$, hence they can not accumulate towards~$\RR_{\ge0}\cdot b$, since~$b$ is an \emph{internal point} of~$\subconem(\vect{\lambda})$.
	\end{proof}
	\begin{lemma}\label{lemma;msqfrsubcone}
	The hypothesis of Theorem~\ref{thm;polyhedwithhyp} for~$m$ squarefree is always satisfied.
	More precisely, if~$m$ is a positive squarefree integer, the subcone~$\subconem$ associated to the accumulation ray $\RR_{\ge0}\cdot V_m$ exists, and can be chosen as
	\be\label{eq;statmsqfrsubcone}
	\subconem=\langle c_T\,:\, \text{\rm$T\in\halfint^+_2$ with $m$ as bottom-right entry} \rangle_{\RR_{\ge0}}.
	\ee
	\end{lemma}
	\begin{proof}
	The points $V_m$ and $P_m$ coincide, since $m$ is squarefree.
	We have already shown, e.g.\ in~\eqref{eq;explformlarge100isinternalinSk}, that Lemma~\ref{lem:relcoefextjac} implies the existence of an arbitrarily large constant $B_m$ and positive constants $\mu_{n,r}^m$ such that
	\be\label{eq;proofsubcone}
	V_m=\sum_{1\le n \le B_m}\sum_{\substack{r\in\ZZ\\ 4nm-r^2>0}}\mu_{n,r}^m c_{T^m_{n,r}}.
	\ee
	
	We define the subcone $\subconem$ of $\overline{\cone_k}$ as the cone generated by the functionals $c_T$ with~${T\in\halfint^+_2}$, not necessarily reduced, such that the bottom-right entry of $T$ is $m$.
	\textcolor{\myblack}{For every fixed matrix~$T$ as above, we may choose~$B_m$ sufficiently large so that~$T$ appears in~\eqref{eq;proofsubcone}.}
	We may enlarge~$B_m$ such that the matrices~$T_{n,r}^m$ appearing in \eqref{eq;proofsubcone} generate over~$\RR$ a space of dimension equal to~$\dim\subconem$. In this way, the point~$V_m$ is \emph{internal} in~$\subconem$.
	
	If $m$ is squarefree, the unique accumulation ray associated to a sequence of rays~$\RR_{\ge0}\cdot c_{T_j}$, where~$T_j=\big(\begin{smallmatrix}
	n_j & r_j/2\\
	r_j/2 & m
	\end{smallmatrix}\big)$ are reduced matrices in~$\halfint^+_2$ with bottom-right entries fixed to~$m$ and of increasing determinant, is the ray $\RR_{\ge0}\cdot V_m$; see Section \ref{subsec;accrays}.
	All functionals~$c_{T_j}$ are contained in~$\subconem$ by definition.
	\end{proof}
	Lemma~\ref{lemma;msqfrsubcone} and our SageMath computations~\cite{zufprog} are sufficient to deduce the polyhedrality of~$\cone_k$ via Theorem~\ref{thm;polyhedwithhyp} for the very first cases of~$k$.
	\begin{prop}
	If~$k\le 38$, then the cone~$\cone_k$ is polyhedral.
	\end{prop}
	\begin{proof}
	Let~$d$ be a positive integer, and let~$D$ be the determinant of some matrix in~$\halfint^+_2$.
	We denote by~$\cone_k(D,d)$ the polyhedral subcone of~$\cone_k$ defined as
	\bes
	\cone_k(D,d)=\big\langle
	c_T : \text{
	$T=\big(\begin{smallmatrix}
	n & r/2\\ r/2 & m
	\end{smallmatrix}\big)\in\halfint^+_2$ such that~$0\le r\le m\le n\le d$ and~$\det T\le D$
	}
	\big\rangle_{\QQ_{\ge0}},
	\ees
	namely the cone generated by all functionals~$c_T$ associated to reduced matrices with determinant at most~$D$ and diagonal entries at most~$d$.
	
	We consider the functionals~$c_T$ as tuples in~$\RR^{\dim M^k_2}$, representing them over a basis of the form~\eqref{eq;basisMk2}.
	We set~$d=30$ and~$D=d^2=900$, and run the SageMath program~\cite{zufprog} to compute~$\cone_k(900,30)$, for all~$k\le38$, obtaining a full-dimensional cone.
	We then check whether the point~$V_m$ lies on the boundary of~$\cone_k(D,d)$, for every~$m\le 30$.
	The results are as in Example~\ref{ex;boundarycontacray}, in particular we deduce that if~$4\le m\le 30$, then~$V_m$ is in the interior of~$\cone_k(900,30)$.
	Since the latter has the same dimension of~$\cone_k$, the same property of~$V_m$ holds for~$\cone_k$ in place of~$\cone_k(900,30)$.
	
	We check that the set~$\{V_m : m\le 30\}$ generates the whole accumulation cone of~$\cone_k$ as follows.
	By Lemma~\ref{lem;VstoPinf}, the points~$V_m$ converge towards~$P_\infty$ when~$m\to\infty$.
	Without loss of generality, we may assume that the basis~\eqref{eq;basisMk2} used to represent the functionals~$c_T$ is made of Siegel modular forms with \emph{real} Fourier coefficients, so that we may choose~$f_1,\dots,f_\ell$ to be (normalized) Hecke forms.
	By~\cite[Th\'eor\`eme~$18$]{DeligneWeilconj1}, we may bound the Fourier coefficients of~$f_j$ by
	\bes
	|c_n(f_j)|\le \sigma_0(n)\cdot n^{(k-1)/2},\qquad\text{for every~$j$.}
	\ees
	It is easy to see that the distance between the point~$V_m$ and~$P_\infty$, with respect to the basis of~$M^k_2(\RR)$ chosen above, is bounded by
	\be\label{eq;boundchecksage}
	d(V_m,P_\infty)\le \sqrt{\frac{\ell\cdot\zeta(1-k)^2}{m^{k-3}}}.
	\ee
	We implemented in~\cite{zufprog} a function that, for fixed~$k$ and~$m$, checks whether the ball centered in~$P_\infty$ and with ray~$d(V_m,P_\infty)$ is contained in the interior of the cone~$\cone_k(D,d)$.
	We use such function to check that if~$k\le 38$, then all points~$V_m$ are contained in the interior of~$\cone_k(900,30)$, for every~$m\ge 30$.
	Hence there is no other point~$V_m$ that may lie on the boundary of~$\cone_k$, for~$k\le 38$, apart from the ones appearing in Example~\ref{ex;boundarycontacray}.
	
	We conclude the proof showing that the hypothesis of Theorem~\ref{thm;polyhedwithhyp} is always satisfied.
	For the points~$V_1$, $V_2$ and~$V_3$, we may choose~$\subconem$ as in Lemma~\ref{lemma;msqfrsubcone}.
	For all remaining generators~$\Qmsqf(\vect{\lambda})$ of accumulation rays of~$\cone_k$ we may choose~$\subconem(\vect{\lambda})=\cone_k(900,30)$, by what we proved above.
	\end{proof}
	
	Lemma \ref{lemma;msqfrsubcone} verifies the hypothesis of Theorem~\ref{thm;polyhedwithhyp} only for~$m$ squarefree.
	It is natural to ask if the analogue statement of Lemma~\ref{lemma;msqfrsubcone} holds also for~$m$ non-squarefree.
	The following example collects some empirical observations deduced with SageMath.
	These suggest that the subcone~\eqref{eq;statmsqfrsubcone} computed for the squarefree cases can not be used to prove the hypothesis of Theorem~\ref{thm;polyhedwithhyp} for~$m$ non-squarefree.
	\begin{ex}\label{ex;tabV4}
	Choose $4$ as non-squarefree integer. We define the subcone $\mathcal{F}_4$ of $\overline{\cone_k}$ as
	\be\label{eq;tabV4}
	\mathcal{F}_4=\langle c_T\,:\, \text{\rm$T\in\halfint^+_2$ with $4$ as bottom-right entry} \rangle_{\RR_{\ge0}},
	\ee
	in analogy with the subcone~$\subconem$ constructed in Lemma \ref{lemma;msqfrsubcone} for squarefree integers $m$.
	The following table, deduced from some empirical observations in SageMath, shows that~$V_4$ could lie in the boundary of~$\mathcal{F}_4$.
	\begin{center}
	\begin{tabular}{ |c||c|c|c|c| }
 	\hline
	 $k$ & Is $V_4$ contained in $\mathcal{F}_4$? & Is $V_4$ internal in $\mathcal{F}_4$? & $\dim\mathcal{F}_4$ & $\dim\overline{\cone_k}$ \\
 	\hline
 	$18$ & yes & yes & $4$ & $4$\\
 	$22$ & yes & yes & $6$ & $6$\\
 	$26$ & yes & yes & $7$ & $7$\\
 	$30$ & yes & yes & $11$ & $11$\\
 	$34$ & yes & no & $14$ & $14$\\
 	$38$ & yes & no & $15$ & $16$\\
 	$42$ & yes & no & $17$ & $22$\\
 	\hline
	\end{tabular}
	\end{center}
	\end{ex}
	With Example~\ref{ex;tabV4} in mind, we conjecture here a property of Fourier coefficients of Jacobi cusp forms sufficient to deduce a correct generalization of Lemma~\ref{lemma;msqfrsubcone} for every~$m$ non-squarefree.
	This conjecture is actually a refinement of Lemma~\ref{lem:relcoefextjac}.
	
	To simplify the exposition, we consider for a moment the case $m=4$, as in Example~\ref{ex;tabV4}.
	The subcone $\mathcal{F}_4$ of $\overline{\cone_k}$, as defined in \eqref{eq;tabV4}, may contain $V_4$ in its boundary if the weight $k$ is large enough.
	If the cone $\subconewo_4(0)$ associated to $V_4$, as hypothesized in Theorem~\ref{thm;polyhedwithhyp}, exists, then it must be a subcone of $\mathcal{F}_4$ of lower dimension.
	In fact, let~$\seq{\RR_{\ge0}\cdot c_{T_j}}{j\in\NN}$ be a sequence of rays associated to reduced matrices $T_j=\big(\begin{smallmatrix}
	n_j & r_j/2\\
	r_j/2 & 4
	\end{smallmatrix}\big)\in\halfint^+_2$ of increasing determinant.
	Clearly, the rays~${\RR_{\ge0}\cdot c_{T_j}}$ lie in~$\mathcal{F}_4$ for every~$j$.
	We know that~${\RR_{\ge0}\cdot c_{T_j}\to\RR_{\ge0}\cdot V_4}$, when~$j\to\infty$, if and only if~$r_j$ is eventually non-divisible by~$2$.
	This follows from Corollary~\ref{cor;ratwithsiegser} applied with~${p=2}$ and~$\lambda_2=0$, since~$V_4=Q_4(0)$.
	This means that~$\subconewo_4$ is a subcone of $\mathcal{F}_4$ which eventually contains the functionals~$c_{T_j}$ as above, with~$r_j$ non-divisible by~$2$.
	
	The argument above generalizes to any $m$ non-squarefree by Corollary~\ref{cor;ratwithsiegsernonprime}.
	This leads us to the following conjecture.
	In Proposition \ref{lemma;conjltokifltdif} we check that this conjecture implies the existence of $\subconem(\vect{\lambda})$, for every $\vect{\lambda}\in\limitmm{m}{k}$.
	
	We recall that if $T$ is a matrix in $\halfint^+_2$ and $t$ is a positive integer, we denote by $\Tdiv{T}{t}$ the matrix arising as in \eqref{not;matTdiv}.
	\begin{conj}\label{conj;forQmnonsqfr}
	Let $k>4$, $k\equiv 2$ mod $4$.
	For every reduced matrix $T$ in $\halfint^+_2$ with bottom-right entry $m$, there exist an arbitrarily large integer $A$ and positive rational numbers $\mu_{n,r}$ such that
	\be\label{eq;forQmbutmnonsqfr}
	\sum_{1\le n\le A}\sum_{r\in S_m(n,T)}\mu_{n,r}c_{n,r}|_{J_{k,m}^{\text{\rm cusp}}(\QQ)}=0,
	\ee
	where the auxiliary set $S_m(n,T)$ is defined as
	\bes
	S_m(n,T)=\Big\{r\in\ZZ\,:\,\text{$\widetilde{T}=\smT\in\halfint^+_2$ and $\frac{a^k_2(\Tdiv{\widetilde{T}}{t_j})}{a^k_2(\widetilde{T})}=\frac{a^k_2(\Tdiv{T}{t_j})}{a^k_2(T)}$ for $j=1,\dots,d$}\Big\}.
	\ees
	\end{conj}
	We remark that Conjecture \ref{conj;forQmnonsqfr} is analogous to Lemma~\ref{lem:relcoefextjac}, but with~\eqref{eq:relcoefextjac} restricted to the \textcolor{\myblack}{indices}~$r$ lying in~$S_m(n,T)$.
	If~$m$ is squarefree, this conjecture coincides exactly with Lemma~\ref{lem:relcoefextjac}.
	If~$m$ is non-squarefree, and $T=\big(\begin{smallmatrix}
	n & r/2\\
	r/2 & m
	\end{smallmatrix}\big)$ is such that~$r$ is not divisible by any divisor $t$ of $m$ with $t^2|m$, by Corollary \ref{cor;ratwithsiegsernonprime} the auxiliary set $S_m(n,T)$ simplifies to
	\bes
	S_m(n,T)=\{r\in\ZZ\,:\,\text{$4nm-r^2>0$ and if $t^2|m$ with $t\neq1$, then $t\nmid r$}\}.
	\ees
	\begin{prop}\label{lemma;conjltokifltdif}
	Let $m$ be a non-squarefree positive integer.
	If Conjecture~\ref{conj;forQmnonsqfr} holds, then there exists a subcone~$\subconem(\vect{\lambda})$ as hypothesized in Theorem~\ref{thm;polyhedwithhyp}, for every~$\vect{\lambda}\in\limitmm{k}{m}$.
	\end{prop}
	\begin{proof}
	Let~$\vect{\lambda}=(\lambda_{t_1},\dots,\lambda_{t_d})\in\limitmm{k}{m}$, and let~$t_0=1$.
	Suppose for simplicity that the convex hull generated by the $d+1$ points~$V_m,V_{m/t_1^2},\dots, V_{m/t_d^2}$ is~$d$-dimensional.
	This is equivalent to
	\bes
	\mathcal{H}\coloneqq\conv_{\RR}(\{ V_{m/t_j^2}\,:\, j=0,\dots,d\})
	\ees
	being a simplex, with the points~$V_{m/t_j^2}$ as vertexes.
	In a simplex, each point can be written as a convex combination of the vertexes in a unique way.
	Hence, if we choose two different tuples~$\vect{\lambda}$ and~$\vect{\lambda}'$ in~$\limitmm{m}{k}$, also the associated points~$\Qmsqf(\vect{\lambda})$ and~$\Qmsqf(\vect{\lambda}')$ are different; see Theorem~\ref{thm;convhull} for further information.
	
	The following two cases prove the result under the assumption that $\mathcal{H}$ is a simplex.
	Eventually, we illustrate how to generalize the proof to the case where $\mathcal{H}$ is not a simplex.
	
	\textbf{First case.} Suppose that $\vect{\lambda}$ is a \emph{non-special tuple of limits}, that is, there exists a reduced matrix $T$ in $\halfint^+_2$ such that $a^k_2(\Tdiv{T}{t_j})/a^k_2(T)=\lambda_{t_j}$ for every $j=1,\dots,d$\textcolor{\myblack}{; recall special limits from Definition~\ref{def;setspeclim}}.
	We prove the existence of~$\subconem(\vect{\lambda})$ under this hypothesis.	
	The idea is analogous to the one used to prove Lemma~\ref{lemma;msqfrsubcone}.
	Let $A$ and $\mu_{n,r}$ be as in Conjecture~\ref{conj;forQmnonsqfr}. For simplicity, we denote by~$T_{n,r}$ the matrix~$\big(\begin{smallmatrix}
	n & r/2\\ r/2 & m
	\end{smallmatrix}\big)$, with~$m$ fixed. Writing the functionals~$c_{T_{n,r}}$ over the basis~\eqref{eq;basisMk2}, we deduce that
	\bas
	\sum_{1\le n\le A}&\sum_{r\in S_m(n,T)}\mu_{n,r}c_{T_{n,r}}
	=\\
	=&\sum_n\sum_r \mu_{n,r} a^k_2(T_{n,r})
	\underbrace{\left(\begin{smallmatrix}
	1 \\
	\const\cdot \sum_{j=0}^d\lambda_{t_j}\alpha_m(t_j,f_1) \\
	\vdots \\
	\const\cdot \sum_{j=0}^d\lambda_{t_j}\alpha_m(t_j,f_\ell) \\
	0 \\
	\vdots \\ 
	0
	\end{smallmatrix}\right)}_{=\Qmsqf(\vect{\lambda})}
	+
	\underbrace{\left(\begin{smallmatrix}
	0 \\  \sum_n\sum_r \mu_{n,r} c^{f_1}_m(n,r) \\ \vdots \\ \sum_n\sum_r \mu_{n,r} c^{f_\ell}_m(n,r) \\ \sum_n\sum_r \mu_{n,r} c_{T_{n,r}}(F_1) \\ \vdots \\ \sum_n\sum_r \mu_{n,r} c_{T_{n,r}}(F_{\ell'})
	\end{smallmatrix}\right)}_{=0}=\\
	=&\Qmsqf(\vect{\lambda})\cdot\sum_{1\le n\le A}\sum_{r\in S_m(n,T)} \mu_{n,r} a^k_2(T_{n,r}).
	\eas
	Define
	\bes
	\xi(A)\coloneqq 1\big/\sum_{1\le n\le A}\sum_{r\in S_m(n,T)} \mu_{n,r} a^k_2(T_{n,r}).
	\ees
	We may rewrite~$\Qmsqf(\vect{\lambda})$ as a linear combination of functionals as
	\be\label{eq;conjQmnonsqfree}
	\Qmsqf(\vect{\lambda})=\sum_{1\le n\le A}\sum_{r\in S_m(n,T)}\xi(A)\mu_{n,r}c_{T_{n,r}}.
	\ee
	\textcolor{\myblack}{Note that the linear coefficients above are \emph{positive}, since so are~$\mu_{n,r}$ and~$a^k_2(T_{n,r})$, respectively by hypothesis and by Lemma~\ref{classsiegeis}~\ref{lem;pt1class}.}

	We prove that we can choose the subcone $\subconem(\vect{\lambda})$ as
	\be\label{eq;proofdefRmlam}
	\subconem(\vect{\lambda})=\left\langle c_{\widetilde{T}}\,:\,\text{$\widetilde{T}=\smT\in\halfint^+_2$ and $r\in S_m(n,T)$}\right\rangle_{\RR_{\ge0}}.
	\ee
	Since the value $A$ in~\eqref{eq;conjQmnonsqfree} can be chosen arbitrarily large, we may suppose that the coefficient extraction functionals appearing in~\eqref{eq;conjQmnonsqfree} generate a vector space over~$\QQ$ with the same dimension as~$\subconem(\vect{\lambda})$.
	In this way, the point~$\Qmsqf(\vect{\lambda})$ is internal in~$\subconem(\vect{\lambda})$.
	Let~$\seq{\RR_{\ge0}\cdot c_{T_j}}{j\in \NN}$ be a sequence of rays which converges to~${\RR_{\ge0}\cdot \Qmsqf(\vect{\lambda})}$, where the reduced matrices~$T_j$ in~$\halfint^+_2$ are of increasing determinant.
	Since~$\vect{\lambda}$ is a non-special tuple of limits in~$\limitmm{m}{k}$ by hypothesis, by Corollary~\ref{cor;ratwithsiegser} and Corollary~\ref{cor;ratwithsiegsernonprime} the functionals~$c_{T_j}$ eventually lie in~$\subconem(\vect{\lambda})$.
	
	\textbf{Second case.}
	We prove the existence of the subcone~$\mathcal{R}_m(\vect{\lambda}')$ associated to a \emph{special tuple of limits}~$\vect{\lambda}'=(\lambda'_{t_1},\dots,\lambda'_{t_d})\in\limitmmsp{m}{k}$, as hypothesized in Theorem~\ref{thm;polyhedwithhyp}.
	
	Since the~$d$-dimensional convex hull~$\mathcal{H}$ is a simplex, if~$\lambda_{t_j}'=0$ for some~$j$, then~$\Qmsqf(\vect{\lambda}')$ lies on the boundary of~$\mathcal{H}$.
	In fact, the point~$\Qmsqf(\vect{\lambda}')$ lies in the interior of the convex hull~${\conv_{\RR}(V_{m/t_j^2}\,:\,\lambda_{t_j}'\neq0)}$; see Lemma~\ref{lemma;decQmlambdas} and Corollary~\ref{cor;ratwithsiegsernonprime}.
	
	Let $\mathcal{F}_m(\vect{\lambda}')$ be the subcone of the~$\RR$-closure $\overline{\conemod}$ defined as
	\bes
	\mathcal{F}_m(\vect{\lambda}')=\Big\langle c_T\,:\,\substack{\text{$T\in\halfint^+_2$ with $m$ as bottom-right entry and such that}\\
	\text{if $\lambda_{t_j}'=0$ for some $j$, then $a^k_2(\Tdiv{T}{t_j})/a^k_2(T)=0$}}\Big\rangle_{\RR_{\ge0}}.
	\ees
	Choose~$\subconem(\vect{\lambda}')$ to be the cone generated by~$\mathcal{F}_m(\vect{\lambda}')$ and all~$V_{m/t_j^2}$ such that~$\lambda_{t_j}'\neq 0$.
	We prove that this subcone of~$\cone_k$ fulfills the properties hypothesized in Theorem~\ref{thm;polyhedwithhyp}.
	Let~$\seq{c_{T_i}}{i\in\NN}$ be a sequence of coefficient extraction functionals associated to reduced matrices~$T_i=\big(\begin{smallmatrix}
	n_i & r_i/2\\
	r_i/2 & m
	\end{smallmatrix}\big)$ in~$\halfint^+_2$ of increasing determinant.
	If~${\RR_{\ge0}\cdot c_{T_i}\to\RR_{\ge0}\cdot \Qmsqf(\vect{\lambda}')}$, then the entries~$r_i$ must be eventually non-divisible by any~$t_j$ such that~$\lambda_{t_j}'=0$ by Corollary~\ref{cor;ratwithsiegsernonprime}.
	This implies that the functionals~$c_{T_i}$ must eventually lie in~$\mathcal{F}_m(\vect{\lambda}')$.
	Therefore, it is enough to prove that~$Q_m(\vect{\lambda}')$ is internal in $\subconem(\vect{\lambda}')$.
	By Corollary~\ref{cor;infmanytup} we deduce that
	\be\label{eq;prooflastpFm}
	\mathcal{F}_m(\vect{\lambda}')=
	\Big\langle \subconem(\vect{\lambda})\,:\,\substack{\text{$\vect{\lambda}=(\lambda_{t_1},\dots,\lambda_{t_d})\in\limitmm{m}{k}\setminus\limitmmsp{m}{k}$ and}\\
	\text{such that if $\lambda'_{t_j}=0$, then $\lambda_{t_j}=0$}}\Big\rangle_{\RR_{\ge0}},
	\ee
	with $\subconem(\vect{\lambda})$ defined as in~\eqref{eq;proofdefRmlam}.	
	We recall that the latter cone contains~$\Qmsqf(\vect{\lambda})$ as internal point, by the \emph{first case} of this proof.
	Since~$\Qmsqf(\vect{\lambda}')$ is internal in~$\conv_{\RR}(V_{m/t_j^2}\,:\,\lambda_{t_j}'\neq0)$, we may choose a point~$W_{\vect{\lambda}}\in\conv_{\RR}(V_{m/t_j^2}\,:\,\lambda_{t_j}'\neq0)$ for every~$\vect{\lambda}\in\limitmm{m}{k}\setminus\limitmmsp{m}{k}$ appearing in~\eqref{eq;prooflastpFm}, such that~$\Qmsqf(\vect{\lambda}')$ is internal in the segment joining~$\Qmsqf(\vect{\lambda})$ with~$W_{\vect{\lambda}}$.
	In fact, also such~$\Qmsqf(\vect{\lambda})$ is contained in~$\conv_{\RR}(V_{m/t_j^2}\,:\,\lambda_{t_j}'\neq0)$.
	In this way, the subcone of~$\subconem(\vect{\lambda}')$ generated by~$\subconem(\vect{\lambda})$ and~$W_{\vect{\lambda}}$ contains~$\Qmsqf(\vect{\lambda}')$ as \textcolor{\myblack}{an} internal point, for every~$\vect{\lambda}$ as above.
	
	Since the cone generated by the union of cones containing~$\Qmsqf(\vect{\lambda}')$ as internal point is a cone that contains~$Q_m(\boldsymbol{\lambda}')$ in its interior, we may deduce that~$Q_m(\vect{\lambda}')$ is an internal point of~$\subconem(\vect{\lambda}')$.
	
	\textbf{$\boldsymbol{\mathcal{H}}$ not a simplex.}
	If~$\mathcal{H}$ is \emph{not} a simplex, we might have that~$\Qmsqf(\vect{\lambda})=\Qmsqf(\vect{\mu})$ for some \emph{different} tuples of limits~$\vect{\lambda}=\vect{\mu}\in\limitmm{m}{k}$.
	When this happens, in each of the previous cases one can substitute the subcone~$\subconem(\vect{\lambda})$ with the \emph{union} of all~$\mathcal{R}_m(\vect{\mu})$ constructed therein such that~$\Qmsqf(\vect{\lambda})=\Qmsqf(\vect{\mu})$.
	\end{proof}
	
	\section{The accumulation rays of the cone of special cycles}\label{sec;acrayschow}
	
	In the previous sections we classified all possible accumulation rays of the modular cone~$\conemod$, writing the generators of these rays as linear combinations of certain points in~$\overline{\mathcal{S}_k}$.
	In this section, we use the classification above to deduce the accumulation rays of the cone of special cycles~$\conecy$ in~$\CH^2(X_\Gamma)\otimes\QQ$.
	In particular, we show that these rays are generated by linear combinations of (rational classes of) special cycles associated to \emph{singular} matrices in~$\halfint_2$.
	These special cycles are the intersection between Heegner divisors and the dual class~$\{\omega^*\}$ of the Hodge bundle.
	Eventually, we rewrite the accumulation rays of~$\conecy$ in term of \emph{primitive Heegner divisors}, which may be considered (up to a factor~$2$) as the irreducible components of the classical Heegner divisors.\\
	
	Let $X_\Gamma$ be an orthogonal Shimura variety associated to \textcolor{\myblack}{an} even unimodular lattice of signature~$(b,2)$, where~$b>2$, and let~$k=1+b/2$. Since the lattice is unimodular, the value~$k$ is an integer satisfying the relations~$k>4$ and $k\equiv 2$ mod $4$.
	We know that there exists a linear map
	\bes
	\psi_\Gamma\colon M^k_2(\QQ)^*\longrightarrow \CH^2(X_\Gamma)\otimes\QQ,\quad c_T\longmapsto \{Z(T)\}\cdot\{\omega^*\}^{2-\rank(T)},
	\ees
	which maps every functional $c_T$ to (the rational class of) the special cycle associated to the matrix $T$; see Section \ref{sec:concoefextrfun} for details.
	For simplicity, we denote by $\psi_\Gamma$ also its extension over~$\RR$.
	As proven with Corollary~\ref{cor;prespoint&acc}, the accumulation rays of~$\conecy$ are images \textcolor{\myblack}{under}~$\psi_\Gamma$ of accumulation rays of~$\conemod$.
	
	As usual, we consider every functional in~$M^k_2(\QQ)^*$ as a vector in $\QQ^{\dim M^k_2}$, writing it over a fixed basis of the form
	\bes
	E^k_2,E^k_{2,1}(f_1),\dots,E^k_{2,1}(f_\ell),F_1,\dots,F_{\ell'};
	\ees
	see the beginning of Section \ref{sec;accrays} for further information.
	
	We want to rewrite the images \textcolor{\myblack}{under} $\psi_\Gamma$ of the points $\Pinf$, $V_s$ and $Q_s(\lambda_{t_1},\dots,\lambda_{t_d})$, defined in Section \ref{sec;accrays} and Section \ref{sec;accumulcone}, as linear combinations of special cycles in $\CH^2(X_\Gamma)\otimes\QQ$.
	In fact, if a ray of the~$\RR$-closure of~$\cone_{X_\Gamma}$ is an accumulation ray of~$\cone_{X_\Gamma}$, then it is generated by~$\psi_\Gamma(\Pinf)$, $\psi_\Gamma(V_s)$ or~$\psi_\Gamma(\Qs)$, for some positive integer~$s$ and some tuple of limits~$(\lambda_{t_1},\dots,\lambda_{t_d})\in\limitmm{m}{k}$, as in Definition~\ref{defi;Qmlambda}.
	We want to make these images \textcolor{\myblack}{under}~$\psi_\Gamma$ explicit.
	Since we already know by Proposition~\ref{lemma;decQmlambdas} how to rewrite every~$\Qs$ as a linear combination of points~$V_{s'}$, for some~$s'>0$, we may restrict our attention only to~$\psi_\Gamma(\Pinf)$ and~$\psi_\Gamma(V_s)$.
	Recall that we denote by~$\{H_s\}$ the divisor class of the~$s$-th Heegner divisor; see Remark~\ref{rem;Heegnerspeccy}.
	\begin{prop}\label{prop;fromCktoCX}
	For every positive integer $s$, the image of the point $V_s$ \textcolor{\myblack}{under} $\psi_\Gamma$ is
	\be\label{eq;propfromCktoCX}
	\psi_\Gamma(V_s)=\frac{\zeta(1-k)}{2g_k(s)}\textcolor{\myblack}{\sum_{t^2|s}\mu(t)\{H_{s/t^2}\}\cdot\{\omega^*\}.}
	\ee
	The image of the point $\Pinf$ \textcolor{\myblack}{under} $\psi_\Gamma$ is~$\{\omega\}^2$.
	\end{prop}
	\begin{cor}\label{cor;fromCktoCX}
	For every positive integer $s$, the image of the ray $\RR_{\ge0}\cdot V_s$ \textcolor{\myblack}{under} $\psi_\Gamma$ is
	\bes
	\RR_{\ge0}\cdot\psi_\Gamma(V_s)=\RR_{\ge0}\cdot\bigg(\textcolor{\myblack}{\sum_{t^2|s}\mu(t)
	\{H_{s/t^2}\}\cdot\{\omega\}}
	\bigg).
	\ees
	The image of the ray $\RR_{\ge0}\cdot \Pinf$ \textcolor{\myblack}{under} $\psi_\Gamma$ is~$\RR_{\ge0}\cdot \{\omega\}^2$.
	\end{cor}
	\begin{proof}[Proof of Corollary \ref{cor;fromCktoCX}]
	Since $k\equiv2$ mod $4$, the value $\zeta(1-k)/2g_k(s)$ is negative for every positive integer $s$. Moreover $\{\omega^*\}=-\{\omega\}$ in $\CH^1(X_\Gamma)=\Pic(X_\Gamma)$.
	The claim follows directly from Proposition \ref{prop;fromCktoCX}.
	\end{proof}
	\begin{proof}[Proof of Proposition \ref{prop;fromCktoCX}]
	First of all, we deduce the image \textcolor{\myblack}{under}~$\psi_\Gamma$ of the point~$P_s$ defined in Section~\ref{sec;accumulcone}. Since $c_s(E^k_1)=2\sigma_{k-1}(s)/\zeta(1-k)$ for every positive integer $s$, see the Fourier expansion \eqref{eq;normeisserelliptic}, by Remark \ref{rem;coefklingcuspeq} we deduce that
	\ba\label{eq;ctrivisPsupfact}
	c_{\left(\begin{smallmatrix}
	s & 0\\ 0 & 0
	\end{smallmatrix}\right)}=&
	\left(a^k_2\left(\begin{smallmatrix}
	s & 0\\ 0 & 0
	\end{smallmatrix}\right),a^k_2\left(f_1,\left(\begin{smallmatrix}
	s & 0\\ 0 & 0
	\end{smallmatrix}\right)\right),\dots, a^k_2\left(f_\ell,\left(\begin{smallmatrix}
	s & 0\\ 0 & 0
	\end{smallmatrix}\right)\right), 0,\dots,0\right)=\\
	=&\left(c_s(E^k_1),c_s(f_1),\dots,c_s(f_\ell),0\dots,0\right)=c_s(E^k_1)\cdot P_s.
	\ea
	This implies that
	\be\label{eq;prooffromCktoCX}
	\psi_\Gamma(P_s)=\frac{\zeta(1-k)}{2\sigma_{k-1}(s)}\left\{Z\left(\begin{smallmatrix}
	s & 0\\ 0 & 0
	\end{smallmatrix}\right)\right\}\cdot\{\omega^*\}.
	\ee
	As we recalled with Remark \ref{rem;Heegnerspeccy}, the rational class~$\left\{Z\left(\begin{smallmatrix}
	s & 0\\ 0 & 0
	\end{smallmatrix}\right)\right\}$ is the Heegner divisor~$\{H_s\}$.
	
	We explained in Proposition \ref{Pmnonsqfr} how to rewrite every point $P_s$ as a linear combination of certain $V_{s'}$ for some positive $s'$.
	The idea is to reverse that formula, writing $V_s$ as a linear combination of certain $P_{s'}$.
	This can be done simply rewriting~$\alpha_s(1,f)$ as
	\bes
	\alpha_s(1,f)=\frac{1}{g_k(s)}\sum_{t^2|s}\mu(t)\sigma_{k-1}(s/t^2)\frac{c_{s/t^2}(f)}{\sigma_{k-1}(s/t^2)},
	\ees
	for every $f\in S^k_1$, from which we deduce that
	\bas
	V_s=&\left(1,\frac{\zeta(1-k)}{2}\alpha_s(1,f_1),\dots,\frac{\zeta(1-k)}{2}\alpha_s(1,f_\ell),0,\dots,0\right)
	=\\
	=&\frac{1}{g_k(s)}\sum_{t^2|s}\mu(t)\sigma_{k-1}(s/t^2)P_{s/t^2}.
	\eas
	This, together with~\eqref{eq;prooffromCktoCX}, implies \textcolor{\myblack}{that}
	\bas
	\textcolor{\myblack}{\psi(V_s)=\frac{1}{g_k(s)}\sum_{t^2|s}\mu(t)\sigma_{k-1}(s/t^2)\psi(P_{s/t^2})=\frac{\zeta(1-k)}{2 g_k(s)}\sum_{t^2|s}\mu(t)\{H_{s/t^2}\}\cdot\{\omega^*\}.}
	\eas
	
	Since the Siegel Eisenstein series $E^k_2$ is normalized, its Fourier coefficient associated to the zero-matrix is~$1$. Moreover, the Fourier coefficient associated to the zero-matrix of any other element of the chosen basis of~$M^k_2$ is trivial;
	see Remark \ref{rem;coefklingcuspeq} for the cases of Klingen Eisenstein series.
	This implies that
	\bes
	\Pinf=c_{\left(\begin{smallmatrix}
	0 & 0\\ 0 & 0
	\end{smallmatrix}\right)},
	\ees
	from which we deduce that $\psi_{\Gamma}(\Pinf)=\{\omega\}^2$.
	\end{proof}
	
	The Heegner divisors are in general reducible.
	If~$\Gamma=\bigO^+(L)$, it is possible to write every Heegner divisor as a sum of its irreducible components via the so-called \emph{primitive} Heegner divisors.
	The remaining part of this section aims to rewrite the generators of the rays~$\RR_{\ge0}\cdot\psi_\Gamma(V_s)$ given by Corollary~\ref{cor;fromCktoCX} in terms of primitive Heegner divisors.
	Eventually, we deduce that the accumulation cone of~$\cone_{X_\Gamma}$ \textcolor{\myblack}{is the cone} generated by the intersections \textcolor{\myblack}{of the} primitive Heegner divisors \textcolor{\myblack}{with} the Hodge class~$\{\omega\}$.\\
%
	
	From now on, we consider only orthogonal Shimura varieties $X_\Gamma$ arising from $\Gamma=\bigO^+(L)$, where $L$ is a even unimodular lattice of signature~$(b,2)$, with~$b>2$.
	
	We avoid to propose here a formal definition of the \emph{primitive Heegner divisors}~$\{\primHeegner{s}\}$ in~$\Pic(X_{\Gamma})$, since the construction, using primitive lattice vectors of~$L$, is similar to the one of the special cycles given in Section~\ref{sec;conesspeccycles}. We refer instead to~\cite[Section~4]{brmo} for details.

	Since~$L$ is unimodular and~$\Gamma=\bigO^+(L)$, \cite[Lemma~4.3]{brmo} implies that every primitive Heegner divisor is twice an irreducible orthogonal Shimura subvariety of~$X_\Gamma$.
	\textcolor{\myblack}{Moreover}~$\{H_s\}$ decomposes in primitive Heegner divisors as
	\be\label{eq;decprimHeegner}
	\{H_s\}=\sum_{t^2|s}\{\primHeegner{s/t^2}\},\qquad\text{for every positive integer $s$},
	\ee
	\textcolor{\myblack}{and conversely the primitive Heegner divisor~$\{\primHeegner{s/t^2}\}$ may be rewritten as}
	\be\label{eq;comlinprimheegner}
	\textcolor{\myblack}{\{\primHeegner{s}\}=\sum_{t^2|s}\mu(t)\{H_{s/t^2}\}}.
	\ee
	\textcolor{\myblack}{These are respectively}~\cite[Section~$4$, $(17)$]{brmo} \textcolor{\myblack}{and \cite[Lemma~$4.2$]{brmo}}.
	\begin{cor}\label{cor;accprimheegner}
	For every positive integer~$s$, the image of the ray~$\RR_{\ge0}\cdot V_s$ \textcolor{\myblack}{under}~$\psi_\Gamma$ is generated by \textcolor{\myblack}{a \emph{primitive} Heegner divisor} intersected with the Hodge class~$\omega$.
	More precisely
	\bes
	\psi_\Gamma(\RR_{\ge0}\cdot V_s)=\RR_{\ge0}\cdot\Big(\textcolor{\myblack}{
	\{\primHeegner{s}\}\cdot\{\omega\}}
	\Big).
	\ees
	\end{cor}
	\begin{proof}
	\textcolor{\myblack}{It is enough to rewrite} the generator of~$\psi_\Gamma(\RR_{\ge0}\cdot V_s)$ given by Corollary~\ref{cor;fromCktoCX} \textcolor{\myblack}{using~\eqref{eq;comlinprimheegner}.}
%
	\end{proof}
	
	\section{Further generalizations}\label{sec;furthgen}
	
	In this section we explain how to use the same pattern of this chapter
	to investigate the geometric properties of the cones of special cycles of higher codimension, via vector valued Siegel modular forms.\\
	
	Let $X$ be an orthogonal Shimura variety associated to a \emph{unimodular} lattice $L$ of signature~$(b,2)$, as in Section \ref{sec;conesspeccycles}.
	Some of the ideas of this chapter
	extend to the cases of special cycles of codimension $g\ge3$, as follows.
	With an analogous argument as in Section~\ref{sec:concoefextrfun}, the rational polyhedrality of the cones in $\CH^g(X)\otimes\QQ$ generated by the special cycles~$\{Z(T)\}$, associated to symmetric half-integral positive semi-definite~${g\times g}$ matrices~$T$ of fixed rank, is implied by the analogous statement on cones of functionals $c_T$ of genus $g$ Siegel modular forms with weight $k=1+b/2$.
	
	Let $M^k_g(\QQ)$ (resp.\ $S^k_g(\QQ)$) be the space of Siegel modular forms (resp.\ Siegel cusp forms) of genus $g$ and weight $k$. It is well-known that $M^k_g(\QQ)$ splits in a direct sum between~$S^k_g(\QQ)$, the space generated by the Siegel Eisenstein series $E^k_g$ of genus $g$, and the spaces of Klingen Eisenstein series associated to Siegel cusp forms of lower genus; see \cite[p.\ 73, Theorem~2]{kli;intro}.
	
	To the best of our knowledge, a clear growth of the coefficients of Klingen Eisenstein series, as in \cite{boda;petnorm} for genus~$2$, is not available in literature.
	For this reason, a generalization of this chapter
	in genus $g\ge3$ seems not yet possible.
	
	Nevertheless, we remark that the cones generated by the functionals~$c_T$ in~$M^k_g(\QQ)^*$ with~${\rank T=1}$ (or~$\rank T=2$) can be deduced \textcolor{\myblack}{with} the results of this notes, following the same idea of the proof of Theorem~\ref{thm;mainthm2}~\ref{pt;conek'}.\\
	
	Another interesting problem is to deduce the geometric properties of the cones of special cycles on orthogonal Shimura varieties~$X$ associated to lattices which are non-unimodular.
	
	Since Kudla's modularity conjecture is proved in \cite{brra;modconj} in full generality, Proposition~\ref{prop;functpsi} may be generalized for coefficient extraction functionals associated to vector valued Siegel modular forms of genus $g$, with respect to the so-called \emph{Weil representation}~$\rho_{L,g}$;
	see~\cite[Section~1.1]{br;borchp} and \cite[Section 2.1]{zh;phd} for the definition of~$\rho_{L,g}$.
	
	The main obstacle to this approach are the properties of the Fourier coefficients of such vector valued modular forms.
	In fact, to the best of out knowledge, not only the growths of the coefficients of the Siegel Eisenstein series and the Klingen Eisenstein series are not yet clear, but also an explicit ``Coefficient Formula'' to compute them is missing, even in the case of genus $2$.\\
	
	For certain non-degenerate quadratic spaces over totally real fields of finite degree, an analogous construction of orthogonal Shimura varieties (and of special cycles) holds; see e.g.~\cite[Section 1.1]{maeda}.
	Kudla's modularity conjecture has been recently proved also for these generalizations, assuming the Beilinson--Bloch conjecture; see \cite[Theorem~1.6]{maeda} and \cite[Theorem~1.1]{kudrem-ma}.
	In this setting, the generating series of special cycles is a \emph{Hilbert--Siegel} modular form, with values in the Chow ring.
	Since Proposition \ref{prop;functpsi} may generalize to these kind of modular forms, it might be interesting to study also cones of coefficient extraction functionals associated to Hilbert--Siegel modular forms.
	
	\section{Examples of convex hulls in $\cone_k$ for fixed $m$}\label{subsec;examplesmfixed}
	
	Let $k$ be an integer such that $k>4$ and $k\equiv 2$ mod $4$.
	To make Theorem~\ref{thm;convhull} and Proposition~\ref{Pmnonsqfr} as clear as possible, in this section we compute explicitly the convex hull in $\overline{\cone_k}$ generated by the points~$\Qmsqf(\vect{\lambda})$, where~$\vect{\lambda}\in\limitmm{k}{m}$, for $m=4$ and $36$.
	As usual, see Section~\ref{sec;accrays}, we represent the coefficient extraction functionals $c_T$ associated to matrices~${T\in\halfint^+_2}$ as vectors in $\QQ^{\dim M^k_2}$ over a basis of the form
	\bes
	E^k_2,E^k_{2,1}(f_1),\dots,E^k_{2,1}(f_\ell),F_1,\dots,F_{\ell'}.
	\ees
	
	\subsection{Case $\boldsymbol{m=4}$}
	Let $\seq{\RR_{\ge0}\cdot c_{T_j}}{j\in\NN}$ be a sequence of rays in~$\overline{\conemod}$ associated to reduced matrices in~$\halfint^+_2$ of the form~$T_j=\big(\begin{smallmatrix}
	n_j & r_j/2\\ r_j/2 & 4
	\end{smallmatrix}\big)$, with increasing determinant and bottom-right entry fixed to $m=4$. As we observed in Section \ref{subsec;accrays}, all accumulation rays in $\overline{\cone_k}$ arising from such sequences are of the form $\RR_{\ge0}\cdot Q_4(\lambda_2)$, where
	\bes
	Q_4(\lambda_2)=\left(\begin{smallmatrix}
	1 \\
	\const\cdot\alpha_4(1,f_1)+\const\cdot\lambda_2\alpha_4(2,f_1) \\ \vdots \\
	\const\cdot\alpha_4(1,f_\ell)+\const\cdot\lambda_2\alpha_4(2,f_\ell)
	\\ 0 \\ \vdots \\0 \end{smallmatrix}\right),\qquad\text{for some $\lambda_2\in\limitmm{4}{k}$,}
	\ees
	where we abbreviate $\zeta=\zeta(1-k)/2$. As usual, we may suppose that the sequence of ratios~$a^k_2(\Tdiv{T_j}{2})/a^k_2(T_j)$ is convergent, and we denote its limit by $\lambda_2$. With the same notation of Proposition \ref{thm;bodapetn}, we compute
	\begin{equation}\label{eq;m=4}
	\begin{split}
	& \alpha_4(1,f)=\frac{c_4(f)-c_1(f)}{\ssigma(4)-1},\\
	& \alpha_4(2,f)=c_1(f)-\alpha_4(1,f),\\
	& \alpha_4(1,f)+\lambda_2\alpha_4(2,f)=(1-\lambda_2)\alpha_4(1,f)+\lambda_2c_1(f),
	\end{split}
	\end{equation}
	for every cusp form~$f\in S^k_1$.
	
	Recall that $V_s=\Big(
	1, \const\cdot\alpha_s(1,f_1), \dots, \const\cdot\alpha_s(1,f_\ell), 0, \dots, 0\Big)^t$ for every positive integer $s$. We deduce from \eqref{eq;m=4} that the points $Q_4(\lambda_2)$ \emph{lie on the segment connecting $V_1$ with $V_4$}. This verifies Theorem \ref{thm;convhull} for $m=4$, since the segment above is the convex hull over $\RR$ generated by $V_1$ and $V_4$.
	More explicitly, these points satisfy the formula
	\bes
	Q_4(\lambda_2)=(1-\lambda_2)V_4+\lambda_2 V_1.
	\ees
	
	By Corollary \ref{cor;ratwithsiegser} and Proposition \ref{prop;infmanylp}, whenever $V_1$ and $V_4$ are \emph{different}, there are infinitely many points $Q_4(\lambda_2)$, which accumulate towards some $Q_4(\lambda'_2)$, where~$\lambda'_2$ is a special limit in $\limitmsp{4}{k}{2}$. The number of such accumulation points is finite; see Remark~\ref{rem;finnuminlimitmspp}.
	Figure~\ref{fig;segment_with_Qm} represents the general case of such arrangement of points.
	
	We recall that $P_s=\Big(1,\const\cdot c_s(f_1)/\ssigma(s),\dots,\const \cdot c_s(f_\ell)/\ssigma(s),0,\dots,0\Big)^t$, for every positive integer~$s$.
	If~$s$ is squarefree, the points~$P_s$ and~$V_s$ coincide.
	The point $P_4$ is \emph{internal} in the segment generated by $V_1$ and $V_4$.
	In fact, it is easy to see \textcolor{\myblack}{by}~\eqref{eq;m=4} that if~$\lambda=1/\ssigma(4)$, then
	\bes
	(1-\lambda)\alpha_4(1,f)+\lambda c_1(f)=\frac{c_4(f)}{\ssigma(4)},\quad \text{for every $f\in S^k_1$}.
	\ees
	This is a direct check of Proposition \ref{Pmnonsqfr} for $s=4$.

	\subsection{Case $\boldsymbol{m=36}$}
	Since this example is similar to the previous one, we omit some details.
	Consider all accumulation rays given by sequences of reduced matrices~$\seq{T_j}{j\in\NN}$ of increasing determinant and bottom-right entries fixed to~${m=36}$.
	These rays are generated by
	\bes
	Q_{36}(\lambda_2,\lambda_3,\lambda_6)=
	\left(\begin{smallmatrix}
	1 \\
	\const\cdot\alpha_{36}(1,f_1)+\const\cdot\lambda_2\alpha_{36}(2,f_1)+\const\cdot\lambda_3\alpha_{36}(3,f_1)+\const\cdot\lambda_6\alpha_{36}(6,f_1) \\ \vdots \\
	\const\cdot\alpha_{36}(1,f_\ell)+\const\cdot\lambda_2\alpha_{36}(2,f_\ell)+\const\cdot\lambda_3\alpha_{36}(3,f_\ell)+\const\cdot\lambda_6\alpha_{36}(6,f_\ell)
	\\ 0 \\ \vdots \\0 \end{smallmatrix}\right),
	\ees
	for some~$(\lambda_2,\lambda_3,\lambda_6)\in\limitmm{36}{k}$.
	Via simple computations we deduce that
	\begin{align*}
	&\alpha_{36}(1,f)=\frac{c_{36}(f)-c_9(f)-c_4(f)+c_1(f)}{\ssigma(36)-\ssigma(9)-\ssigma(4)+1},\\
	&\alpha_{36}(2,f)=\alpha_{9}(1,f)-\alpha_{36}(1,f),\\
	&\alpha_{36}(3,f)=\alpha_4(1,f)-\alpha_{36}(1,f),\\
	&\alpha_{36}(6,f)=c_1(f)-\alpha_4(1,f)-\alpha_9(1,f)+\alpha_{36}(1,f),\\
	&\alpha_{36}(1,f)+\lambda_2\alpha_{36}(2,f)+\lambda_3\alpha_{36}(3,f)+\lambda_6\alpha_{36}(6,f)=\\
	&\qquad=\Big(1-\lambda_2-\lambda_3+\lambda_6\Big)\alpha_{36}(1,f)+\Big(\lambda_2-\lambda_6\Big)\alpha_9(1,f)+\Big(\lambda_3-\lambda_6\Big)\alpha_4(1,f)+\lambda_6 c_1(f),
	\end{align*}
	for every~$f\in S^k_1$.
	Hence, it is clear that
	\bes
	Q_{36}(\lambda_2,\lambda_3,\lambda_6)=\Big(1-\lambda_2-\lambda_3+\lambda_6\Big)V_{36}+\Big(\lambda_2-\lambda_6\Big)V_9+\Big(\lambda_3-\lambda_6\Big)V_4+\lambda_6V_1.
	\ees
	
	By Lemma \ref{lem;lemratnew} and Corollary \ref{cor;ratwithsiegsernonprime}, it follows that
	\bes
	\lambda_3,\lambda_2<1\quad\text{and}\quad\lambda_6\le\lambda_3,\lambda_2.
	\ees
	The inequality $1-\lambda_2-\lambda_3+\lambda_6\ge 0$ is less trivial, but was proved in~\eqref{eq;thmconvhull} \textcolor{\myblack}{by means of multiplicative arithmetic functions}.
	This implies that~$Q_{36}(\lambda_2,\lambda_3,\lambda_6)$ is contained in the convex hull generated by~$V_1$, $V_4$, $V_9$ and~$V_{36}$.
	\begin{rem}
	Suppose that the convex hull generated by $V_1$, $V_4$, $V_9$ and $V_{36}$ is a simplex of dimension~$3$.
	Since every point in a simplex can be written as a convex sum of the vertexes of the simplex in a unique way, we deduce that for \emph{different} values~${(\lambda_2,\lambda_3,\lambda_6)\in\limitmm{m}{k}}$ we have \emph{different} points~$Q_{36}(\lambda_2,\lambda_3,\lambda_6)$.
	Under this hypothesis, Corollary~\ref{cor;infmanytup} implies that the points~$Q_{36}(\lambda_2,\lambda_3,\lambda_6)$ accumulate towards \emph{an infinite number of points} of~$\conv_{\RR}(\{V_1,V_4,V_9,V_{36}\})$.
	\end{rem}
	
	The point $P_{36}$ lie in the convex hull given by $V_1$, $V_4$, $V_9$ and~$V_{36}$.
	In fact, it is easy to check that
	\bes
	P_{36}=\Big(1-\lambda-\lambda'+\lambda''\Big)V_{36}+\Big(\lambda-\lambda''\Big)V_9+\Big(\lambda'-\lambda''\Big)V_4+\lambda''V_1.
	\ees
	with $\lambda=1/\ssigma(4)$, $\lambda'=1/\ssigma(9)$ and $\lambda''=1/\ssigma(36)$.
	This is a direct check of Proposition~\ref{Pmnonsqfr} for~$s=36$.
	
	\printbibliography
\end{document}